%% file: bare_conf.tex
\documentclass[conference]{IEEEtran}

\ifCLASSINFOpdf
\else
\fi
\usepackage{amsmath}
\usepackage{comment}
\usepackage{cite}
\usepackage{amsmath,amssymb,amsfonts}
\usepackage{algorithmic}
\usepackage{graphicx}
\usepackage{pgfplots}
\usepackage{textcomp}
\usepackage{xcolor}
\usepackage{multirow}
\usepackage{pifont}
\usepackage{url}
\usepackage{tikz}
\usetikzlibrary{patterns}

\usepackage{acronym}
\usepackage{multirow}
\usepackage[flushleft]{threeparttable} % to use notes for the table
\usepackage[T1]{fontenc} % allows the direct use of >, <, = etc.

\usepackage{subcaption}

\usepackage{enumitem}

\usepackage[utf8]{inputenc}
\usepackage{eurosym}

\usepackage{booktabs} % To thicken table lines

\usepackage{url}
\usepackage{threeparttable}
\usepackage{makecell}
\usepackage{pifont}
\usepackage{graphicx}
\usepackage{pgfplots}
\pgfplotsset{compat=newest}
\usepackage{textcomp}
\usepackage{xcolor}
\usepackage{siunitx}

\usepackage{url}

\usepackage[]{fancyhdr} % 
\newcommand{\changefont}{\fontsize{9}{9}\selectfont}
\IEEEoverridecommandlockouts
\begin{document}

\title{The Value of Flexibility in a Carbon Neutral Power System}
\author{\IEEEauthorblockN{Elena Raycheva\textsuperscript{1}\textsuperscript{2}, Jared B. Garrison\textsuperscript{3}, Christian Schaffner\textsuperscript{2}, Gabriela Hug\textsuperscript{1}}
\IEEEauthorblockA{\textit{\textsuperscript{1}EEH - Power Systems Laboratory, \textsuperscript{2}ESC - Energy Science Center}, \textsuperscript{3}FEN - Research Center for Energy Networks}
ETH Z\"urich, Z\"urich, Switzerland \\
Emails: \{raycheva, hug\}@eeh.ee.ethz.ch, garrison@fen.ethz.ch, schaffner@esc.ethz.ch}

\maketitle
\thispagestyle{fancy}
\pagestyle{fancy}

%\thispagestyle{fancy}
%\pagestyle{fancy}

\begin{abstract}
In this paper, we use a formulation of the generation expansion planning problem with hourly temporal resolution, to investigate the impact of the availability of different sources of flexibility on a carbon-neutral central European power system (with a focus on Switzerland) for the year 2040. We assess the role of flexible generation, load shifting and imports on the investment and operation of existing and newly built units. Our results show that including load shifting as part of the optimization could reduce the need for investments in both RES and conventional technologies. The combination of newly built flexible generators (gas turbines) and load shifting increases the flexibility of the simulated power system and results in the overall lowest system costs. The reduction of cross-border transmission capacity between Switzerland and its neighbors has a significant impact on the domestic operation and investments but could also affect the surrounding countries.

\end{abstract}

\begin{IEEEkeywords}
cross border congestion management, flexible generation, generation expansion planning, load shifting, RES integration
%The author shall provide up to 5 keywords (in alphabetical order) to help identify the major topics of the paper. The thesaurus of IEEE indexing keywords is posted at %\url{http://www.ieee.org/organizations/pubs/ani_prod/keywrd98.txt.}
\end{IEEEkeywords}


% \fi
\IEEEpeerreviewmaketitle

\section{Introduction}
The integration of large shares of intermittent Renewable Energy System (RES) technologies driven by net-zero Green House Gas (GHG) emissions targets and policies related to phase out of conventional thermal and nuclear generation provides many opportunities, but also challenges. One of the main challenges is related to maintaining reliable operation of the future electric power system, which will be characterized by an increased need for flexibility to balance the fluctuations of weather-dependent generation. The main objective of this work is to integrate different sources of flexibility in a Generation Expansion Planning (GEP) formulation and study their impact on the investment and operation of existing and newly built generation and storage capacities in a carbon neutral power system. The sources of flexibility that we focus on are 1) dispatchable generation (conventional, hydro and battery storage technologies), 2) load shifting and 3) imports/exports. 

In~\cite{Palmintier2016, Ma2013} the authors investigate the operational aspect of flexibility in expansion planning, focusing on detailed operation of conventional generators using Unit Commitment (UC) constraints in the formulation of the optimization problem. These studies, however, do not include other technologies which could provide flexibility such as hydro power and Battery Energy Storage Systems (BESS). Furthermore, they do not assess the role of flexible loads on the investments and operation of the units. On the other hand, in~\cite{Maranon-Ledesma2019,Gils2016} the authors include a detailed Demand Response (DR) representation within the capacity expansion planning but omit UC constraints and the latter does not model ramp limits. While modeling very large test systems (\cite{Maranon-Ledesma2019} has a European scope and~\cite{Gils2016} models Germany), both~\cite{Maranon-Ledesma2019} and~\cite{Gils2016} do so by aggregating generators. In addition, the transmission lines are also aggregated and power flows between different zones are constrained using a transport model. In~\cite{Gils2016}, the authors treat the aggregated German power system as an island and do not model the interconnections with the surrounding countries. 

To bridge the gap between the aforementioned studies, the present work proposes a GEP formulation which includes a variety of technologies providing flexibility such as conventional thermal units, for which we include ramp limits, hydro power as well as BESS. In addition, flexible loads are modeled as part of the optimization. The GEP formulation is static (i.e. spans a single target year) and deterministic, however, we model the target year with hourly resolution in one shot to capture short-term fluctuations of intermittent RES but also the seasonal operation of hydro storages. The expansion planning includes DC power flow constraints and is conducted for a real-size power system (i.e. the Swiss power system), considering also the cross-border connections with the neighbouring countries as well as an aggregated representation of their projected generation portfolio under a carbon-neutral scenario for the year 2040. This allows to capture the impact of imports and exports on the investment and operational decisions in Switzerland.

The remainder of the paper is organized as follows: the problem formulation is detailed in Section~\ref{pf}, Section~\ref{input data} presents the test system and corresponding input data used. Section~\ref{results} describes the simulated scenarios and presents the results and Section~\ref{conclusion} concludes the paper.

%This paper (Analyzing Demand Response in a Dynamic Capacity Expansion Model for the European Power Market) addresses the problem dynamic model with 168h-336h per simulated year, no dc power flow, one node per country

%Economic potential for future demand response in Germany – Modeling
%approach and case study --aggregated generators, no detail of the grid (islanded system) + no ramp limits 

%4) important to maintain chronological accuracy and high temporal resolution (for the daily load shifting constraints it's good to model the days without only using representative hours because as we saw the operation within a day can be quite different and our load shifting constraints are set over a day)

\section{Problem Formulation}\label{pf}
We use the Centralized Investments module (CentIv), a core module of the interconnected energy systems modeling platform Nexus-e~\cite{GJORGIEV2022118193} to perform the simulations. The goal of CentIv, which is based on previous work~\cite{Raycheva2020}, is to co-optimize the capacity investment and operational decisions of units at the transmission system level. This includes large-scale thermal generators (nuclear, gas, biomass, etc.), hydro power plants (run-of-river (RoR), dam and pump storages), large-scale RES (wind, PV) and grid-scale batteries. Mathematically, CentIv minimizes the sum of the annualized investment costs and the operating costs of all existing and candidate generation and storage technologies over the planning horizon $T$ (one year):
{%\small
%\vspace{-2pt}
\begin{equation}
\small
%\begin{align}
\begin{split}
\label{eq:obj fun centiv}
\min \ \underbrace{\sum_{d \in D} \alpha_d^\text{inv} C_d^\text{inv} u_d^\text{inv}}_\text{(i)} + \underbrace{\sum_{j \in J} \sum_{t \in T} (C_{j}^\text{voc} + C_{j}^\text{fuel} + C_{j}^\text{emi}) p_{j,t}}_\text{(ii)} \\
+ \underbrace{\sum_{k \in K} \sum_{t \in T} C_{k}^\text{voc} p^\text{dis}_{k,t}}_\text{(iii)} + 
\underbrace{\sum_{r \in R} \sum_{t \in T} C_{r}^\text{voc} p_{r,t}}_\text{(iv)}
+ \underbrace{\sum_{n \in N} \sum_{t \in T} C^\text{ls} ls_{n,t}}_\text{(v)}
%\end{align}
\end{split}
\end{equation}}\normalsize
\noindent where (i) are the investment costs of building candidate generators or storages $d$, (ii-iv) are the total operating costs of each thermal unit $j$, energy storage unit $k$, and renewable generator $r$ for each time step $t$ and (v) indicates the load shedding cost at transmission node $n$. To derive the operational costs, we multiply a constant variable operating cost, $C^{voc}$, by the power produced by each thermal unit $p_{j,t}$, each hydro or battery storage device $p_{k,t}^{dis}$, and each renewable generator $p_{r,t}$. For thermal generators, the fuel cost $C_{j}^{fuel}$ and $CO_2$ emissions cost $C_{j}^{emi}$ are also accounted for. The load shedding cost is $C^{ls}$\footnote{The cost of shedding any load, $C^{ls}$, is based on a value of lost load of 10'000 EUR/MWh~\cite{MORALESESPANA2022122544}.}. The investment cost (including fixed O\&M cost) for unit $d$ from the candidate units $D$ is $C_{d}^{inv}$. The investment decision variable is $u_{d}^{inv}$ and the annuity factor, which annualizes the total investment cost for the simulated year, is $\alpha_d^{inv}$. The optimization objective \eqref{eq:obj fun centiv} is subject to four sets of constraints related to: a) investments, b) operation, c) reserve provision, and d) the transmission system (i.e. DC power flow). Next, we briefly describe the changes that were made for the current formulation compared to the previous work in~\cite{Raycheva2020}. 

\subsubsection{Investment and Operation Constraints}\label{investment changes}
In this work we simulate each hour of the target year. First, this allows us to capture short-term fluctuations of intermittent RES and demand but also battery storage and pumps. Furthermore, by simulating the full year in one shot, it's possible to better represent the seasonal operation of dams. To speed up CentIv and be able to perform simulations with such high temporal resolution, all investment decisions are linearized, i.e., $u_{d}^{inv}$ is a continuous variable with decisions for investment in thermal units as given by~\eqref{cont_cand}. The formulation for candidate thermal generators (i.e., $j \in J^D$) is simplified by removing the Unit Commitment (UC) and investment constraints from the previous formulation (i.e., (2)-(9) and (19) in~\cite{Raycheva2020}). Instead, the generated power in each time step is limited by each generator's rated capacity~\eqref{eq:pmax_cand} and ramping constraints allow for the incorporation of reserve provision~\eqref{eq:ramp1_cand}-\eqref{eq:ramp2_cand}:
\par
{\small
\vspace{-0.3cm}
\begin{subequations} \label{conventional_simplified_inv}
\begin{align}
%\text{s.t.} \qquad
& 0 \le p_{j,t} \le P_{j}^{max} u_{j}^{inv}, \forall j \in J^D, \forall t \label{eq:pmax_cand}\\
& p_{j,t-1} - p_{j,t} + (r_{j,t}^{{SCR}\downarrow} + r_{j,t}^{{TCR}\downarrow}) \le R_j^{D} u_{j}^{inv}, \forall j \in J^D, \forall t \label{eq:ramp1_cand}\\
& p_{j,t} + (r_{j,t}^{{SCR}\uparrow} + r_{j,t}^{{TCR}\uparrow}) - p_{j,t-1} \le R_j^{U} u_{j}^{inv}, \forall j \in J^D, \forall t \label{eq:ramp2_cand}\\
&0 \leq u_{j}^{inv}\leq 1,\forall j \in J^{D} \label{cont_cand}%\quad \textrm{and} \quad u_{j}^{inv}=1, \forall j \notin J^{C} 
\end{align}
\end{subequations}}\normalsize
\noindent where $P_{j}^{max}$ is the maximum power output, $R_{j}^{U/D}$ are the ramp-up/ramp-down limits and $r_{j,t}^{{SCR}\uparrow\downarrow}$, $r_{j,t}^{{TCR}\uparrow\downarrow}$ are the variables for the unit's contribution towards upward/downward secondary (SCR) and tertiary (TCR) reserves. Ramping constraints are included since~\cite{Schwele2020} shows that their omission could be the most distorting simplification of the UC in GEP studies. The existing thermal units are also constrained by the equations in~\eqref{conventional_simplified_inv} but with the investment variable $u_{d}^\text{inv}$ set to $1$. The constraints for candidate battery storages are identical to (10)-(16) and (20) in~\cite{Raycheva2020}. Similar to the thermal generators, the binary investment variable for candidate storages (i.e: $k \in K^D$) is relaxed:
\par
{\small
%\vspace{-0.3cm}
\begin{equation} \label{storages}
%\text{s.t.} \qquad
0 \leq u_{k}^{inv}\leq 1,\forall k \in K^{D}
\end{equation}}\normalsize
The constraints for candidate RES technologies are identical to (21) in~\cite{Raycheva2020}. Additionally, in the case of a required RES production target\footnote{In the context of the simulated power system, we impose that such a target is to be covered by non-hydro RES. This means that we count towards the target only generation from existing or newly built PV, wind and biomass.}, $B^{RES}$, the following constraint is added:
\par
{\small
%\vspace{-0.3cm}
\begin{equation} \label{res target}
%\text{s.t.} \qquad
\sum_{r \in R^{PV, wind}} \sum_{t \in T} p_{r,t} + \sum_{j \in J^{Bio}} \sum_{t \in T} p_{j,t} \ge B^{RES}
\end{equation}}\normalsize
\subsubsection{Transmission System Constraints}\label{investment changes}
The transmission system constraints~\eqref{eq:nodal all} begin with the active power balance at each bus node $n \in N$~\eqref{nodal balance} where $P_{n,t}^{D}$ is the nodal hourly demand, $ls_{n,t}$ refers to the load shedding variable, and the remaining terms correspond to the power output of each generator and storage system. The nodal active power $p_{n,t}$ is the sum of the active power flows of all lines $l \in L$ connected to $n$ as shown by~\eqref{nodal balance 2}. The active power flow $p_{l}$ of a single line is constrained by~\eqref{line flow1} and \eqref{line flow2} where $B_{l}$ is the admittance, $\delta_n$, $\delta_i$ are the voltage angles at the start and end nodes and $P_{l}^{max}$ is the line limit. Load shedding is allowed at each bus with an associated demand and constrained by~\eqref{load shedding} which does not allow load shedding to be more than the nodal hourly load:
\par
{\small
\vspace{-0.3cm}
\begin{subequations} \label{eq:nodal all}
\begin{align}
&p_{n,t} = P_{n,t}^{D}-ls_{n,t}+\sum_{k\in K_{n,t}}p_{k,t}^{ch}-\sum_{j\in J_{n,t}} p_{j,t} \nonumber\\
&-\sum_{k\in K_{n,t}}p_{k,t}^{dis}-\sum_{r\in R_{n,t}}p_{r,t}, \quad \forall n, \forall t \label{nodal balance}\\
&p_{n,t} = \sum_{i\in l(n,i)} p_{l(n,i),t}, \forall n, \forall t \label{nodal balance 2}\\
&p_{l(n,i),t} = B_{l}(\delta_{n,t} - \delta_{i,t}), \forall l(n,i), \forall t \label{line flow1}\\
&-P_{l}^{max} \le p_{l(n,i),t} \le P_{l}^{max},\forall l(n,i), \forall t\label{line flow2}\\
&ls_{n,t} \le max(0,P_{n,t}^{D}), \forall n,\forall t,  \quad ls_ {n,t} \ge 0 \label{load shedding}
\end{align}
\end{subequations}}\normalsize

In this work, we represent the hourly nodal electricity demand profile $P_{n,t}^{D}$ as the combination of several different profiles, namely heat pump demand $P_{n,t}^{HP}$, electric mobility demand $P_{n,t}^{EM}$, hydrogen demand $P_{n,t}^{H_2}$ and the base demand $P_{n,t}^{base}$, i.e.: 
\par
{\small
\begin{equation}
%\text{s.t.} \qquad
P_{n,t}^{D} = P_{n,t}^{HP} + P_{n,t}^{EM} + P_{n,t}^{H_2} + P_{n,t}^{base}, \forall n, \forall t \label{demand profiles}
\end{equation}}\normalsize
To facilitate load shifting in CentIv we introduce four new variables. The first two of these, $e^{up/down}_{n,t}$, enable the hourly upward/downward shifting of the electric mobility load. The other two variables, $l^{up/down}_{n,t}$, allow for upward/downward shifting of the combination of all other loads (i.e., HPs, $H_2$, and the base demand). For both of these up/down load shift variables, two constraints are used to limit the maximum hourly power shift and two more to limit the maximum daily energy shift. In each hour, the upward and downward shifting is constrained as shown in~\eqref{eq:tso_emob_power}-\eqref{eq:tso_emob_power2} and~\eqref{eq:tso_dsm_power} where $E^{hr,max}$ is the maximum hourly power shift allowed for shifting the electric mobility load and $L^{hr,max}$ is the maximum hourly power shift allowed for shifting of the other loads. Similarly, the maximum energy shifted within a day is constrained by~\eqref{eq:tso_emob_energy} and~\eqref{eq:tso_dsm_energy} where $t_0$ indicates the starting time of each simulated day and $E^{day,max}$/$L^{day,max}$ are the daily energy shifting limits. Equations~\eqref{eq:tso_emob_constantdaily} and~\eqref{eq:tso_dsm_constantdaily} ensure that over each day, the up and down shifts are balanced: 
\par
{\small
\vspace{-0.3cm}
\begin{subequations} \label{eq:shifting all}
\begin{align}
& 0  \le \; \; e^{up}_{n,t}\; \le E^{hr,max}, \forall n, \forall t  \label{eq:tso_emob_power} \\
& 0  \le \; \; e^{down}_{n,t} \; \le min(E^{hr,max}, P_{n,t}^{EM}),  \forall n, \forall t \label{eq:tso_emob_power2}\\
& \sum_{t=t_0}^{t_0+24}(e^{up}_{n,t}\!+\!e^{down}_{n,t})\le
E^{day,max} \label{eq:tso_emob_energy} \\
& \sum_{t=t_0}^{t_0+24}(e^{up}_{n,t}\!-\!e^{down}_{t}) = 0\label{eq:tso_emob_constantdaily} \\
& 0  \le \; \; l^{up/down}_{n,t}\; \le L^{hr,max}, \forall n, \forall t\label{eq:tso_dsm_power} \\
%& 0  \le \; \; l^{down}_{n,t} \; \le min(L^{hr,max}, P_{n,t}^{remain}),  \forall n, \forall t \label{eq:tso_dsm_power2} \\
& \sum_{t=t_0}^{t_0+24}(l^{up}_{n,t}\!+\!l^{down}_{n,t})\le
L^{day,max} \label{eq:tso_dsm_energy} \\
& \sum_{t=t_0}^{t_0+24}(l^{up}_{n,t}\!-\!l^{down}_{t}) = 0\label{eq:tso_dsm_constantdaily}
\end{align}
\end{subequations}}\normalsize
It is important to note that the linear formulation in~\eqref{eq:shifting all} allows for simultaneous up and down shifting in each hour\footnote{To impose mutual exclusivity, it is necessary to introduce additional binary variables~\cite{MORALESESPANA2022122544}, resulting in a MILP formulation characterized by a significant increase in solver run-time.}. Once the four newly introduced and constrained variables are included in~\eqref{nodal balance}, the nodal balance equation becomes: 
\par
{\small
\vspace{-0.3cm}
\begin{equation} \label{nodal balance new}
\begin{split}
p_{n,t} &= P_{n,t}^{D}+e^{up}_{n,t}+l^{up}_{n,t}-e^{down}_{n,t}-l^{down}_{n,t}-ls_{n,t}
+\sum_{k\in K_{n,t}}p_{k,t}^{ch}\\&-\sum_{j\in J_{n,t}} p_{j,t}-\sum_{k\in K_{n,t}}p_{k,t}^{dis}-\sum_{r\in R_{n,t}}p_{r,t}, \quad \forall n, \forall t 
\end{split}
\end{equation}}\normalsize

\section{Input Data}\label{input data}

%IMPORTANT: We use all data from Calliope including the e-mobility profile !
The updated problem formulation of CentIv is applied to a representation of the central European electricity market, which is shown in Fig.~\ref{fig:2025 grid} and described in this section. This system consists of a detailed representation of Switzerland (CH) and an aggregated representation of Germany (DE), France (FR), Italy (IT) and Austria (AT) for the year 2040. Much of the input data used were taken from current or future projections of the regions modeled; however some data were also utilized from projections created by the Euro-Calliope modeling framework~\cite{TRONDLE20201929}. Section~\ref{non-ch} details the input data related to the surrounding countries, while Section~\ref{ch} describes the input data used for CH. 

\subsection{Aggregated surrounding countries (DE, FR, IT, AT)}\label{non-ch}
\begin{figure}[htbp]
  \centering
    \includegraphics[scale=0.158]{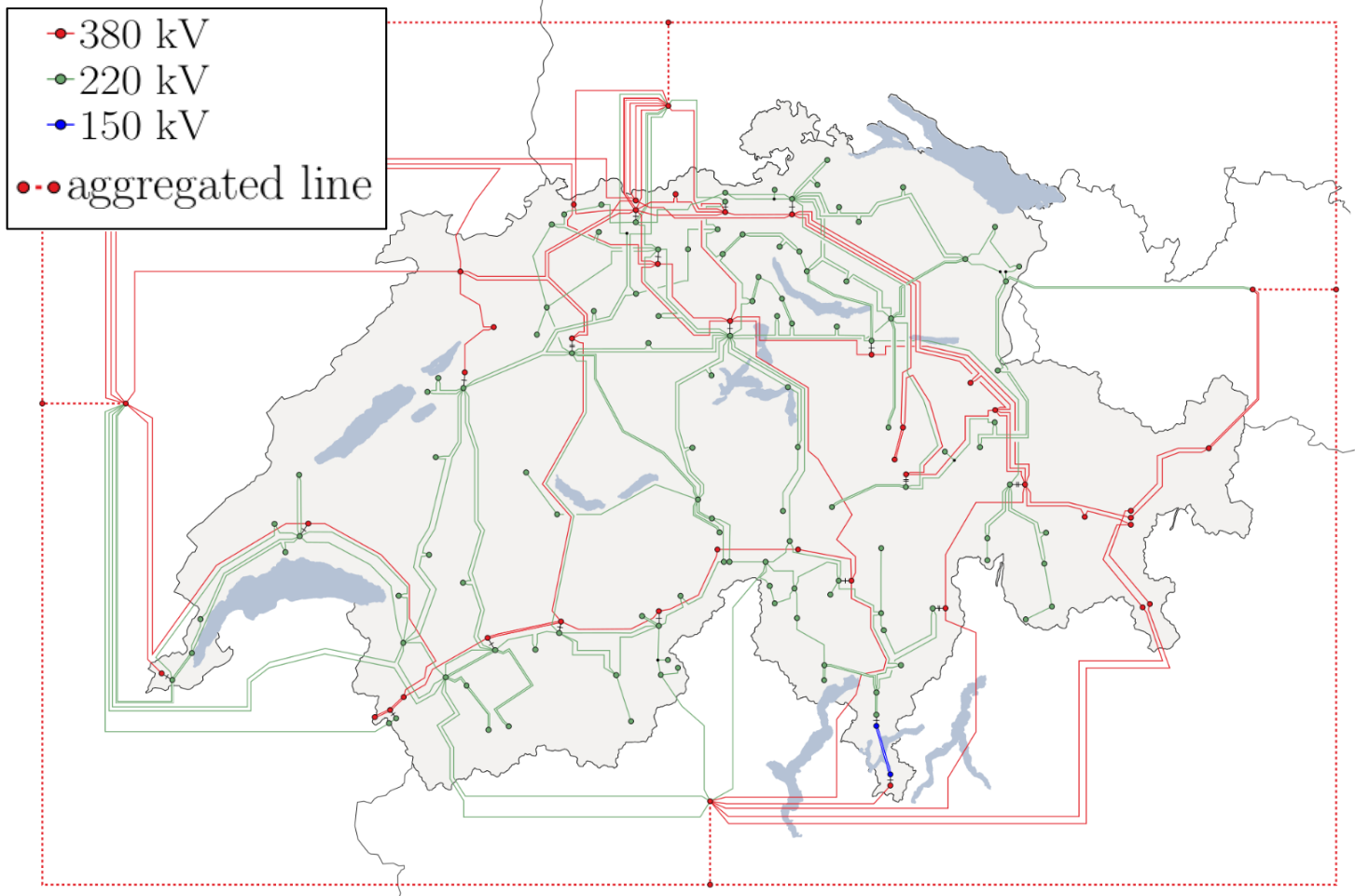}
    \vspace{-0.1cm}
    \caption[Modeled 2040 transmission network]{Overview of the modeled transmission system in 2040.}
     \label{fig:2025 grid}
\vspace{-0.2cm}
\end{figure}

%includes a detailed representation of the \gls{ch} transmission grid and an aggregated representation of the transmission grid of the four neighboring countries - \gls{de}, \gls{fr}, \gls{it}, and \gls{at}, with data from Swissgrid~\cite{Swissgrid2025} and the \gls{entsoe}~\cite{ENTSOE2018c,ENTSOE2018b}.
%The line limits of the aggregated lines between Switzerland and the neighboring countries are modified to have transfer capacities that reflect the market-based limits (i.e., \gls{ntc} or \gls{fb} limit). Analogously, the aggregated lines connecting the neighboring countries also use modified limits to reflect the market-based transfer capacities. We gathered the data for these limits on market-based transfer capacities from Swissgrid~\cite{Swissgrid2015a} and recent \gls{entsoe} studies~\cite{ENTSOE-ERAA2021,ENTSOE-MAF2020}.
The tie lines connecting DE, FR, IT and AT, along with the lines connecting these four countries to the CH border\footnote{As shown in Fig.~\ref{fig:2025 grid}, all CH cross-border lines are included but they all connect to a common border node that then connects to each neighboring country.}, are an aggregated representation of the transmission grid of the four neighboring countries. The line parameters are set by a network reduction process~\cite{Fortenbacher2018} that is applied to the physical European grid data from~\cite{ENTSOE2018c,ENTSOE2018b}. The line limits are set to reflect the market-based transfer capacities using data from~\cite{ENTSOE-ERAA2021,ENTSOE-MAF2020}. The generators and storages in DE, FR, IT and AT are aggregated to a single unit per technology type with capacity values summarized in Table~\ref{calliope capacities} and taken from Euro-Calliope model results that represent the year 2040 where the European continent has reached a nearly carbon-neutral energy system~\cite{Pickering2022}\footnote{The results in~\cite{Pickering2022} are different from the ones used in the present work due to the 2040 baseline for the rest of Europe and the higher resolution of Switzerland and the transmission constraints applied.}. The solar and wind production time series for the surrounding countries are also taken from the Euro-Calliope results~\cite{Pickering2022}. Additionally, to account for imports and exports from these four countries to their other neighbors, which are not included in this work, we utilize the cross border flow results of the Euro-Calliope analysis~\cite{Pickering2022}.
%Operational cost data for the units are from~\ref{WHRE?}. 
To represent the variable operating costs of all EU generators, we use data from~\cite{IEA2020b}. Note that, several VOM costs were adjusted as part of a calibration process performed for a separate study. This calibration process aimed to adjust generator costs and network reactances to achieve a dispatch that is more consistent with the historical 2018 conditions. Parameters that were focused on include the annual generation by generator type for all non-CH countries, the monthly generation by generator type for CH, the annual cross border flows, the annual average wholesale prices, and the monthly behavior of Swiss hydro storage units. The fuel prices for all EU generators were taken from~\cite{ENTSO-TYNDP2020}. Since all gas units in the neighboring countries are assumed to consume carbon-neutral synthetic natural gas, priced based on~\cite{Kiani2021}, no other $CO_2$ price is necessary. 
As mentioned earlier, the hourly demand profile in this work is a combination of four profiles, namely 1) e-mobility demand, 2) heat pump demand, 3) hydrogen demand and 4) base demand. The annual values for each demand type for the surrounding countries are summarized in Table~\ref{calliope non-ch demands}. The hourly profiles for each of these demands are also from Euro-Calliope~\cite{Pickering2022}.

% The grid

\begin{table}[!htbp]
	\centering
	\caption{Installed Capacities (2040) in DE, FR, IT, AT from Euro-Calliope}\label{calliope capacities}
	\renewcommand{\arraystretch}{1.2}
	\begin{tabular}{||ccccc||}%{P{1.1cm}P{2.1cm}P{1.7cm}P{1.1cm}P{0.7cm}}
		\hline
		\textbf{Gen type} & \textbf{DE} & \textbf{FR} & \textbf{IT} & \textbf{AT} \\
		\hline
		Wind [GW] & 42.6 & 172.0 & 50.4 & 52.3
		\\
		PV [GW] & 267.5 & 180.2 & 289 & < 0.1
		\\
		Biomass [GW] & 6.4 & 2.9 & 2.0 & 0.6
		\\
		\hline
		Hydro Dam [GW] & 0.2 & 9.5 & 5.5 & 5.2
		\\
		Hydro Pump [GW] & 6.4 & 5.1 & 7.9 & 3.6
		\\
		Hydro RoR [GW] & 3.5 & 6.3  & 5.9 & 4.5
		\\
		\hline
		GasCC [GW] & 0.3 & 7.0 & 14.7 & 0
		\\
		Batt [GW] & 0.1 & 0.4 & 12.3 & < 0.1
		\\
		\hline
	\end{tabular}
\end{table}

\begin{table}[!htbp]
	\centering
	\caption{Demand (2040) from Euro-Calliope}\label{calliope non-ch demands}
	\renewcommand{\arraystretch}{1.2}
	\begin{tabular}{||cccccc||}%{P{1.1cm}P{2.1cm}P{1.7cm}P{1.1cm}P{0.7cm}}
		\hline
		\textbf{Demand type} &  \textbf{DE} & \textbf{FR} & \textbf{IT} & \textbf{AT} &\textbf{CH}\\
		\hline
		E-mobility [TWh] & 339.7 & 288.3  & 199.7 & 41.6 & 30.4
		\\
		Heat Pump [TWh] & 0.2 & 57.0 & 87.7  & < 0.1 & 1.4
		\\
		Hydrogen [TWh] & 0.1 & 0.4 & 11.2 & < 0.1 & 0.3
		\\
		Base [TWh] & 457.7 & 349.8  & 290.1 & 55.9 & 49.2
		\\
		\hline
	\end{tabular}
\end{table}

%explain how the limits on shifting are determined 
The limits on the maximum daily energy shift for the e-mobility demand for each country are set using the assumption that the total annual shiftable energy can be at most 10\% of the total annual e-mobility demand. Assuming each day of the year shifts the same maximum amount allowed, this 10\% assumption translates directly into a daily maximum quantity (i.e., each day is given the same shifting allowance). Once these daily energy limits are determined, the maximum hourly power up/down shift is calculated such that the daily energy shift could be spread over 10 hours for each upward and downward shifting. The calculation of these limits is based on the e-mobility profile for each country from Euro-Calliope.
Finally, the limits on the maximum hourly power up/down shift for the other loads in each country are set equal to the DSM capacities defined for 2040 in~\cite{ENTSO-TYNDP2020}. Once these power limits are applied, the maximum daily energy shift is calculated based on the maximum up/down power shift being allowed in up to 3 hours per day each (i.e, three times the DSM capacities equals the maximum energy that can be shifted down or up within any one day).

\subsection{Detailed CH}\label{ch}
While the neighboring countries are represented with aggregated transmission network elements, the Swiss network is modeled in full detail using data from~\cite{Swissgrid2025} and shown in Fig.~\ref{fig:2025 grid}. 
For existing Swiss generator capacities and locations, we use data from~\cite{BFE2018a,BFE2019b,BFE2019a,BFE2019c} and previous studies~\cite{Abrell2019b}. The generating capacities represent those existing in 2020, of which we assume all hydro, biomass, wind, and PV to also remain in place until 2040. However, the phase out of nuclear is assumed to occur prior to 2040 along with removal of all conventional gas and oil units in Switzerland.
In addition to the existing hydro units, twenty-eight new hydro units are included that represent a mixture of pump, dam, and RoR generators with a total additional capacity of 2.6 GW. Of these additions, 1.8 GW represent three new hydro pump units that are planned to be operational between 2029 and 2037~\cite{BFE_EP2050plus}. The rest are based on projects that have been previously or are now in planning stages~\cite{BFE2019_Hydro}.
To represent the variable operating costs and fuel costs of all Swiss generators (existing and new) we use data from~\cite{Bauer2019,Bauer2017}. The costs of biomass reflect current waste incineration subsidies, which we expect to continue in the future. The VOM cost for each technology type is assumed to stay the same until 2040; however, the fuel and $CO_2$ portions of the total variable operating cost will change based on the assumed trajectories for the prices of each fuel and the price of $CO_2$ in future years.
% All: CO2 price from TYNDP2020 (high end)
Since some of the candidate gas units (i.e., those with CCS) will still emit a limited amount of $CO_2$, a price based on~\cite{ENTSO-TYNDP2020} is also included. All non-CCS gas units are assumed to consume a carbon-neutral synthetic methane (similar to what is assumed for the neighboring countries); however, the fuel price for the Swiss units of this type is set based on the average value provided for gas-to-methane from the recent Swiss study on Power-to-X~\cite{Kober2019}.
In total, in CH we model 303 existing and 82 candidate units with cost parameters summarized in Table~\ref{centiv candidate units}. 
The production time series for existing wind generators are gathered from~\cite{Abrell2019b} that included detailed assessments of the RES potentials and generation profiles. 
Since the centralized cost minimization perspective of CentIv does not capture the decisions that the end-consumers\footnote{For a TSO-DSO coordinated generation expansion planning investigation, the reader is referred to~\cite{CentIvDistIv_PowerTech2021}.} face when investing in rooftop PV (which is expected to play an increasingly important role in the future Swiss power system~\cite{GJORGIEV2022118193,Prognos}), we incorporate a non-hydro RES target of 25 TWh\footnote{This ensures that the PV generation is in the same range as the projections from~\cite{Prognos}.} for the year 2040 and introduce PV candidate units at 48 system nodes. Rooftop PV is modeled by using cantonal irradiation profiles from~\cite{Meteoswiss} and by applying a very low variable cost since most of this generation will not be seen by the transmission system but will instead be used to satisfy household demands. 

%Now discuss the Swiss existing generators and the demand ...
Similar to the neighboring countries, the hourly demand profile for Switzerland is a combination of e-mobility demand, heat pump demand, hydrogen demand and base demand. The annual totals for each demand type in CH are summarized in Table~\ref{calliope non-ch demands}. The hourly profiles for each of these demands are again from Euro-Calliope~\cite{Pickering2022}.

% Also add comments about sythetic gas versus gas-CCS candidates.

%The candidate units
\begin{table}[!htbp]
	\centering
	\caption{Cost parameters of candidate units in CH (2040) based on~\cite{Bauer2017,Bauer2019,PNNL2020} and own assumptions}\label{centiv candidate units}
	\renewcommand{\arraystretch}{1.1}
	\begin{tabular}{||cccc||}%{P{1.1cm}P{2.1cm}P{1.7cm}P{1.1cm}P{0.7cm}}
		\hline
		\textbf{Unit type} & \begin{tabular}[c]{@{}c@{}} \textbf{Inv. Cost} \\ \textbf{[kEUR/MW/a]}\end{tabular}&
		\begin{tabular}[c]{@{}c@{}} \textbf{Tot. Var. Cost} \\ \textbf{[EUR/MWh]}\end{tabular} &
		\begin{tabular}[c]{@{}c@{}} \textbf{Capacity}\\ \textbf{[MW]}\end{tabular}\\
		\hline
		Biomass & 125 & 1 & 240
		\\
		%\hline
		Wind & 151 & 36.4 &1'960
		\\
		PV & 71 & 1 & 50'000
		\\
		Battery (100MW-4h) & 204 & 0.5 & 1'100
		\\
		Gas CC (CCS) & 135 & 159.5 & 5'500
		\\
		Gas CC (Syn) & 76.5 & 308.8 & 5'500
		\\
		\hline
	\end{tabular}
	\vspace{0.1cm}
	\begin{tablenotes}
      \scriptsize
      \item Investment costs include fixed costs. The total variable cost includes fuel and $CO_2$ costs for gas. Biomass costs are subject to subsidies, while wind costs are not. For GasCC (CCS) we introduce 11 candidates: 11x500MW, for GasCC (Syn): 11x500MW, for Biomass: 12x20MW. The hourly ramp rates of all gas and biomass candidates are assumed to equal $40\%$ of the corresponding $P^{max}$ value~\cite{Garrison2014AGU}. The reserve contribution attributes per technology are summarized in Table 1 in~\cite{Raycheva2020}.
      \end{tablenotes}
\end{table}

\section{Results}\label{results}
Table~\ref{simulated scenarios} defines the four simulated scenarios for the target year 2040. The optimization problem is implemented in Pyomo~\cite{bynum2021pyomo} and solved with Gurobi~\cite{gurobi}. Section~\ref{gas and shifting} discusses the impact of load shifting and building new flexible generators on the investment and operational decisions in Switzerland. In Section~\ref{reduced grid} we show the results for the same four scenarios but impose a 70\% reduction of the cross-border limits between Switzerland and the surrounding countries. This reduction stems from the European Union's (EU) decision to reserve 70\% of the capacity of the grid elements of its member states for "internal trade". Last but not least, Section~\ref{sensitivities} includes a sensitivity analysis to show the impact of rising gas prices on investment decisions for one of the previously simulated scenarios and Section~\ref{limitations} lists the main limitations of the present work.

\begin{table}[!t]
	\centering
	\caption{Simulated Scenarios (2040)}
	\vspace{-0.2cm}
	\label{simulated scenarios}
	\renewcommand{\arraystretch}{1.1}
	\begin{tabular}{||p{3.0cm}p{1.3cm}p{2.8cm}||}
	%\begin{tabular}{||ccc||}
		\hline
		\textbf{Scenario Name}& \textbf{Abbr.}& \textbf{Description}\\
		\hline
		"No Gas No Shift" & "NGNS" & No gas candidate units + no load shifting allowed
		\\
		"No Gas With Shift" & "NGWS" & No gas candidate units + load shifting allowed
		\\
		"With Gas No Shift" & "WGNS" & Gas candidate units + no load shifting allowed
		\\
		"With Gas With Shift" & "WGWS" & Gas candidate units + load shifting allowed
        \\
		\hline
	\end{tabular}
	\begin{tablenotes}
      \scriptsize
      \item Note: the abbreviations in the second column are used throughout the text.
      \end{tablenotes}
\vspace{-0.2cm}
\end{table}

\subsection{Load Shifting and Gas}\label{gas and shifting}
Table~\ref{investments1 2040} summarizes the investments made per technology type in 2040 under the first four scenarios, the relative change in the objective function value of each scenario compared to the baseline scenario ("NGNS") as well as the total annual load shedding in (CH) and outside (non-CH) Switzerland. In all scenarios the entire biomass candidate capacity is built due to the overall low investment costs reflecting on-going subsidies assumed to remain in the future.
In terms of the RES target, 20.9 TWh come from new and existing PV. The production from existing and newly built biomass units contributes 4.0 TWh and 0.1 TWh are from existing wind generation. The latter is constant across the scenarios. The results show that having more flexibility (i.e. new gas and/or load shifting capability) allows for better utilization of PV and thus less investments required to satisfy the RES target\footnote{In all four scenarios the target of 25 TWh is satisfied but not exceeded.}. Similarly, the investments in grid-batteries are also lower. Apart from leading to cost savings (decreasing objective function value) and less investments in battery storage, load shifting could potentially reduce the amount of gas flexibility required in the system. Comparing the "WGNS" and "WGWS" scenarios, we see a reduction (by 260 MW) in new gas installations. In both scenarios gas turbines are used as peaker power plants. Building new gas power plants instead of solely relying on load shifting reduces the amount of load shedding in Switzerland and leads to the overall lowest objective function value which means that under the given assumptions, investing in gas turbines is economically viable.

\begin{table}[!t]
\renewcommand{\arraystretch}{1.1}
\begin{threeparttable}
	\centering
	\caption{New investments in CH (2040)}
	\vspace{-0.2cm}
	\label{investments1 2040}
	\renewcommand{\arraystretch}{1.2}
	\begin{tabular}{||p{2.2cm}p{1.9cm}p{1.1cm}p{1.7cm}||}
		\hline
		\textbf{Scenario Name} & \textbf{Investment [MW]} & \textbf{Obj. f-n [\%]} &  \textbf{Load Shed [GWh/a]} 
		\\
		\hline
		"NGNS" & Biomass: 240 & -- & Non-CH: 0
		\\
		  & PV: 19'172 &               & CH: 14.8
		\\
		       & Batt: 200 & &
		\\
		\hline
		"WGNS" & Biomass: 240 & $< 1\% \downarrow$ & Non-CH: 0
		\\
		 & PV: 18'097&               & CH: 5.6
		\\
		       & Batt: 100 & &
		\\
		       & Gas (Syn): 500 & &
		\\
		       & Gas (CCS): 255 & &
		\\
		\hline
		"NGWS"& Biomass: 240 & $1.1\% \downarrow$ & Non-CH: 0
		\\
		& PV: 18'141 &               & CH: 7.2
		\\
		       & Batt: 180 & &
		\\
		\hline
     	"WGWS" & Biomass: 240 & $1.24\% \downarrow$ & Non-CH: 0
		\\
		 & PV: 17'807 &               & CH: 2.6
		\\
		       & Batt: 100 & &
		\\
		       & Gas (Syn): 495 & &
		\\
		\hline
	\end{tabular}
\end{threeparttable}

\begin{tablenotes}
  \item Note: The objective function value spans the cost of operation and investments in Switzerland as well as operation in the surrounding countries. The percent change in the objective function value in the third column is calculated using the total objective function value.
\end{tablenotes}
\end{table}

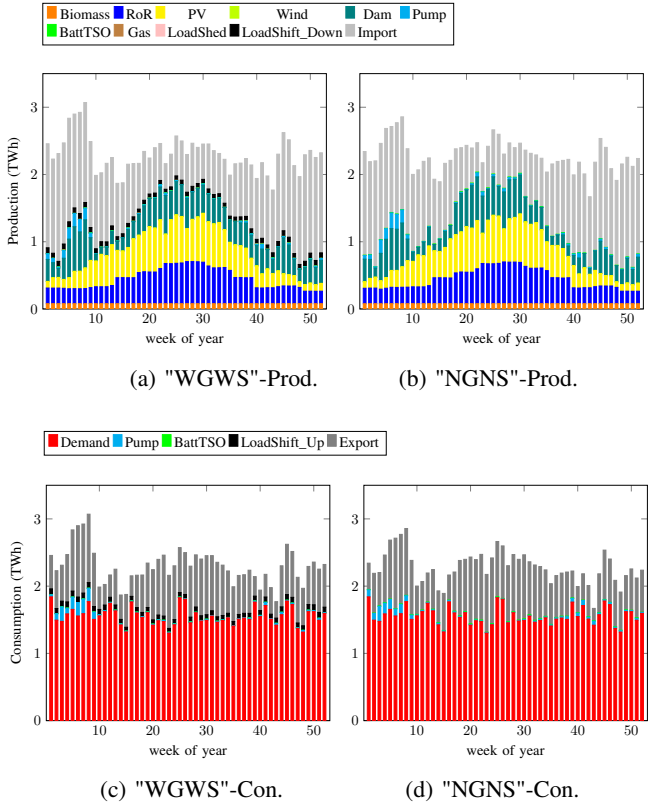
\begin{figure}[!htbp]
  \centering
  \hspace{-0.05cm}
  \subfloat["WGWS"-Prod.]{\input{week_producer_wgws}}
  \hspace{-1.85cm}
  \subfloat["NGNS"-Prod.]{\input{week_producer_ngns}}
  \vspace{0.5cm}
  \subfloat{
  \subfloat["WGWS"-Con.]{\input{week_consumer_wgws}}
  \hspace{-0.90cm}
  \subfloat["NGNS"-Con.]{\input{week_consumer_ngns}}
  }
\caption{Weekly production and consumption per technology type in 2040: "WGWS" (left) vs. "NGNS" (right).}\label{fig:2040_weeklyoperation}
\end{figure}

\begin{figure}[!tb]
  \centering
  \hspace{-0.05cm}
  \subfloat["WGWS"-Prod.]{\input{hour6840plus24_producer_wgws}}
  \hspace{-1.85cm}
  \subfloat["NGNS"-Prod.]{\input{hour6840plus24_producer_ngns}}
  \vspace{0.5cm}
  \subfloat{
  \subfloat["WGWS"-Con.]{\input{hour6840plus24_consumer_wgws}}
  \hspace{-0.90cm}
  \subfloat["NGNS"-Con.]{\input{hour6840plus24_consumer_ngns}}
  }
\caption{Hourly production and consumption per technology type (a day in the first week of October): "WGWS" (left) vs "NGNS" (right).}\label{fig:2040_dailyoperation}
\end{figure}
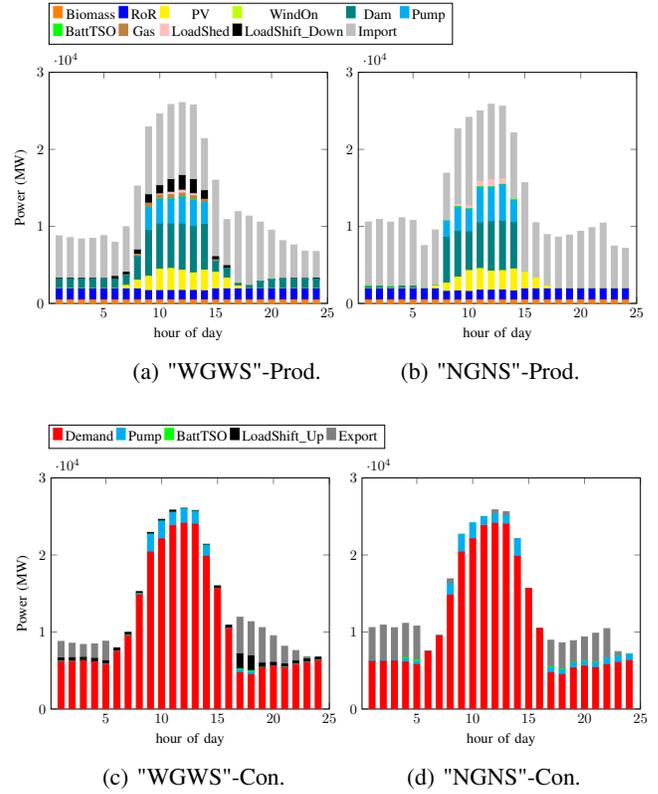

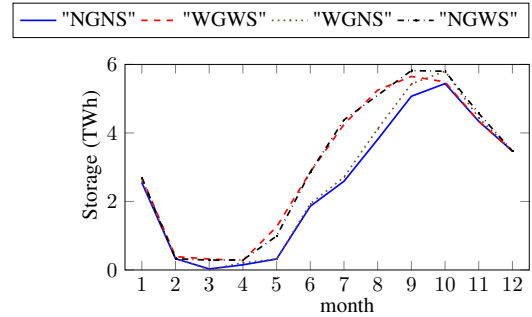
\begin{figure}[!htbp]
\centering
  \input{dam_levels}
  \caption{End-of-month dam storage levels in CH (2040).}
  \label{fig:damlevels}
\end{figure}

Figure~\ref{fig:2040_weeklyoperation} compares the operation of the Swiss generation and storage technologies in the baseline scenario with the operation in the most favorable of the four scenarios, namely "WGWS". The four plots show the weekly stack of generation and consumption in Switzerland in 2040\footnote{In addition to the electricity generation per technology type, we include the imports, downward load shifting as well as load shedding as "production". In addition to the demand and pump load, the consumption stack includes upward load shifting and exports. The total height of each weekly bar in the consumption and generation plots in a given scenario must be equal (i.e. supply equals demand).}. In general, we see similar operation in both scenarios. Pumps are mainly utilized in Jan-Feb to take advantage of hours with very high wind generation in the surrounding countries. The energy is stored and used in subsequent hours due to the lack of base load capacity in the system and low PV output during winter. In order to understand better what drives the investments, Figure~\ref{fig:2040_dailyoperation} shows the hourly production and consumption in a day in October, the month characterized by the highest load shedding in the system. In the highest demand hours, we see load shedding (pink) in both scenarios. However, due to the dispatch of gas this load shedding is lower in the "WGWS" scenario. Furthermore, load shifting is correctly employed to further reduce the need for load shedding or building more gas capacity by increasing the downward shifting in the hours with peak demand and increasing the upward shifting during the rest of the day to satisfy the requirement for equal upward and downward daily shifted energy. An interesting observation that can be made is the lack of production from hydro dams in the off-peak hours in the "NGNS" case. This is due to the lower dam storage levels at the end of September in the "NGNS" scenario as seen in Figure~\ref{fig:damlevels} which plots the end-of-month storage levels in the Swiss dams over the year for all four scenarios. The difference in dam operation whereby more water is stored in the period May-October in the two scenarios with shifting is due to the inherently higher flexibility of the system. 

\begin{comment}
\begin{figure}[!tb]
  %\hspace{1cm}
  %\centering
  %\subfloat[Generation (2020)\label{subfig:2020 prod}]{\input{figs/hour4512plus24_producer_2020}}
  %\hspace{1cm}
  %\subfloat[Consumption (2020)\label{subfig:2020 consum}]{\input{figs/hour4512plus24_consumer_2020}}\\
  %\vspace{0.2cm}
  \subfloat{\input{hour6840plus24_producer_wgws}}
  \subfloat{\input{hour6840plus24_producer_ngns}}\\
  \vspace{0.3cm}
  \subfloat{\input{hour6840plus24_consumer_wgws}}
  \subfloat{\input{hour6840plus24_consumer_ngns}}
  \caption{Hourly generation and consumption per technology type (a day in the first week of October): "WGWS" (left) vs "NGNS" (right).}\label{fig:2040_dailyoperation}
\end{figure}
\end{comment}

\subsection{Reduced Cross-border Transmission Capacity}\label{reduced grid}
As shown in Figure~\ref{fig:2040_weeklyoperation} and Figure~\ref{fig:2040_dailyoperation}, the electricity imports and exports relative to the domestic generation and demand are considerable, highlighting the unique position of Switzerland as a transit country in terms of electricity supply. Furthermore, the ability to export and import due to a reliable grid connectivity is considered as a key source of flexibility in a given power system~\cite{Ulbig2015}. The goal of this section is to investigate the impact of reducing the amount of flexibility stemming from imports/exports on the generation portfolio in the country. 

Similar to Table~\ref{investments1 2040}, Table~\ref{investments 2040 reduced grid} summarizes the investments in the case of reduced (by 70\%) cross-border transfer limits, denoted by (R) and compares them to the baseline scenario ("NGNS") from Section~\ref{gas and shifting}. All four new simulated scenarios show an increase in the objective function. In the runs without gas candidate units, in addition to all wind candidate capacity being built, the amount of PV installations increases by 79\% ("NGNS (R)") and 82\% ("NGWS (R)"). While before the driver for the PV installations was the RES target, here the driver is the lack of available import capacity. Despite the higher RES production, the amount of load shedding both in and outside Switzerland is higher. This means that the reduction in transfer capacity also has a negative impact on the neighbouring countries. Load shifting is beneficial for the reduction of the investments in both gas turbines and battery storages. In the scenarios with gas candidate capacities we observe similar amounts of PV installations as in Table~\ref{investments1 2040}. However, the investment in gas turbines has increased and the operation of these units has changed. 

\begin{table}[!b]
\renewcommand{\arraystretch}{1.1}
\begin{threeparttable}
	\centering
	\caption{New investments in CH (2040) - Reduced Cross-Border Transmission Capacity}
	\vspace{-0.2cm}
	\label{investments 2040 reduced grid}
	\renewcommand{\arraystretch}{1.2}
	\begin{tabular}{||p{2.0cm}p{2.2cm}p{1.0cm}p{1.7cm}||}
		\hline
		\textbf{Scenario Name} & \textbf{Investment [MW]} & \textbf{Obj. f-n [\%]} &  \textbf{Load Shed [GWh/a]} 
		\\
		\hline
		"NGNS" & Biomass: 240 & -- & Non-CH: 0
		\\
		       & PV: 19'172 &               & CH: 14.8
		\\
		       & Batt: 200 & &
		\\
		\hline
		\hline
		"NGNS (R)" & Biomass: 240 & $21.7\% \uparrow$ & Non-CH: 73.3
		\\
		       & PV: 34'111 &               & CH: 17.9
		\\
		       & Wind: 1'960 & &
		\\
		       & Batt: 1'100 & &
		\\
		\hline
		"WGNS (R)" & Biomass: 240 & $18.3\% \uparrow$ & Non-CH: 48.1
		\\
		       & PV: 18'882&               & CH: 5.9
		\\
		       & Wind: 103 & &
		\\
		       & Batt: 200 & &
		\\
		       & Gas (Syn): 1'465 & &
		\\
		       & Gas (CCS): 1'685 & &
		\\
		\hline
		"NGWS (R)" & Biomass: 240 & $16.7\% \uparrow$ & Non-CH: 2.7
		\\
		       & PV: 32'956&               & CH: 9.4
		\\
		       & Wind: 1'960 & &
		\\
		       & Batt: 200  & &
		\\
		\hline
		"WGWS (R)" & Biomass: 240 & $14.6\% \uparrow$ & Non-CH: 2.7
		\\
		       & PV: 17'282 &               & CH: 3.4
		\\
		       & Wind: 103 & &
		\\
		       & Batt: 100 & &
		\\
		        & Gas (Syn): 565 & &
		\\
		       & Gas (CCS): 1'665 & &
		\\
		\hline
	\end{tabular}
\end{threeparttable}
\begin{tablenotes}
  \item Note: The designation (R) pertains to the scenarios with reduced cross-border transmission capacity.
\end{tablenotes}
\end{table}

\begin{figure}[!tb]
  \centering
  \hspace{-0.05cm}
  \subfloat["WGWS (R)"-Prod.]{\input{week_producer_red_wgws}}
  \hspace{-1.85cm}
  \subfloat["NGWS (R)"-Prod.]{\input{week_producer_red_ngws}}
  \vspace{0.5cm}
  \subfloat{
  \subfloat["WGWS (R)"-Con.]{\input{week_consumer_red_wgws}}
  \hspace{-0.90cm}
  \subfloat["NGWS (R)"-Con.]{\input{week_consumer_red_ngws}}
  }
\caption{Weekly production and consumption per technology type in 2040: "WGWS (R)" (left) vs. "NGWS (R)" (right).}\label{fig:2040_weeklyoperation_reduced}
\end{figure}
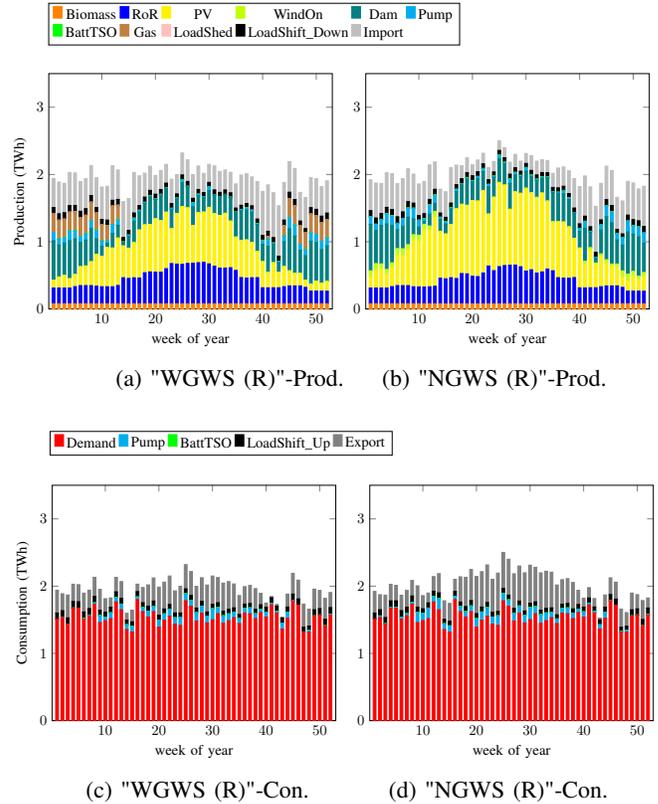

Figure~\ref{fig:2040_weeklyoperation_reduced} compares the weekly generation and consumption per technology type for the "NGWS (R)" and "WGWS (R)" scenarios. In both scenarios we see lower volumes of imports and exports (compared to Figure~\ref{fig:2040_weeklyoperation}) driven by the reduction of the cross-border capacity limits. In "WGWS (R)" the gas turbines with CCS are operated continuously during winter (Jan-Mar and Nov-Dec). These units have higher investment costs, but their total variable costs are significantly lower compared to the syngas turbines which are used as peaker power plants. In the "NGWS (R)" scenario we see increased electricity production from PV in winter due to the very high invested capacity and also wind generation. It is interesting to note that while wind is being dispatched in the winter months, there is very little production during summer. In summer, wind is curtailed due to the oversupply of PV which has a lower total variable cost than wind. In fact the driving factor for investments in both technologies is the reduction in cross-border capacity in winter. The large PV generation in summer results in a net export position of Switzerland in this period. Looking at the generation and consumption of pumped hydro, we see a similar trend in both scenarios - Swiss pumps are partially used to shift water to weeks in the periods of Jan-Feb and Nov-Dec. In both "WGWS (R)" and "NGWS (R)" there is enough available low-cost electricity in summer which is stored in the reservoirs and used during hours with lower available generation and import capacity. While such operation is rather untypical for pumped hydro, it shows that the pumped storage capacity (2.2 TWh) can be used in a variety of ways. 

\begin{comment}
\begin{figure}[!tb]
  \subfloat{\input{hour6840plus24_producer_red_ngws}}
  \hspace{-2cm}
  \subfloat{\input{hour6840plus24_producer_red_wgws}}\\
  \vspace{0.3cm}
  \subfloat{\input{hour6840plus24_consumer_red_ngws}}
  \subfloat{\input{hour6840plus24_consumer_red_wgws}}
  \caption{Hourly generation and consumption per technology type (a day in the first week of October): "NGWS (R)" (left) vs "WGWS (R)" (right).}\label{fig:2040_dailyoperation_reduced}
\end{figure}
\end{comment}

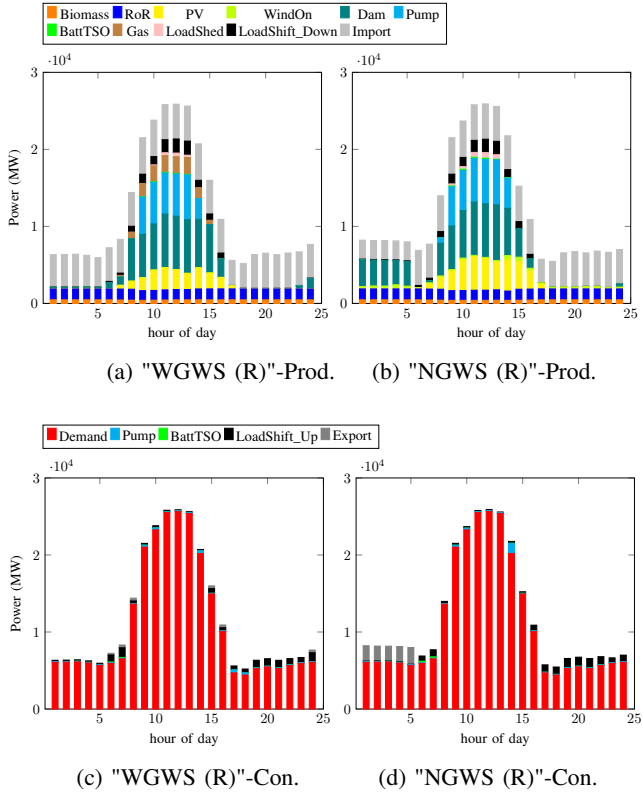
\begin{figure}[!tb]
  \centering
  \hspace{-0.05cm}
  \subfloat["WGWS (R)"-Prod.]{\input{hour6840plus24_producer_red_wgws}}
  \hspace{-1.85cm}
  \subfloat["NGWS (R)"-Prod.]{\input{hour6840plus24_producer_red_ngws}}
  \vspace{0.5cm}
  \subfloat{
  \subfloat["WGWS (R)"-Con.]{\input{hour6840plus24_consumer_red_wgws}}
  \hspace{-0.90cm}
  \subfloat["NGWS (R)"-Con.]{\input{hour6840plus24_consumer_red_ngws}}
  }
\caption{Hourly production and consumption per technology type (a day in the first week of October): "WGWS (R)" (left) vs "NGWS (R)" (right).}\label{fig:2040_dailyoperation_reduced}
\end{figure}

Figure~\ref{fig:2040_dailyoperation_reduced} plots the hourly production and consumption during the same day in October as in Section~\ref{gas and shifting} but now for the "NGWS (R)" and "WGWS (R)" scenarios. Due to the higher investments in PV, there is more PV production in "NGWS (R)" than in Figure~\ref{fig:2040_dailyoperation}. Similarly, due to the higher installed capacity in "WGWS (R)", there is more production from gas. The shifting pattern and the amount of energy shifted are similar. In both scenarios, Switzerland is only importing electricity without simultaneously exporting in the majority of the hours during the day, a direct consequence of the reduced cross-border capacities.

\subsection{Sensitivity Analysis: Rising Gas Prices}\label{sensitivities}
Since many of the model inputs are subject to great uncertainty, in this section we briefly discuss the impact of rising gas prices on the investment decisions in the "WGWS (R)" scenario. Table~\ref{investments 2040 gas sensitivity} summarizes the previously used metrics (i.e. change in objective function value, load shedding and invested capacity per technology type) in case of a 25\%, 50\%, 75\% and 100\% increase in the projected price of gas in 2040. The first important conclusion is that even in case of a 100\% increase in gas prices, the objective function value is still lower than in the "NGWS (R)", and thus this means it makes sense to build gas units even in a high-gas price scenario. This indicates the system's need for flexible capacity\footnote{Whether this flexibility should come from gas turbines or a different technology is subject to a greater discussion.}. When looking at the PV investments, we see that a 25\% increase in the gas price leads to more than 60\% higher installed capacity, highlighting the important role that PV could potentially play in the future. As expected, increasing the gas price leads to a significant reduction in the Gas (CCS) installed capacity as well as  reduction in the total annual production. 
\begin{table}[!tb]
\renewcommand{\arraystretch}{1.1}
\begin{threeparttable}
	\centering
	\caption{New investments in CH (2040) - Rising Gas Prices "WGWS (R)"}
	\vspace{-0.2cm}
	\label{investments 2040 gas sensitivity}
	\renewcommand{\arraystretch}{1.2}
	\begin{tabular}{||p{2.5cm}p{1.0cm}p{1.0cm}p{1.0cm}p{1.0cm}||}
		\hline
		\textbf{Metric} & \textbf{25\%} & \textbf{50\%} &  \textbf{75\%} &  \textbf{100\%}
		\\
		\hline
		Obj. f-n [\%] & $15.2 \uparrow$ & $15.6 \uparrow$ & $15.8 \uparrow$  & $15.9 \uparrow$
		\\
		Biomass [MW] & 240 & 240  & 240 & 240
		\\
		PV [MW] & 28'244 & 30'900 & 32'044 & 33'012
		\\
		Wind [MW] & 1037  & 1960  & 1960 & 1960
		\\
		Battery [MW] & 100  & 69 & 31 & 100
		\\
		Gas (Syn) [MW] &345  & 510 & 570 & 685
		\\
		Gas (CCS) [MW] & 1'105& 590& 500 & 165
		\\
		\hline
	\end{tabular}
\end{threeparttable}
\end{table}

\subsection{Limitations}\label{limitations}
An important limitation of this work is the removal of the UC constraints in the problem formulation. As the amount of intermittent RES generation in the system increases, the minimum up/downtime constraints may become binding. This means that the investments in newly built gas units in the present analysis might be under estimated. Nevertheless, the current formulation is sufficient to show mutual influences between different sources of flexibility and how they could potentially change investment and operational decisions. As already shown in the previous section, the scenario results can differ significantly depending on the input assumptions, therefore we stress the importance of performing sensitivity analyses. Last but not least, in the context of the simulated Swiss power system, we do not model hydro cascades in detail which, in particular, could lead to different operation of hydro pumps. 

\section{Conclusion}\label{conclusion}
In this work a reformulated version of CentIv, the transmission-level generation expansion planning module of the Nexus-e modeling platform, is used to investigate the value of flexibility in a carbon-neutral central European power system in 2040. The power system consists of a detailed representation of Switzerland, for which we establish operational and capacity investment decisions, and an aggregated representation of the four neighboring countries, for which we model the operation only and take the capacities as an input from Euro-Calliope. We simulate two sets of scenarios in order to understand the impact of flexible generation, load shifting and imports on the investment and operation of existing and newly built units. 

Our results show that modeling load shifting within a  generation expansion planning model could reduce the need for investments in both renewable and conventional technologies. In the context of Switzerland, this means that while increasing the number of electric vehicles, parallel efforts should be made to plan the efficient control of these loads. The combination of gas turbines (either running on synthetic gas or with CCS) and load shifting increases the flexibility of the system and results in the overall lowest system costs. While this could be one possible transition path towards carbon-neutrality for Switzerland, it is associated with certain risks such as increasing gas prices which could lead to an almost tenfold reduction in newly built units as presented in the sensitivity analysis. Nevertheless, our results show that there is an inherent need for flexibility in the simulated power system largely dominated by intermittent RES. The reduction of cross-border transmission capacity between Switzerland and the surrounding countries has a profound impact on the domestic operation and investments but could also have an impact on the surrounding countries as witnessed by the increased amount of load shedding. Future work will focus on implementing a discomfort cost associated with load shifting as part of the objective function in order to reduce the amount of simultaneous up and down shifting in the current formulation.

\section*{Acknowledgment}

The authors would like to thank the members of the Nexus-e team and in particular Pranjal Jain for many valuable discussions. Special thanks also goes to Dr. Brynmor Pickering for providing the Euro-Calliope data used in this work.



%\bibliographystyle{IEEEtran}
%\begin{thebibliography}{1}


%\end{thebibliography}

\bibliographystyle{ieeetr}  
\bibliography{BibFile1, BibFile2}

\end{document}

%% file: week_producer_wgws.tex
\pgfplotsset{compat=1.11,
/pgfplots/ybar legend/.style={
/pgfplots/legend image code/.code={%
\draw[##1,/tikz/.cd,yshift=-0.25em]
(0cm,0cm) rectangle (3pt,0.8em);},
   },}
\begin{tikzpicture}[scale=0.55]
\begin{axis}[ybar stacked,
bar width=1.9230769230769231pt,
xlabel=week of year,ylabel=Production (TWh),xmin=0.1,xmax=53,ymin=0,ymax=3.5,
legend style={legend columns=6,at={(1.45,1.3)},},
]
\addplot +[orange] coordinates {
(1,0.07844399361462707)
(2,0.07869512047684153)
(3,0.0786948300670229)
(4,0.07869543754744984)
(5,0.07869543913685662)
(6,0.07869538501608657)
(7,0.07864599034180855)
(8,0.07869536998393421)
(9,0.07869419043235512)
(10,0.07869389189189747)
(11,0.07869333887437978)
(12,0.07869477296887645)
(13,0.07869157376647753)
(14,0.078691855188881)
(15,0.07869451692199829)
(16,0.07869493280780097)
(17,0.07869241077255469)
(18,0.07869372496303555)
(19,0.07869246196037784)
(20,0.07868878595792134)
(21,0.07867446763163875)
(22,0.07865496951360706)
(23,0.07866636003998784)
(24,0.07868109550645795)
(25,0.07869188656968405)
(26,0.07869377834487054)
(27,0.07869369305871009)
(28,0.07869318109770769)
(29,0.07869197652412682)
(30,0.07866592359066654)
(31,0.07869097722135604)
(32,0.07869054430049556)
(33,0.078691155062026)
(34,0.07869470777646058)
(35,0.07869473454473996)
(36,0.0786947323909455)
(37,0.07869472629383313)
(38,0.07868912145327088)
(39,0.07869488249669951)
(40,0.07869497263538229)
(41,0.07868212957866115)
(42,0.07869484555973154)
(43,0.07869483797390839)
(44,0.07869466973627819)
(45,0.0786949631080349)
(46,0.07868725578310592)
(47,0.07869480678694057)
(48,0.0786946360603708)
(49,0.07869487067513987)
(50,0.07869474633944692)
(51,0.07869465956117977)
(52,0.07869488812267247)
};
\addplot +[blue] coordinates {
(1,0.23334527375772482)
(2,0.23221285843288522)
(3,0.23241946222057952)
(4,0.2289534337990989)
(5,0.22217715535503454)
(6,0.22344453383254492)
(7,0.22230024182367533)
(8,0.22521076373715707)
(9,0.23887986753411403)
(10,0.25199773732899605)
(11,0.2517458882831064)
(12,0.2521887602000406)
(13,0.2704184994881418)
(14,0.38631908761399425)
(15,0.38736704401852007)
(16,0.38653874886848955)
(17,0.38634806103320424)
(18,0.46289937377905793)
(19,0.4750945339467663)
(20,0.4715693487806958)
(21,0.4711557646287558)
(22,0.522486745264678)
(23,0.5979271461705097)
(24,0.5954145425127023)
(25,0.6006171844754351)
(26,0.6056528688530306)
(27,0.6258761071605534)
(28,0.6255905629645255)
(29,0.6237722657514935)
(30,0.6171686349609697)
(31,0.560830738630574)
(32,0.5328514362333115)
(33,0.5343852014737829)
(34,0.538078372344016)
(35,0.49543061210384)
(36,0.388802489270064)
(37,0.3888056203555822)
(38,0.3861259025436409)
(39,0.38739431152885323)
(40,0.2362531651798377)
(41,0.2336225287575507)
(42,0.23728811175066802)
(43,0.2368741664594409)
(44,0.25287625257881596)
(45,0.2640170229032321)
(46,0.263136261538554)
(47,0.26408449660779987)
(48,0.24210107836433759)
(49,0.1858711885914115)
(50,0.1859912672180395)
(51,0.1848542299316238)
(52,0.18604430162634766)
};
\addplot +[yellow] coordinates {
(1,0.0962122671986389)
(2,0.15022772063825265)
(3,0.15698489359421441)
(4,0.13038188830246614)
(5,0.16239721294807896)
(6,0.25550706127339007)
(7,0.2577115305289467)
(8,0.30290845688930257)
(9,0.4030721099086636)
(10,0.37904555329688955)
(11,0.48016516665936204)
(12,0.4617184279315878)
(13,0.603831051535777)
(14,0.4048143961962835)
(15,0.39363149641947554)
(16,0.44721318604618193)
(17,0.5769228437176188)
(18,0.553575203469401)
(19,0.5938767085830073)
(20,0.6732510216217031)
(21,0.6508485840131875)
(22,0.7230882336761387)
(23,0.43466202011345956)
(24,0.6504714103897026)
(25,0.7178040890520584)
(26,0.6896086233240677)
(27,0.48940786532060876)
(28,0.6268229666274531)
(29,0.6634942699777626)
(30,0.7251210294682204)
(31,0.6571875320976357)
(32,0.6524112026158998)
(33,0.6723875075519645)
(34,0.5611907702398109)
(35,0.41551367241027837)
(36,0.48511045436144745)
(37,0.471404737478961)
(38,0.4369368056329522)
(39,0.28623698245506735)
(40,0.3242932967472126)
(41,0.20781767268128448)
(42,0.2777298938511674)
(43,0.21045225681838975)
(44,0.2332269076478765)
(45,0.18025729895930265)
(46,0.16612085718442432)
(47,0.1530037986077441)
(48,0.08224630493641485)
(49,0.09592635477033964)
(50,0.11823154814372007)
(51,0.098020289659032)
(52,0.11454290029489224)
};
\addplot +[lime] coordinates {
(1,0.002764685156966848)
(2,0.003208589195879739)
(3,0.002233671255379744)
(4,0.002046259371758591)
(5,0.002107663130521154)
(6,0.0031287101801399995)
(7,0.0028617221446365754)
(8,0.0022775982415542215)
(9,0.0024765298440130507)
(10,0.0030448731189999054)
(11,0.0027509459392174706)
(12,0.0020095299659644104)
(13,0.0035889505445861006)
(14,0.0026095197697015527)
(15,0.0021236858979048067)
(16,0.0028775612804597217)
(17,0.0023184162086771674)
(18,0.0022356078977023458)
(19,0.0025089811921416203)
(20,0.0033280627457970575)
(21,0.0025904095678097217)
(22,0.00240147512528954)
(23,0.002576687769883009)
(24,0.0026270681320351005)
(25,0.0023233825691111793)
(26,0.0025240695792263486)
(27,0.0020249682758517054)
(28,0.0028592704143882365)
(29,0.003068111816783907)
(30,0.0018349956095493809)
(31,0.002166050801204393)
(32,0.001811126097269328)
(33,0.001482265090384765)
(34,0.0022199457682151647)
(35,0.0023650488688418208)
(36,0.0014762549567927756)
(37,0.0030526525491118657)
(38,0.001896801036024362)
(39,0.0020575175354547496)
(40,0.0019005144598472625)
(41,0.0016149014074457418)
(42,0.0017410029175936748)
(43,0.002274395724837109)
(44,0.002145425819360588)
(45,0.0016470301182762415)
(46,0.0009318523245390772)
(47,0.0011495317699653433)
(48,0.002858648794789306)
(49,0.0029386413763281822)
(50,0.002216432121439509)
(51,0.0015599071892319785)
(52,0.001967067591240791)
};
\addplot +[teal] coordinates {
(1,0.3263264410569231)
(2,0.22524029685668331)
(3,0.132358839586943)
(4,0.3146247683129701)
(5,0.5737390436691191)
(6,0.6623956215393834)
(7,0.5849937771620819)
(8,0.7172400256508963)
(9,0.33701508020416454)
(10,0.10268539572566085)
(11,0.10199899519311374)
(12,0.12670356731403232)
(13,0.19393640641709478)
(14,0.13519887722657678)
(15,0.1968892157793299)
(16,0.2669625073957054)
(17,0.25206449070940523)
(18,0.308393594351863)
(19,0.38389837755122597)
(20,0.41695622687879313)
(21,0.43954680067472635)
(22,0.48862025675061127)
(23,0.574721036658335)
(24,0.4200335369026875)
(25,0.4938614207436802)
(26,0.471061937302692)
(27,0.45040295752456666)
(28,0.40427565722939884)
(29,0.4206602976767479)
(30,0.42098294583583173)
(31,0.40574702549062436)
(32,0.36545821784454036)
(33,0.3595087245470083)
(34,0.3293218190077081)
(35,0.32723983135097223)
(36,0.34577378444010354)
(37,0.3387675121032848)
(38,0.37440226378971886)
(39,0.3440881233752166)
(40,0.3010180061644429)
(41,0.37000309736551984)
(42,0.2864761751622779)
(43,0.21166572232751113)
(44,0.35480825306538155)
(45,0.5334582576391946)
(46,0.4435871422181953)
(47,0.244678131104354)
(48,0.18785661042806592)
(49,0.2544282837270793)
(50,0.34868399801379013)
(51,0.25335914334845666)
(52,0.342826728136612)
};
\addplot +[cyan] coordinates {
(1,0.08430023068013807)
(2,0.0626721442558648)
(3,0.02010893463541492)
(4,0.11590522199735286)
(5,0.16663850470443542)
(6,0.19759439619288205)
(7,0.18449854677974012)
(8,0.17276611638892342)
(9,0.06872821119377648)
(10,0.00330129473619901)
(11,0.0018547473959284364)
(12,0.0010590722966360604)
(13,0.00777735696209957)
(14,0.0008815512976093912)
(15,0.0018489782522000322)
(16,0.02031138028078294)
(17,0.005611352487394724)
(18,0.014776213156071383)
(19,0.014802800116643004)
(20,0.0008976259105317487)
(21,0.010864395382108928)
(22,0.03315546058128582)
(23,0.03484845118729022)
(24,0.005705332950352261)
(25,0.016112173189487013)
(26,0.002803702276400483)
(27,0.017230238667783646)
(28,0.007276871317333397)
(29,0.004837062511654211)
(30,0.016085181627098766)
(31,0.007783766349825836)
(32,0.005046697669692132)
(33,0.003705548653421899)
(34,0.0011357527395053988)
(35,0.0010943822153123804)
(36,0.0023920196396789144)
(37,0.02367069311056563)
(38,0.026808356006358215)
(39,0.020137758817585796)
(40,0.018536051367140732)
(41,0.06988961572624723)
(42,0.004827246353341658)
(43,0.016767093745491714)
(44,0.002320206649397039)
(45,0.019558668757435792)
(46,0.025421942995608297)
(47,0.010738128142139201)
(48,0.0005256545220068915)
(49,0.008150365695876865)
(50,0.011944656331265291)
(51,0.019203732495331884)
(52,0.027963211444171493)
};
\addplot +[green] coordinates {
(1,0.0008649707650464384)
(2,0.0014179908366404956)
(3,0.00042424945395333914)
(4,0.0029311530730986514)
(5,0.0014990889690258615)
(6,0.0004978481470722376)
(7,0.001140512483658664)
(8,0.0011840619799927858)
(9,0.002068346030472696)
(10,0.001245105099826798)
(11,0.0013832967284815778)
(12,0.0005032777340392021)
(13,0.0019549870060457117)
(14,0.0019064674878518353)
(15,0.0004309594052164053)
(16,0.000930995540844742)
(17,0.001550482424200957)
(18,0.0010996278488140492)
(19,0.0012958872303914478)
(20,0.0028794485851078755)
(21,0.002944565267751271)
(22,0.002682915064778378)
(23,0.0014205563603602687)
(24,0.00268307770653173)
(25,0.002058639740794751)
(26,0.0010107032407833058)
(27,0.0011652812744141277)
(28,0.000814235720239391)
(29,0.002232147665639718)
(30,0.003335670540743197)
(31,0.002298520583313339)
(32,0.0020031286428343247)
(33,0.001496189784600775)
(34,0.0003676555964647336)
(35,0.00014136367676511745)
(36,0.00021110525881771206)
(37,0.00017799170248968509)
(38,0.0012962518765905196)
(39,0.001197407557012554)
(40,0.0009458606017452352)
(41,0.002206807593663909)
(42,0.0004967750617934618)
(43,0.001873782325612986)
(44,0.00047711619491990414)
(45,0.001331741786067188)
(46,0.0011360939592106151)
(47,0.0004101065426119564)
(48,0.0011224445785649976)
(49,0.00040831739223995086)
(50,0.0004176890950881112)
(51,0.0009800489792286997)
(52,0.0009163430423009739)
};
\addplot +[brown] coordinates {
(1,0.0004224601139346182)
(2,3.747235707620193e-07)
(3,3.818780605938081e-07)
(4,0.0015885693768924703)
(5,0.0010031082788352305)
(6,3.530669104809525e-07)
(7,0.0002543492016615851)
(8,0.0006669754094572968)
(9,4.137376095195798e-07)
(10,4.265794744958353e-07)
(11,4.2814420559856666e-07)
(12,4.096310774041239e-07)
(13,4.375678847514172e-07)
(14,4.3472183624057854e-07)
(15,4.271669324802953e-07)
(16,0.00016797649286656775)
(17,4.500439265433072e-07)
(18,0.0012993242722634166)
(19,0.00024927618769575205)
(20,4.516625662695049e-07)
(21,0.0001647484654688296)
(22,4.5875288708328197e-07)
(23,4.31630117026444e-07)
(24,4.515205855001029e-07)
(25,0.0011117729678092074)
(26,4.476809408906927e-07)
(27,4.4625216688872665e-07)
(28,4.265643910529813e-07)
(29,4.4746222630251606e-07)
(30,4.617023374378971e-07)
(31,4.50283504139581e-07)
(32,4.587074786264975e-07)
(33,4.5454092728670024e-07)
(34,4.4446100487458824e-07)
(35,4.3208676057670757e-07)
(36,4.2283021601365097e-07)
(37,4.2568859081747977e-07)
(38,0.0027705985930949082)
(39,0.0006446075116207451)
(40,0.0027843033852495425)
(41,0.009671913216910667)
(42,4.228905936165292e-07)
(43,0.0008971453266099139)
(44,3.9143106966547405e-07)
(45,0.00020822839691186796)
(46,0.00291754397316028)
(47,3.9032446661714505e-07)
(48,3.905280144850238e-07)
(49,0.0003631493134742471)
(50,3.792292287950369e-07)
(51,0.00012361742351190498)
(52,3.843555844757908e-07)
};
\addplot +[pink] coordinates {
(1,3.705257392496185e-07)
(2,3.5570077256862247e-07)
(3,3.5409279092172153e-07)
(4,3.574124496247595e-07)
(5,3.6121122523123474e-07)
(6,3.5712436877471177e-07)
(7,3.5719398091112885e-07)
(8,3.636078343273299e-07)
(9,3.505741552921645e-07)
(10,3.523221552560042e-07)
(11,3.532417364508134e-07)
(12,3.583121389654002e-07)
(13,3.5351939017476435e-07)
(14,3.4469726441826634e-07)
(15,3.415693628841696e-07)
(16,3.5548988795089515e-07)
(17,3.4627076576190255e-07)
(18,6.311800496898043e-05)
(19,3.4715383866417053e-07)
(20,3.4198894996007114e-07)
(21,3.425871006525272e-07)
(22,3.417627304708942e-07)
(23,3.3840317721192997e-07)
(24,3.3712130684539004e-07)
(25,3.516657528756846e-07)
(26,3.49664239162291e-07)
(27,3.4194243841162165e-07)
(28,3.4456374804054895e-07)
(29,3.398684403565898e-07)
(30,3.3997029212007094e-07)
(31,3.4198437914747315e-07)
(32,3.393202648094463e-07)
(33,3.432403463940165e-07)
(34,3.425879823915792e-07)
(35,3.406857186186931e-07)
(36,3.4323845553936517e-07)
(37,3.4475626503167604e-07)
(38,0.000331133794998804)
(39,3.559122977114902e-07)
(40,0.0003382040738176738)
(41,0.0012140584461161916)
(42,3.490304614929513e-07)
(43,3.4765347866367915e-07)
(44,3.5536222587140044e-07)
(45,3.6380061958819263e-07)
(46,0.0006061666835910091)
(47,3.482462946155109e-07)
(48,3.481877154512384e-07)
(49,3.616912794515074e-07)
(50,3.6047120460874783e-07)
(51,3.554148489895952e-07)
(52,3.588224959248688e-07)
};
\addplot +[black] coordinates {
(1,0.0874287292239934)
(2,0.07546814374773748)
(3,0.08479316627479636)
(4,0.09093827327165895)
(5,0.08826845820289818)
(6,0.08182803373243003)
(7,0.08747447594767427)
(8,0.08271319225958755)
(9,0.08003905843352568)
(10,0.07554499040171568)
(11,0.08178845591381685)
(12,0.07511909427862964)
(13,0.075385080790608)
(14,0.06422707995454856)
(15,0.05747047387617017)
(16,0.06870189998226182)
(17,0.0668246496729698)
(18,0.062208647133551266)
(19,0.06661702407256498)
(20,0.061752866596565655)
(21,0.06695109708915782)
(22,0.060109110361230234)
(23,0.055767961496245695)
(24,0.0601696582133687)
(25,0.06637277288662474)
(26,0.06435189361482874)
(27,0.07166846540394198)
(28,0.07737079845752895)
(29,0.06937479296840766)
(30,0.06436046166426349)
(31,0.07033321155132892)
(32,0.06445391369783385)
(33,0.06254162760207627)
(34,0.0617848205754058)
(35,0.06027959117007782)
(36,0.06410324601372251)
(37,0.0652603968371895)
(38,0.06787984484733961)
(39,0.0713748690475478)
(40,0.07084041055094574)
(41,0.07661441114044432)
(42,0.0645938504204396)
(43,0.07592729032189519)
(44,0.07900512806479429)
(45,0.08417036124046011)
(46,0.08509241504368542)
(47,0.0712005916615962)
(48,0.06673670789106828)
(49,0.0828673382346842)
(50,0.08436768205585282)
(51,0.08493960466311892)
(52,0.08181425132723011)
};
\addplot +[lightgray] coordinates {
(1,1.5465494634905117)
(2,1.4004710853168199)
(3,1.6033059293115708)
(4,1.5036257426088342)
(5,1.5399308149411464)
(6,1.39549152507652)
(7,1.503881439106736)
(8,1.486563592262522)
(9,1.2782090733657634)
(10,1.091050265898926)
(11,1.0245845297331382)
(12,1.1682016610625139)
(13,1.0164738286918618)
(14,0.7937572904099564)
(15,0.7579998366215719)
(16,0.880757038443733)
(17,0.7923150322145909)
(18,0.6804510123301568)
(19,0.7260359597262964)
(20,0.5358628566699638)
(21,0.6746769783965947)
(22,0.5498028408053004)
(23,0.3804176438493893)
(24,0.48641881167283735)
(25,0.5939497676055875)
(26,0.5829209010011256)
(27,0.5609618372438653)
(28,0.6413178081120563)
(29,0.5371096913036174)
(30,0.5253792469962015)
(31,0.6687146871089684)
(32,0.6519179619767346)
(33,0.5982223656668764)
(34,0.6531392076507269)
(35,0.610008829782976)
(36,0.7926016913760072)
(37,0.7997872550908546)
(38,0.8600887853564559)
(39,0.9496936441337164)
(40,0.9054272879367093)
(41,1.1242721831551152)
(42,1.0261672373039037)
(43,0.9336882846850079)
(44,1.298960238878704)
(45,1.4592344959370853)
(46,1.4479653573341769)
(47,1.3074157269918445)
(48,1.3056979617587814)
(49,1.5834978996545535)
(50,1.527242519277322)
(51,1.5311971124945627)
(52,1.487759353356918)
};
\legend{Biomass,RoR,PV,Wind,Dam,Pump,BattTSO,Gas,LoadShed,LoadShift\_Down,Import}
\end{axis}
\end{tikzpicture}

%% file: week_producer_ngns.tex
\pgfplotsset{compat=1.11,
/pgfplots/ybar legend/.style={
/pgfplots/legend image code/.code={%
\draw[##1,/tikz/.cd,yshift=-0.25em]
(0cm,0cm) rectangle (3pt,0.8em);},
   },}
\begin{tikzpicture}[scale=0.55]
\begin{axis}[ybar stacked,
bar width=1.9230769230769231pt,
xlabel=week of year,xmin=0.1,xmax=53,ymin=0,ymax=3.5,
legend style={legend columns=4,at={(1.0,-0.1)},},
]
\addplot +[orange] coordinates {
(1,0.07799019996639965)
(2,0.07870799997184073)
(3,0.07870799996701071)
(4,0.07802690000909898)
(5,0.07843789997180742)
(6,0.07870799997204068)
(7,0.07870799996943084)
(8,0.07870119996799765)
(9,0.07870799995352828)
(10,0.07870799995669628)
(11,0.07870369994862128)
(12,0.07870799995195996)
(13,0.07870799994497445)
(14,0.0787079999588509)
(15,0.07870799996565066)
(16,0.07870799996573966)
(17,0.07870799994804138)
(18,0.07870799995121046)
(19,0.07870799995399963)
(20,0.07862359997615326)
(21,0.07842135601327598)
(22,0.07856129997380297)
(23,0.07866319987890934)
(24,0.07860089995346599)
(25,0.07868499994537818)
(26,0.0787079999478233)
(27,0.0787079999593128)
(28,0.07870799994868072)
(29,0.07870799994998948)
(30,0.07858099994675305)
(31,0.0787079999441195)
(32,0.07866089997374608)
(33,0.07863650006292248)
(34,0.07870799994738582)
(35,0.07870799995809465)
(36,0.07870799996211071)
(37,0.07870799995557436)
(38,0.07839933151361766)
(39,0.07826129996851118)
(40,0.07865354262419419)
(41,0.07846396870140616)
(42,0.07870799996884675)
(43,0.07868239996933687)
(44,0.07870799996706178)
(45,0.07865619997156004)
(46,0.07853699997177582)
(47,0.07870799997013128)
(48,0.07870799996986111)
(49,0.0785120999738757)
(50,0.07862859997141068)
(51,0.0784789002376018)
(52,0.07832769997062307)
};
\addplot +[blue] coordinates {
(1,0.23041708726744664)
(2,0.23214796823108932)
(3,0.2322104582214536)
(4,0.2168431793030008)
(5,0.23029613684701633)
(6,0.2439912842137077)
(7,0.24210565768446896)
(8,0.24463421057423404)
(9,0.24560836810316752)
(10,0.2509459933212059)
(11,0.24991644891993448)
(12,0.2502791446058807)
(13,0.2684369404266588)
(14,0.3854174496627461)
(15,0.3859852329958835)
(16,0.38543917664636385)
(17,0.3828487613136261)
(18,0.45799214271724814)
(19,0.4716492764401464)
(20,0.46835223074507937)
(21,0.4683137946221618)
(22,0.5194703044676853)
(23,0.5977855986598042)
(24,0.5938156186344268)
(25,0.5936897551095567)
(26,0.5993667614187175)
(27,0.6211503704351828)
(28,0.6178284478866451)
(29,0.618825401970852)
(30,0.6148769533462796)
(31,0.5549891527840856)
(32,0.5273304776522751)
(33,0.5279550396641705)
(34,0.5311271654624908)
(35,0.4917231643864205)
(36,0.3866761811803313)
(37,0.3858535886135813)
(38,0.38214108314508266)
(39,0.3861943464389913)
(40,0.23586716418238443)
(41,0.23325553380600378)
(42,0.2376380648293649)
(43,0.2366937745666649)
(44,0.25293085684168254)
(45,0.2636572700456565)
(46,0.26307968204469395)
(47,0.26374355170608144)
(48,0.24267834834055468)
(49,0.18610908309843394)
(50,0.18607307253355146)
(51,0.1849720835931597)
(52,0.18612949915550037)
};
\addplot +[yellow] coordinates {
(1,0.0946112666082795)
(2,0.1416351462487153)
(3,0.1597543996016873)
(4,0.1228540679097609)
(5,0.1568284116191075)
(6,0.24856932458842998)
(7,0.25134966959402183)
(8,0.30156380743664735)
(9,0.39374734851011767)
(10,0.37647513324763743)
(11,0.4719419805229784)
(12,0.4615047086537069)
(13,0.5855563489148183)
(14,0.39943164764828964)
(15,0.38906164631720896)
(16,0.4490349983757473)
(17,0.5777202050011682)
(18,0.5672380644208342)
(19,0.5925753432624181)
(20,0.6768723156769828)
(21,0.6524672829579851)
(22,0.7161997257161122)
(23,0.447977419022935)
(24,0.6452800106365411)
(25,0.7192117916883334)
(26,0.7003446216659898)
(27,0.4881312461382481)
(28,0.6414672083226772)
(29,0.6618501808328953)
(30,0.7166898175341533)
(31,0.6587549371571478)
(32,0.6512649203117431)
(33,0.675922986218299)
(34,0.5659159323412937)
(35,0.42480569877336544)
(36,0.47845080176845184)
(37,0.47840670507503186)
(38,0.4371850185165046)
(39,0.3159663008468955)
(40,0.3305305781078485)
(41,0.20617693340285254)
(42,0.2940006435455512)
(43,0.1983667909549288)
(44,0.24067778605776227)
(45,0.1827672239233647)
(46,0.15415034145543957)
(47,0.1510643526613047)
(48,0.07974071549427088)
(49,0.09434933950480076)
(50,0.11964455580235599)
(51,0.09613888868203668)
(52,0.11744162042527338)
};
\addplot +[lime] coordinates {
(1,0.0027649799873574457)
(2,0.003213559986857248)
(3,0.0021822800072691623)
(4,0.00191084999132851)
(5,0.0017980599576243588)
(6,0.0027759899943995492)
(7,0.002270999565753524)
(8,0.002036680002482838)
(9,0.002244219959174135)
(10,0.003051059985000259)
(11,0.002757139984853621)
(12,0.0020156799852694314)
(13,0.003595139985075135)
(14,0.002615679984864566)
(15,0.002129809985241013)
(16,0.0028836599851789274)
(17,0.00232463998475641)
(18,0.002241229985571366)
(19,0.002515169985926838)
(20,0.003334369985395047)
(21,0.002595599985929455)
(22,0.0024077599857166604)
(23,0.002582959985821841)
(24,0.00263335998551012)
(25,0.0023288499858611857)
(26,0.0025303099855074947)
(27,0.0020312099859010288)
(28,0.0028655399855840825)
(29,0.003074369985623419)
(30,0.001841259985754767)
(31,0.002172299985268709)
(32,0.001817369985170565)
(33,0.0014883799855877622)
(34,0.0022261199852828865)
(35,0.0023712299850731696)
(36,0.0014823599852800387)
(37,0.003058789984831882)
(38,0.0018997899842536315)
(39,0.0020602799845462903)
(40,0.0019086099846766704)
(41,0.001617009985204159)
(42,0.0017471199848064442)
(43,0.002279479984591462)
(44,0.002151729984599463)
(45,0.001652299985400139)
(46,0.0009301599902511316)
(47,0.001148779985839985)
(48,0.0028649899845022064)
(49,0.0029440299849098845)
(50,0.002220759985135361)
(51,0.0015643799849536092)
(52,0.0019726899855900225)
};
\addplot +[teal] coordinates {
(1,0.3356434657351825)
(2,0.2812869186134796)
(3,0.13689015037659744)
(4,0.4074916204466707)
(5,0.5696012249088828)
(6,0.6210412045084166)
(7,0.6121329699107497)
(8,0.6522326260798914)
(9,0.3195579113058483)
(10,0.13843107712707312)
(11,0.11464606845015811)
(12,0.2042110276489503)
(13,0.2794558650974731)
(14,0.09360761095402036)
(15,0.17897967295854122)
(16,0.23817984504767195)
(17,0.1982309347121893)
(18,0.46333989690997374)
(19,0.5731922432533139)
(20,0.5441035565046769)
(21,0.638809088924707)
(22,0.6484392089703235)
(23,0.5523037502647098)
(24,0.4632301635012637)
(25,0.5721801615169541)
(26,0.44205702268719615)
(27,0.5512422853373823)
(28,0.5803329491965465)
(29,0.5666696991712007)
(30,0.5802830999122496)
(31,0.36626942093299536)
(32,0.26963086900833716)
(33,0.316035741232437)
(34,0.3322467534161062)
(35,0.2229649837628018)
(36,0.1304569693797002)
(37,0.15588877798774914)
(38,0.20210123007006636)
(39,0.15456963326576087)
(40,0.14799560485281457)
(41,0.23891285614420488)
(42,0.21251327030027356)
(43,0.083761306021927)
(44,0.29029493558199293)
(45,0.4826393546429766)
(46,0.38782329450130315)
(47,0.292602021829448)
(48,0.24049842983307923)
(49,0.2191754111477407)
(50,0.3674829953332512)
(51,0.22640908235496232)
(52,0.3877204855150551)
};
\addplot +[cyan] coordinates {
(1,0.05433724044883633)
(2,0.0660711108747641)
(3,0.017033165663994435)
(4,0.19869472547726258)
(5,0.1567832849297836)
(6,0.2327680167825851)
(7,0.21989376695928303)
(8,0.19337825936846373)
(9,0.062397633872964896)
(10,0.0009647632042994864)
(11,0.00025725449721531287)
(12,0.0008274569538292408)
(13,0.005619322110507681)
(14,2.28994507392072e-09)
(15,0.0053618330098246245)
(16,0.010280520924317044)
(17,0.0018163044969361339)
(18,0.013634806367326406)
(19,0.016224222939166453)
(20,0.004435409385008468)
(21,0.022511001473790298)
(22,0.05978757678697532)
(23,0.061767801317071935)
(24,0.008233459137040244)
(25,0.032338533502811816)
(26,0.010343793758268438)
(27,0.046427491562161326)
(28,0.03413684614859708)
(29,0.021318326628678524)
(30,0.014799132378157194)
(31,0.0025172555094577433)
(32,0.0038574510286833123)
(33,0.0030870305150039437)
(34,0.0016816953162541715)
(35,0.007051416924588522)
(36,7.500034577713385e-05)
(37,0.0010179784888099106)
(38,0.024825012795552156)
(39,0.028031578638363027)
(40,0.021752383442491488)
(41,0.07111697209674626)
(42,0.0026092070223037766)
(43,0.01635399961338433)
(44,0.001677961647182815)
(45,0.030319080785523863)
(46,0.028007915805847747)
(47,0.011632471920310649)
(48,1.6777868451140938e-10)
(49,0.008685356748179843)
(50,0.014217853248556811)
(51,0.023352580086145647)
(52,0.04094544025733702)
};
\addplot +[green] coordinates {
(1,0.005883454513128393)
(2,0.007166747636216096)
(3,0.00508791745106191)
(4,0.007409975580510194)
(5,0.007497084115136914)
(6,0.007424569367087867)
(7,0.006865094956343023)
(8,0.006640609588215706)
(9,0.009276669737675871)
(10,0.008082438436542598)
(11,0.0096718522285903)
(12,0.0086520015665431)
(13,0.009644984272903168)
(14,0.007985858322674688)
(15,0.006063710800440251)
(16,0.0061323039388724845)
(17,0.009791863655036254)
(18,0.00939718689684448)
(19,0.008671155060348101)
(20,0.010627755740942069)
(21,0.010851860913652306)
(22,0.012240490889952578)
(23,0.006305676214112226)
(24,0.010407466282167214)
(25,0.01027407442194888)
(26,0.009788101001569368)
(27,0.00871046230169823)
(28,0.009064471555479235)
(29,0.00958501730943188)
(30,0.01120781868939258)
(31,0.010044553468969921)
(32,0.011673116311121777)
(33,0.010905125110748143)
(34,0.009351602542703431)
(35,0.0074504259441417475)
(36,0.007157150686749189)
(37,0.010095456563812165)
(38,0.009530191218626173)
(39,0.006562964768093544)
(40,0.0061580987453648895)
(41,0.006729374527352764)
(42,0.006452054640872639)
(43,0.005849540639463197)
(44,0.006646447594676143)
(45,0.006458784682111562)
(46,0.00623458636444093)
(47,0.005591371491245378)
(48,0.007827439657711975)
(49,0.005650725787724882)
(50,0.0061983588493199305)
(51,0.009813935979998065)
(52,0.006325069678548803)
};
\addplot +[pink] coordinates {
(1,0.0002746691873924894)
(2,5.222759721293153e-13)
(3,5.208116873517755e-13)
(4,0.00039464352441488686)
(5,0.0004809932931772097)
(6,5.238368804115806e-13)
(7,5.253721719559702e-13)
(8,5.7380874577501636e-05)
(9,5.455667921900609e-13)
(10,5.404577868451059e-13)
(11,5.669937418766259e-13)
(12,5.594261531845023e-13)
(13,5.576354003845819e-13)
(14,5.377773708194422e-13)
(15,5.25500650811337e-13)
(16,5.373283330468014e-13)
(17,5.669039632350378e-13)
(18,0.0004000998224354856)
(19,5.626647249046736e-13)
(20,5.691251009531229e-13)
(21,0.00015174254706791575)
(22,5.646582386444236e-13)
(23,5.242398469893477e-13)
(24,5.488989725578657e-13)
(25,0.00037872691441679003)
(26,5.782119186555538e-13)
(27,5.39988974962441e-13)
(28,5.511308999601516e-13)
(29,5.553961191640515e-13)
(30,5.551917167883549e-13)
(31,5.628245098819614e-13)
(32,5.667337034112193e-13)
(33,5.728564545237444e-13)
(34,5.61395576896161e-13)
(35,5.422028264577416e-13)
(36,5.361766385155664e-13)
(37,5.448603017530071e-13)
(38,0.002168818927364943)
(39,0.0003981516621821519)
(40,0.0013939775330921655)
(41,0.004893009619163826)
(42,5.234589298134929e-13)
(43,0.00038909594366914635)
(44,5.257815969032025e-13)
(45,0.0003313477975561124)
(46,0.002753664737318523)
(47,5.19596157722208e-13)
(48,5.192117238524636e-13)
(49,0.00021881792068847)
(50,1.7394184164605715e-05)
(51,7.925533388837089e-05)
(52,8.358173477575047e-06)
};
\addplot +[lightgray] coordinates {
(1,1.5399470661811234)
(2,1.3804255524400284)
(3,1.5705414309318664)
(4,1.4289885611455044)
(5,1.483133291187048)
(6,1.2783173089031752)
(7,1.3566144905778375)
(8,1.3780426963298495)
(9,1.271444953107255)
(10,1.1446386799967527)
(11,1.1411513340996384)
(12,1.1893871228118646)
(13,1.0157160152433893)
(14,0.961777737589288)
(15,0.8455843560668943)
(16,0.9951708219293992)
(17,0.9585638656828834)
(18,0.7890205964045122)
(19,0.6552965685148873)
(20,0.6466645485242939)
(21,0.520294669513993)
(22,0.43532830693297064)
(23,0.427076219557783)
(24,0.6174633119539913)
(25,0.6548546166332613)
(26,0.7551726926229588)
(27,0.5076757914873855)
(28,0.5049351267659292)
(29,0.43879878167321895)
(30,0.4476885561938697)
(31,0.6885282945321991)
(32,0.7704877882393453)
(33,0.676962952889056)
(34,0.7202652305665498)
(35,0.8013199390119693)
(36,1.0737240984191463)
(37,1.081348733951355)
(38,1.0716713232662538)
(39,1.2528987238249965)
(40,1.169784008029687)
(41,1.3424013647520483)
(42,1.174756803980629)
(43,1.0467703503153618)
(44,1.3100162886898508)
(45,1.4892887942275812)
(46,1.4819628779323137)
(47,1.2818876868774924)
(48,1.2598287477257255)
(49,1.5609903799183245)
(50,1.4815362911706702)
(51,1.5002128405755883)
(52,1.4192261391208298)
};
%\legend{Biomass,RoR,PV,Wind,Dam,Pump,BattTSO,LoadShed,Import}
\end{axis}
\end{tikzpicture}

%% file: week_consumer_wgws.tex
\pgfplotsset{compat=1.11,
/pgfplots/ybar legend/.style={
/pgfplots/legend image code/.code={%
\draw[##1,/tikz/.cd,yshift=-0.25em]
(0cm,0cm) rectangle (3pt,0.8em);},
   },}
\begin{tikzpicture}[scale=0.55]
\begin{axis}[ybar stacked,
bar width=1.9230769230769231pt,
xlabel=week of year,ylabel=Consumption (TWh),xmin=0.1,xmax=53,ymin=0,ymax=3.5,
legend style={legend columns=5,at={(1.2,1.23)},},
]

\addplot +[red] coordinates {
(1,1.840712)
(2,1.494129)
(3,1.476114)
(4,1.5875089999999998)
(5,1.654178)
(6,1.556578)
(7,1.590569)
(8,1.769512)
(9,1.505817)
(10,1.5583049999999998)
(11,1.6209859999999998)
(12,1.745997)
(13,1.6329259999999999)
(14,1.423814)
(15,1.319497)
(16,1.766367)
(17,1.601767)
(18,1.535702)
(19,1.6053510000000002)
(20,1.4205539999999999)
(21,1.485781)
(22,1.4742389999999999)
(23,1.3022369999999999)
(24,1.4263379999999999)
(25,1.828263)
(26,1.805257)
(27,1.454972)
(28,1.6049339999999999)
(29,1.4817000000000002)
(30,1.490828)
(31,1.5557379999999998)
(32,1.4642439999999999)
(33,1.48929)
(34,1.5291329999999999)
(35,1.4078249999999999)
(36,1.509258)
(37,1.525248)
(38,1.505753)
(39,1.763019)
(40,1.547852)
(41,1.711466)
(42,1.5098209999999999)
(43,1.419563)
(44,1.575134)
(45,1.7789400000000002)
(46,1.7257399999999998)
(47,1.37124)
(48,1.316717)
(49,1.620353)
(50,1.620754)
(51,1.49248)
(52,1.594722)
};
\addplot +[cyan] coordinates {
(1,0.03050522692801044)
(2,0.09546379120227413)
(3,0.2147267192771876)
(4,0.08408929845851376)
(5,0.18234101271435588)
(6,0.1994672697821658)
(7,0.21755262453261687)
(8,0.19901117895645395)
(9,0.10911235085971394)
(10,0.018318300267598883)
(11,0.014883862887790446)
(12,0.012841358146267016)
(13,0.009401471850975028)
(14,0.015463037065431059)
(15,0.010816716051448962)
(16,0.012211049701422288)
(17,0.011049439429084075)
(18,0.009885829809389136)
(19,0.0092612012708635)
(20,0.013130377063851839)
(21,0.014492777493316608)
(22,0.017196979496576066)
(23,0.012230975443198419)
(24,0.012921736300259081)
(25,0.00740586475630164)
(26,0.005163607874005853)
(27,0.011378260435619446)
(28,0.006708193399453777)
(29,0.0066608148815868636)
(30,0.015969326547011217)
(31,0.0057042629038050255)
(32,0.01381053408190866)
(33,0.009926792075071196)
(34,0.005271728221201146)
(35,0.0055175542041110545)
(36,0.005256522221211305)
(37,0.004955828824607528)
(38,0.020096307846630298)
(39,0.015581233594499533)
(40,0.017001920615486634)
(41,0.047752262982243734)
(42,0.007862350121127657)
(43,0.02726574147972476)
(44,0.03637363808695385)
(45,0.010692589881659522)
(46,0.021415630592706205)
(47,0.00942926444940438)
(48,0.025825495289547223)
(49,0.022113066312907135)
(50,0.01708749711947621)
(51,0.030071018920569986)
(52,0.009799548565284569)
};
\addplot +[green] coordinates {
(1,0.0009385665005160499)
(2,0.001711446858067248)
(3,0.0006919691413955965)
(4,0.0030826426964854522)
(5,0.0016664554744931029)
(6,0.0005010276263910388)
(7,0.0013336850359873386)
(8,0.0010041860644339367)
(9,0.0023554981590133885)
(10,0.0014026864382115573)
(11,0.001532111420583769)
(12,0.0006169128228132258)
(13,0.001952677468280802)
(14,0.0019137786665738993)
(15,0.0007901520577546818)
(16,0.0006298099491615686)
(17,0.0019418779892747139)
(18,0.001171317124758966)
(19,0.001255657943380372)
(20,0.003115728591833792)
(21,0.003119845690018987)
(22,0.0033052634925557723)
(23,0.0011434189212175766)
(24,0.003146540993659355)
(25,0.0023699879865153157)
(26,0.0011324741471997527)
(27,0.0012254025514473515)
(28,0.0008652735928911241)
(29,0.0021651903815359543)
(30,0.003867356753707832)
(31,0.002205176345168774)
(32,0.0021524500479777288)
(33,0.001907977395788309)
(34,0.000326864999060772)
(35,0.00016374840625791952)
(36,0.0002151808178264028)
(37,0.00011519605653557199)
(38,0.0013428551126046138)
(39,0.001373735081686867)
(40,0.0012104285135309044)
(41,0.002350488639694114)
(42,0.0004379789005823442)
(43,0.0018157103875824788)
(44,0.0009255343903510652)
(45,0.0013192964124483138)
(46,0.0011059717287396943)
(47,0.0004423321993989797)
(48,0.0014477232343825682)
(49,0.0003714406592268236)
(50,0.0005588598964779375)
(51,0.000969703614074021)
(52,0.0007764582033004055)
};
\addplot +[black] coordinates {
(1,0.08742872922399343)
(2,0.07546814444090268)
(3,0.08479316627479636)
(4,0.09093827327165893)
(5,0.08826845820289822)
(6,0.08182803384795755)
(7,0.08747447594767428)
(8,0.08271319225958755)
(9,0.08003905843352568)
(10,0.07554499040171564)
(11,0.08178845591381675)
(12,0.07511909427862964)
(13,0.0753850809061355)
(14,0.0642270799545485)
(15,0.057470473876170214)
(16,0.06870189998226185)
(17,0.06682464967296992)
(18,0.062208647133551294)
(19,0.06661702407256498)
(20,0.061752866596565606)
(21,0.06695109708915774)
(22,0.06010911036123008)
(23,0.055767961496245674)
(24,0.06016965821336874)
(25,0.06637277288662478)
(26,0.06435189361482871)
(27,0.07166846540394198)
(28,0.07737079845752902)
(29,0.06937479296840761)
(30,0.06436046166426349)
(31,0.0703332115513289)
(32,0.06445391369783382)
(33,0.06254162760207634)
(34,0.0617848205754058)
(35,0.06027959117007781)
(36,0.06410324601372253)
(37,0.06526039683718948)
(38,0.06787984484733955)
(39,0.0713748690475478)
(40,0.07084041055094577)
(41,0.07661441114044416)
(42,0.06459385042043958)
(43,0.07592729032189519)
(44,0.07900512806479425)
(45,0.08417036124046016)
(46,0.08509241504368543)
(47,0.07120059166159629)
(48,0.06673670789106831)
(49,0.0828673382346843)
(50,0.08436768205585289)
(51,0.08493960466311894)
(52,0.08181425167381275)
};
\addplot +[gray] coordinates {
(1,0.4970741901929239)
(2,0.5628421255586564)
(3,0.5349986849315425)
(4,0.704071717683551)
(5,0.9100027514658677)
(6,1.06020932146768)
(7,1.0268329840475563)
(8,1.0179857858453876)
(9,0.791859150374443)
(10,0.33303873580300625)
(11,0.30577554204455715)
(12,0.3316243923388123)
(13,0.5323931219164132)
(14,0.3629888349316781)
(15,0.4878824600080185)
(16,0.3052466485164225)
(17,0.4810653937521041)
(18,0.5567274783932312)
(19,0.6605872994461187)
(20,0.7466338899250822)
(21,0.8280732580492549)
(22,0.9061522786951982)
(23,0.7896291031095802)
(24,0.7996292118812203)
(25,0.6684916400408858)
(26,0.6227241235991223)
(27,0.7581878987809403)
(28,0.7751436822672704)
(29,0.8433404300502791)
(30,0.8779095716179356)
(31,0.8197724762257849)
(32,0.8099839542514639)
(33,0.7487548109033608)
(34,0.6294172499324572)
(35,0.51698277052461)
(36,0.5803334201071338)
(37,0.5740427599053047)
(38,0.6421536827586652)
(39,0.29017144844991566)
(40,0.30412713936254115)
(41,0.33742598229757087)
(42,0.39530055716264423)
(43,0.24454340768881944)
(44,0.6110764714015054)
(45,0.747456011628791)
(46,0.6822486982300933)
(47,0.6790636953116325)
(48,0.5571136866422157)
(49,0.5674417529588944)
(50,0.635023066348206)
(51,0.6444722009814856)
(52,0.6354173569070966)
};
\legend{Demand,Pump,BattTSO,LoadShift\_Up,Export}
\end{axis}
\end{tikzpicture}

%% file: week_consumer_ngns.tex
\pgfplotsset{compat=1.11,
/pgfplots/ybar legend/.style={
/pgfplots/legend image code/.code={%
\draw[##1,/tikz/.cd,yshift=-0.25em]
(0cm,0cm) rectangle (3pt,0.8em);},
   },}
\begin{tikzpicture}[scale=0.55]
\begin{axis}[ybar stacked,
bar width=1.9230769230769231pt,
xlabel=week of year,xmin=0.1,xmax=53,ymin=0,ymax=3.5,
legend style={legend columns=4,at={(1.0,-0.1)},},
]
\addplot +[red] coordinates {
(1,1.840712)
(2,1.494129)
(3,1.476114)
(4,1.5875089999999998)
(5,1.654178)
(6,1.556578)
(7,1.590569)
(8,1.769512)
(9,1.505817)
(10,1.5583049999999998)
(11,1.6209859999999998)
(12,1.745997)
(13,1.6329259999999999)
(14,1.423814)
(15,1.319497)
(16,1.766367)
(17,1.601767)
(18,1.535702)
(19,1.6053510000000002)
(20,1.4205539999999999)
(21,1.485781)
(22,1.4742389999999999)
(23,1.3022369999999999)
(24,1.4263379999999999)
(25,1.828263)
(26,1.805257)
(27,1.454972)
(28,1.6049339999999999)
(29,1.4817000000000002)
(30,1.490828)
(31,1.5557379999999998)
(32,1.4642439999999999)
(33,1.48929)
(34,1.5291329999999999)
(35,1.4078249999999999)
(36,1.509258)
(37,1.525248)
(38,1.505753)
(39,1.763019)
(40,1.547852)
(41,1.711466)
(42,1.5098209999999999)
(43,1.419563)
(44,1.575134)
(45,1.7789400000000002)
(46,1.7257399999999998)
(47,1.37124)
(48,1.316717)
(49,1.620353)
(50,1.620754)
(51,1.49248)
(52,1.594722)
};
\addplot +[cyan] coordinates {
(1,0.10351631201866039)
(2,0.09886956277787934)
(3,0.22390034484292437)
(4,0.15911703742097394)
(5,0.1467387281037824)
(6,0.10307113071407716)
(7,0.13849532097375866)
(8,0.08682188359421728)
(9,0.06309062175043272)
(10,0.01380713498748654)
(11,0.014794401524153494)
(12,0.009161890117237234)
(13,0.015090331416821758)
(14,0.02516432885741448)
(15,0.01464116135169448)
(16,0.020632662408233908)
(17,0.016025697106723545)
(18,0.014735826112717221)
(19,0.00672369223013952)
(20,0.008585694053178682)
(21,0.011424730756620058)
(22,0.008675834618027142)
(23,0.006331477327369169)
(24,0.007142621867264001)
(25,0.005167975319756762)
(26,0.002859352009462328)
(27,0.007715291091653982)
(28,0.006523267825264304)
(29,0.004446777251647427)
(30,0.012147985216913066)
(31,0.009118085057101915)
(32,0.017153564685745533)
(33,0.018517611179910295)
(34,0.012210864215070689)
(35,0.01519087456416279)
(36,0.01825239417369648)
(37,0.019316903914729087)
(38,0.038779129557031486)
(39,0.05465850896958291)
(40,0.05114459416973221)
(41,0.07511334861040819)
(42,0.02807227686984032)
(43,0.059752689557147444)
(44,0.023725480669497125)
(45,0.016105430862615264)
(46,0.027259400330207257)
(47,0.01859798823869492)
(48,0.019841142618748205)
(49,0.02076715640933756)
(50,0.015138885307927055)
(51,0.028855824975602167)
(52,0.01192569472353132)
};
\addplot +[green] coordinates {
(1,0.0070376202557756815)
(2,0.0075567831640914835)
(3,0.005917716113341476)
(4,0.007495065047269893)
(5,0.008284879590530562)
(6,0.00817268551766029)
(7,0.007518055578302813)
(8,0.007384310101975403)
(9,0.009851033786940567)
(10,0.008816258298316651)
(11,0.010358854166284003)
(12,0.009725678456868683)
(13,0.010465477725983649)
(14,0.008186421847133382)
(15,0.006876428813076501)
(16,0.006419216435728406)
(17,0.011041518750853765)
(18,0.00972920420180357)
(19,0.009259882126466931)
(20,0.011948281886751594)
(21,0.011350019901464016)
(22,0.01347264830397228)
(23,0.007093832400082966)
(24,0.01075028919780084)
(25,0.011872909937880004)
(26,0.010193342682083831)
(27,0.009357024461827313)
(28,0.009987795755661171)
(29,0.010447479891301242)
(30,0.011848434000108607)
(31,0.011317426454538522)
(32,0.012375345632149315)
(33,0.012015106829404314)
(34,0.010038634361581823)
(35,0.00851825414622675)
(36,0.007766005517323759)
(37,0.010754271445308985)
(38,0.010143506802798347)
(39,0.0072208876035848805)
(40,0.006678256652010679)
(41,0.0072754209523167215)
(42,0.006989201302619894)
(43,0.006363802754637865)
(44,0.007534858044204509)
(45,0.00680823346496434)
(46,0.006676517463745335)
(47,0.0062000000440064935)
(48,0.008393494578764064)
(49,0.006374880445036063)
(50,0.006645292235489394)
(51,0.010585882136800322)
(52,0.006809519889806769)
};
\addplot +[gray] coordinates {
(1,0.3906034976190472)
(2,0.5900996580599008)
(3,0.4964757412635533)
(4,0.7084934209176479)
(5,0.8756547791336278)
(6,1.0457738820969869)
(7,1.033358272664681)
(8,0.993569276524522)
(9,0.8042264490112537)
(10,0.42036875198829654)
(11,0.4229065229604645)
(12,0.43070057360280883)
(13,0.5882508068518986)
(14,0.47237923570501605)
(15,0.5508596719337887)
(16,0.37241044796821177)
(17,0.581170358935969)
(18,0.8218049931597765)
(19,0.777497405052501)
(20,0.9919258105974987)
(21,0.8858606462928089)
(22,0.9760471908004316)
(23,0.858800315172561)
(24,0.9754333790182337)
(25,0.8186376244592192)
(26,0.7800016083953908)
(27,0.8320325416526716)
(28,0.8478935262280984)
(29,0.9022355203778335)
(30,0.9511432187684783)
(31,0.785810402801491)
(32,0.8209499821914298)
(33,0.7711710376678211)
(34,0.6901400010003151)
(35,0.6048607300349467)
(36,0.6214541620354047)
(37,0.6390588552595983)
(38,0.6552461630758439)
(39,0.400044882823515)
(40,0.38836911667916646)
(41,0.38971225347055066)
(42,0.4635426860990617)
(43,0.18346724569592446)
(44,0.5767096676499842)
(45,0.733916691732473)
(46,0.6438036050077905)
(47,0.6903402481580269)
(48,0.5671950339748443)
(49,0.5091402072286655)
(50,0.6134817035333598)
(51,0.5891002397142853)
(52,0.6246397876672521)
};
%\legend{Demand,Pump,BattTSO,Export}
\end{axis}
\end{tikzpicture}

%% file: hour6840plus24_producer_wgws.tex
\pgfplotsset{compat=1.11,
/pgfplots/ybar legend/.style={
/pgfplots/legend image code/.code={%
\draw[##1,/tikz/.cd,yshift=-0.25em]
(0cm,0cm) rectangle (3pt,0.8em);},
   },}
\begin{tikzpicture}[scale=0.54]
\begin{axis}[ybar stacked,
bar width=4.166666666666667pt,
xlabel=hour of day, ylabel=Power (MW),xmin=0.1,xmax=25,ymin=0,ymax=30000,
legend style={legend columns=6,at={(1.45,1.3)},},
]

\addplot +[orange] coordinates {
(1,468.42163475058186)
(2,468.4216387775112)
(3,468.4216484399053)
(4,468.4215968369787)
(5,468.4215463890967)
(6,468.42157312432397)
(7,468.42174780014733)
(8,468.42073949177495)
(9,465.8761088410698)
(10,463.8414585818237)
(11,466.4313940903431)
(12,468.46705729802073)
(13,468.46735647571046)
(14,468.4470537670033)
(15,468.4220809205542)
(16,468.42163422730334)
(17,468.41590212388974)
(18,468.41592108537054)
(19,468.4213591453252)
(20,468.42144161248865)
(21,468.42148974218514)
(22,468.42159544762603)
(23,468.42166969818413)
(24,468.4217105155034)
};
\addplot +[blue] coordinates {
(1,1413.1301193550323)
(2,1413.1302026278724)
(3,1413.1321882697482)
(4,1413.1313351142558)
(5,1413.130544107149)
(6,1413.1232054751795)
(7,1413.1224644325525)
(8,1413.1027739280664)
(9,1201.4693298410227)
(10,1201.1527348694483)
(11,1209.515506513889)
(12,1213.8703005990005)
(13,1213.8873675291356)
(14,1172.5160090210245)
(15,1413.118911342434)
(16,1413.1197994915)
(17,1410.6656229124899)
(18,1410.6919354334066)
(19,1413.1463491028514)
(20,1413.1311159343609)
(21,1413.1293912064395)
(22,1413.131000425705)
(23,1413.13146878656)
(24,1413.1307298088122)
};
\addplot +[yellow] coordinates {
(1,9.998015683787491)
(2,10.553280718884096)
(3,9.930159500394888)
(4,9.292618492426831)
(5,9.01631312429509)
(6,15.25876015686314)
(7,417.242785403821)
(8,1109.9650726648615)
(9,1804.0724226813477)
(10,2749.74207253527)
(11,2817.8965637920214)
(12,2593.9841046870292)
(13,2267.1673700529595)
(14,2669.772216481978)
(15,2122.009600993991)
(16,1389.6664544618327)
(17,359.6136872647152)
(18,9.706763531176149)
(19,7.7693243060569435)
(20,8.146875933222336)
(21,8.25537501805213)
(22,8.585259122327697)
(23,8.748230359623946)
(24,9.384562932199767)
};
\addplot +[lime] coordinates {
(1,8.533338577956629)
(2,11.01329630530232)
(3,10.633391073108514)
(4,15.733200493205072)
(5,10.16315144058563)
(6,4.144025548322277)
(7,5.723991972358934)
(8,6.883862079903061)
(9,6.474067459997745)
(10,7.924005714439032)
(11,4.741504671769782)
(12,5.591657525679198)
(13,6.672120705929259)
(14,17.754116528061395)
(15,17.30330156209264)
(16,12.073290216318341)
(17,6.34879723118745)
(18,6.918659377833618)
(19,7.943176759813862)
(20,7.743284653122831)
(21,9.293302310155632)
(22,9.293501657218352)
(23,4.384088055582535)
(24,6.793787915118815)
};
\addplot +[teal] coordinates {
(1,1255.901291721153)
(2,1253.2254604072966)
(3,1253.753502694626)
(4,1245.963708481213)
(5,1234.9566257665742)
(6,1258.2517766892252)
(7,1285.350699471445)
(8,3065.05060856469)
(9,5987.455875359132)
(10,5890.863250209742)
(11,5813.745852656658)
(12,6021.50831174476)
(13,6030.767873360681)
(14,5930.084735309916)
(15,1423.9806126366539)
(16,1272.7699231399083)
(17,296.72489193801573)
(18,378.552430149888)
(19,931.7364008897906)
(20,1182.456498519507)
(21,1225.9366623119324)
(22,1238.8380244010298)
(23,1237.2528410399323)
(24,1233.7794734284103)
};
\addplot +[cyan] coordinates {
(1,4.573711640564774)
(2,4.576005096956667)
(3,4.579521187392088)
(4,4.577647168608568)
(5,4.576399830826567)
(6,4.576860676333484)
(7,4.57922288765042)
(8,4.650710081378769)
(9,2953.85158849454)
(10,3177.358105791422)
(11,3230.8349566445354)
(12,3442.3959201240314)
(13,3405.158818606282)
(14,2858.336645050698)
(15,4.542423156612481)
(16,4.580768818958239)
(17,0.7086469020297641)
(18,0.7084589501346922)
(19,4.447351776319536)
(20,4.56512785688399)
(21,4.584489735347559)
(22,4.5900407674966495)
(23,4.588104122044337)
(24,4.585476359483282)
};
\addplot +[green] coordinates {
(1,0.0839037328393624)
(2,0.0846762032074742)
(3,0.0856065821509383)
(4,0.0867874554127303)
(5,0.0881650635604336)
(6,0.088650420419412)
(7,0.0865706986525835)
(8,0.107260486298522)
(9,64.167945687973)
(10,77.8781191885306)
(11,79.1329903702771)
(12,80.1970702135286)
(13,82.6247924879659)
(14,0.0037861887952651)
(15,0.0828112402513618)
(16,0.0914993905997549)
(17,0.0305396896394577)
(18,0.0301178030349295)
(19,0.0949434235963027)
(20,0.0876952958598783)
(21,0.0872367099122928)
(22,0.0872094458453239)
(23,0.0873781082903193)
(24,0.0876079365129386)
};
\addplot +[brown] coordinates {
(1,0.0025213699876196)
(2,0.00251871023931263)
(3,0.00250803998033709)
(4,0.00247157887610827)
(5,0.00247622935231794)
(6,0.0022992795376985)
(7,102.962595030215)
(8,298.980879219143)
(9,494.99993928352)
(10,494.999948264483)
(11,494.99994815701)
(12,494.999948900067)
(13,494.99994936126)
(14,388.100092892221)
(15,192.08144516412)
(16,0.00324136860271921)
(17,0.0023686449639538)
(18,0.00244928766288562)
(19,0.00253718081500679)
(20,0.00254245673068411)
(21,0.00254592796746012)
(22,0.00253811347393883)
(23,0.00253053154563864)
(24,0.00252539726205012)
};
\addplot +[pink] coordinates {
(1,0.00199086626778932)
(2,0.00199236767034938)
(3,0.00199629520630682)
(4,0.00198616898448216)
(5,0.00196766235106524)
(6,0.00205252295428626)
(7,0.00215469198856127)
(8,0.00233641786074931)
(9,31.1469473279957)
(10,148.485713402082)
(11,313.775913527243)
(12,379.841493783465)
(13,260.152764756778)
(14,0.00245020768226525)
(15,0.00235147043989311)
(16,0.00217456300700039)
(17,0.00189730389044831)
(18,0.00189221612421773)
(19,0.00194069789595419)
(20,0.00195542337524682)
(21,0.00194700519621239)
(22,0.00197163716827946)
(23,0.00198951401073853)
(24,0.00200347615606521)
};
\addplot +[black] coordinates {
(1,87.5071345554948)
(2,87.853055105944)
(3,88.2587869476419)
(4,85.6134931240356)
(5,82.9066613094527)
(6,333.252594062651)
(7,343.73309172035)
(8,569.194509730234)
(9,1087.89736998888)
(10,1065.4194644943)
(11,1632.77704052498)
(12,1900.84121991601)
(13,1835.77252047043)
(14,1118.70313627405)
(15,422.61962305945)
(16,339.619211936395)
(17,6.86701498365378)
(18,7.05133652348014)
(19,67.1747249964131)
(20,79.5211706548082)
(21,88.9104766997345)
(22,90.2402971919203)
(23,91.4027854849078)
(24,126.391924068125)
};
\addplot +[lightgray] coordinates {
(1,5531.470711416315)
(2,5289.844038738399)
(3,5132.361156724161)
(4,5212.921894821433)
(5,5575.219744552898)
(6,4444.864934824325)
(7,5915.588054173368)
(8,8305.53076863073)
(9,8828.32674529217)
(10,9341.97658446321)
(11,9755.291219966002)
(12,9467.940057236849)
(13,9704.513322357423)
(14,6766.602571864465)
(15,9932.260005712222)
(16,5981.854867226998)
(17,9364.313448100405)
(18,9058.080172323853)
(19,7665.587903770292)
(20,6347.557624401243)
(21,4925.936831279692)
(22,4361.092434682723)
(23,3544.5695787299746)
(24,3476.8458265622303)
};
\legend{Biomass,RoR,PV,WindOn,Dam,Pump,BattTSO,Gas,LoadShed,LoadShift\_Down,Import}
\end{axis}
\end{tikzpicture}

%% file: hour6840plus24_producer_ngns.tex
\pgfplotsset{compat=1.11,
/pgfplots/ybar legend/.style={
/pgfplots/legend image code/.code={%
\draw[##1,/tikz/.cd,yshift=-0.25em]
(0cm,0cm) rectangle (3pt,0.8em);},
   },}
\begin{tikzpicture}[scale=0.54]
\begin{axis}[ybar stacked,
bar width=4.166666666666667pt,
xlabel=hour of day,xmin=0.1,xmax=25,ymin=0,ymax=30000,
legend style={legend columns=4,at={(1.0,-0.1)},},
]
\addplot +[orange] coordinates {
(1,468.49999981869604)
(2,468.4999998186883)
(3,468.4999998187356)
(4,468.49999981029373)
(5,468.4999998104847)
(6,468.4999998256997)
(7,468.4999998264688)
(8,468.4999993620769)
(9,445.49999998079466)
(10,462.70000009295495)
(11,468.4999999728508)
(12,468.49999998811114)
(13,468.4999999893613)
(14,444.00000008393283)
(15,468.49999980012484)
(16,468.49999982580584)
(17,468.4999998066405)
(18,468.49999980647834)
(19,468.4999998067402)
(20,468.4999998069702)
(21,468.4999998070372)
(22,468.49999980733793)
(23,468.49999980757235)
(24,468.49999980772026)
};
\addplot +[blue] coordinates {
(1,1416.0955135061265)
(2,1416.0955135040247)
(3,1416.0955135075897)
(4,1416.0955104404559)
(5,1416.0955105724156)
(6,1416.095514891311)
(7,1416.0955149119507)
(8,1095.7745117259774)
(9,1168.2666231465566)
(10,1086.3065934215465)
(11,1243.5889180822787)
(12,1241.674671661211)
(13,1271.7836316008822)
(14,1162.7338696277145)
(15,1416.095514929636)
(16,1416.0955149026854)
(17,1416.0955078132577)
(18,1416.0955078109646)
(19,1416.0955078270197)
(20,1416.0955078405955)
(21,1416.095507843149)
(22,1416.0955078642826)
(23,1416.0955078767881)
(24,1416.0955078807317)
};
\addplot +[yellow] coordinates {
(1,12.698468793750216)
(2,13.318124793205026)
(3,12.760263072918796)
(4,12.289953673879225)
(5,11.486663442829744)
(6,17.638491910152)
(7,393.3606223084425)
(8,1054.5293244489853)
(9,1770.1620462745718)
(10,2680.522478188294)
(11,2776.0189883652633)
(12,2467.0649252770477)
(13,2484.6522492740974)
(14,2850.9662048503114)
(15,2140.065690146526)
(16,1378.9565898453413)
(17,366.74858767666484)
(18,12.321157431825707)
(19,10.489251521839678)
(20,10.75712942723911)
(21,10.983783972487613)
(22,11.319784690796144)
(23,11.534852458697168)
(24,11.928189710729068)
};
\addplot +[lime] coordinates {
(1,8.569999905944872)
(2,11.049999905300849)
(3,10.669999905436859)
(4,15.769999893521803)
(5,10.199999894855894)
(6,4.179999916727811)
(7,5.7599999155607975)
(8,6.91999992735302)
(9,6.479999985133559)
(10,7.92999997821835)
(11,4.939999658157604)
(12,5.819999980276966)
(13,6.949999993216737)
(14,17.069999994888974)
(15,17.339999913328594)
(16,12.109999913334002)
(17,6.389999891530467)
(18,6.959999891050873)
(19,7.979999890400014)
(20,7.779999890587468)
(21,9.329999889908999)
(22,9.32999989008912)
(23,4.419999894437319)
(24,6.829999891634344)
};
\addplot +[teal] coordinates {
(1,241.59666490615842)
(2,287.71925243322175)
(3,146.97026396907734)
(4,323.514228599102)
(5,307.33494880780245)
(6,10.986467676022702)
(7,11.006332169609912)
(8,5930.99999553078)
(9,5948.210191314135)
(10,5053.604586068201)
(11,5940.8395789575725)
(12,6417.118371863088)
(13,6435.800000801472)
(14,6010.1645554001425)
(15,11.152034414423122)
(16,11.046006395275098)
(17,1.4440992553390771e-05)
(18,1.4441254568568063e-05)
(19,1.4440996227239468e-05)
(20,1.444118862245648e-05)
(21,1.4439982952517556e-05)
(22,1.4444765818053995e-05)
(23,1.444296040808294e-05)
(24,1.444136068532332e-05)
};
\addplot +[cyan] coordinates {
(1,8.450323587352117e-07)
(2,8.44788186112417e-07)
(3,8.449219853458776e-07)
(4,6.62433347275341e-07)
(5,6.675954719756686e-07)
(6,1.1634744194542797e-06)
(7,1.1627950415180121e-06)
(8,2135.3110155592444)
(9,3076.358679175454)
(10,2895.303776514852)
(11,4557.706651491562)
(12,4439.738051897484)
(13,4727.385081159848)
(14,2846.1862370318167)
(15,1.1567369429501152e-06)
(16,1.163329389637396e-06)
(17,5.804508704544342e-07)
(18,5.804907105256116e-07)
(19,5.803870934373719e-07)
(20,5.803606686227868e-07)
(21,5.802898641165452e-07)
(22,5.801311272153229e-07)
(23,5.798887511430138e-07)
(24,5.798819123263613e-07)
};
\addplot +[green] coordinates {
(1,147.7873024458848)
(2,79.7809788099733)
(3,141.0717175817245)
(4,1.7874609412649622e-07)
(5,1.823726946715048e-07)
(6,7.16008367597566e-07)
(7,2.99821863764259e-07)
(8,2.8186196099989997e-08)
(9,173.1531607612194)
(10,124.089915045993)
(11,164.31416459140848)
(12,107.3413346940245)
(13,68.93226874946093)
(14,146.2636577132975)
(15,3.00408783655157e-07)
(16,7.52950136739243e-07)
(17,2.84413327698394e-07)
(18,2.8250082355286e-07)
(19,2.8090321100775896e-07)
(20,2.80019584600682e-07)
(21,2.79592696350246e-07)
(22,2.78721394610559e-07)
(23,2.78273465390145e-07)
(24,2.7815191943785e-07)
};
\addplot +[pink] coordinates {
(1,3.05986748182077e-09)
(2,3.06026562817116e-09)
(3,3.06152611626332e-09)
(4,3.05674089587884e-09)
(5,3.05042390184004e-09)
(6,3.08268463857049e-09)
(7,3.11724199546217e-09)
(8,3.20350120637137e-09)
(9,290.652119282977)
(10,436.205136920205)
(11,698.815069019375)
(12,900.514939565386)
(13,681.532921115793)
(14,31.287683661131)
(15,3.20079407641823e-09)
(16,3.12012119844897e-09)
(17,3.02919291306435e-09)
(18,3.02864335704459e-09)
(19,3.03787983355135e-09)
(20,3.04361856777174e-09)
(21,3.04255864922541e-09)
(22,3.05092692417364e-09)
(23,3.05637541829023e-09)
(24,3.06023348186317e-09)
};
\addplot +[lightgray] coordinates {
(1,8275.309516691914)
(2,8617.252319191319)
(3,8363.482398515838)
(4,8881.852437939197)
(5,8547.308010743129)
(6,5588.39940025868)
(7,7245.0774063656345)
(8,6217.225878822725)
(9,9804.517180285646)
(10,11407.6257496995)
(11,9134.218004531413)
(12,9806.786421090555)
(13,9465.005465909953)
(14,8625.62779184404)
(15,11609.646636276158)
(16,7190.091763549868)
(17,6671.913179506813)
(18,6699.519089825928)
(19,6967.205385707525)
(20,7459.12999772716)
(21,7932.670802498666)
(22,8532.935720920834)
(23,5553.315261402675)
(24,5249.735172082937)
};
%\legend{Biomass,RoR,PV,WindOn,Dam,Pump,BattTSO,LoadShed,Import}
\end{axis}
\end{tikzpicture}

%% file: hour6840plus24_consumer_wgws.tex
\pgfplotsset{compat=1.11,
/pgfplots/ybar legend/.style={
/pgfplots/legend image code/.code={%
\draw[##1,/tikz/.cd,yshift=-0.25em]
(0cm,0cm) rectangle (3pt,0.8em);},
   },}
\begin{tikzpicture}[scale=0.54]
\begin{axis}[ybar stacked,
bar width=4.166666666666667pt,
xlabel=hour of day,ylabel=Power (MW),xmin=0.1,xmax=25,ymin=0,ymax=30000,
legend style={legend columns=5,at={(1.2,1.23)},},
]
\addplot +[red] coordinates {
(1,6169.999999999999)
(2,6194.000000000003)
(3,6233.999999999999)
(4,6094.0)
(5,5784.000000000002)
(6,7481.0)
(7,9515.000000000005)
(8,14793.0)
(9,20373.999999999996)
(10,22080.000000000004)
(11,23784.000000000007)
(12,24098.0)
(13,23991.00000000002)
(14,19824.999999999996)
(15,15637.999999999989)
(16,10451.999999999996)
(17,4725.999999999999)
(18,4468.000000000001)
(19,5333.999999999999)
(20,5569.999999999998)
(21,5375.000000000001)
(22,5764.0)
(23,6037.999999999996)
(24,6280.999999999998)
};
\addplot +[cyan] coordinates {
(1,30.623030892621344)
(2,30.67880072299677)
(3,30.69346660717257)
(4,30.74183805014826)
(5,30.806349984147296)
(6,30.753516946085774)
(7,30.873944177154655)
(8,151.62925885753938)
(9,2321.8494197890623)
(10,2309.5287651950894)
(11,1746.1727570720213)
(12,1878.4400915261983)
(13,1672.9668937747429)
(14,1449.5915463155268)
(15,31.613673340091133)
(16,30.595370476562806)
(17,388.51276429527377)
(18,388.49634482499374)
(19,42.810370966977914)
(20,31.46842745559034)
(21,30.94367692489555)
(22,30.84646372753915)
(23,30.920899874371777)
(24,31.019618201555502)
};
\addplot +[green] coordinates {
(1,5.78469830169201)
(2,6.33473621868139)
(3,7.08044437593582)
(4,8.27600704136041)
(5,10.6697060914579)
(6,14.296438055827)
(7,20.5302308645885)
(8,1.0184710489186)
(9,0.00021098820297101)
(10,0.000213204612846406)
(11,0.000212725349106596)
(12,0.000211576371148749)
(13,0.000216788829675435)
(14,0.00431991336378873)
(15,3.21787603170435)
(16,1.79318507680011)
(17,99.9617041539478)
(18,99.9617160150207)
(19,1.17942520389751)
(20,1.24572849130414)
(21,1.25223277073599)
(22,1.24496683240936)
(23,1.23549889477653)
(24,1.22536413432749)
};
\addplot +[black] coordinates {
(1,425.317056431309)
(2,425.605117127969)
(3,425.024407813467)
(4,427.432465111473)
(5,431.350015297699)
(6,415.935734012125)
(7,390.408159474711)
(8,296.240747622256)
(9,229.887655667482)
(10,230.111419273204)
(11,288.968861277039)
(12,93.1957790503643)
(13,106.216085760507)
(14,115.725893544049)
(15,323.59057412091)
(16,397.813265521955)
(17,1941.28653146848)
(18,1940.96182021477)
(19,575.037413931854)
(20,444.579218075136)
(21,433.474377857128)
(22,428.357646304802)
(23,426.828773461311)
(24,426.179625403346)
};
\addplot +[gray] coordinates {
(1,2147.898563365427)
(2,1882.086486310712)
(3,1684.361126297164)
(4,1895.295408871861)
(5,2541.656503442231)
(6,0.0)
(7,0.0)
(8,0.0)
(9,0.0)
(10,0.0)
(11,0.0)
(12,0.0)
(13,0.0)
(14,0.0)
(15,0.0)
(16,0.0)
(17,4757.93077341105)
(18,4442.739230948124)
(19,4613.297781285752)
(20,3464.340938059013)
(21,2303.888439733265)
(22,1369.833775367192)
(23,275.6044715396172)
(24,0.0)
};
\legend{Demand,Pump,BattTSO,LoadShift\_Up,Export}
\end{axis}
\end{tikzpicture}

%% file: hour6840plus24_consumer_ngns.tex
\pgfplotsset{compat=1.11,
/pgfplots/ybar legend/.style={
/pgfplots/legend image code/.code={%
\draw[##1,/tikz/.cd,yshift=-0.25em]
(0cm,0cm) rectangle (3pt,0.8em);},
   },}
\begin{tikzpicture}[scale=0.54]
\begin{axis}[ybar stacked,
bar width=4.166666666666667pt,
xlabel=hour of day,xmin=0.1,xmax=25,ymin=0,ymax=30000,
legend style={legend columns=4,at={(1.0,-0.1)},},
]
\addplot +[red] coordinates {
(1,6169.999999999999)
(2,6194.000000000003)
(3,6233.999999999999)
(4,6094.0)
(5,5784.000000000002)
(6,7481.0)
(7,9515.000000000005)
(8,14793.0)
(9,20373.999999999996)
(10,22080.000000000004)
(11,23784.000000000007)
(12,24098.0)
(13,23991.00000000002)
(14,19824.999999999996)
(15,15637.999999999989)
(16,10451.999999999996)
(17,4725.999999999999)
(18,4468.000000000001)
(19,5333.999999999999)
(20,5569.999999999998)
(21,5375.000000000001)
(22,5764.0)
(23,6037.999999999996)
(24,6280.999999999998)
};
\addplot +[cyan] coordinates {
(1,138.79999742385178)
(2,138.79999822498905)
(3,138.79999713680576)
(4,389.19999433380957)
(5,389.1999939005373)
(6,24.799876007210493)
(7,24.79987601387043)
(8,1519.300000979423)
(9,2309.30000019578)
(10,2074.2882359185237)
(11,1204.9413746593548)
(12,1218.4239400436527)
(13,1093.8711767558882)
(14,2309.300000196865)
(15,24.79987600344494)
(16,24.79987600521669)
(17,625.6171042320217)
(18,641.0115569225046)
(19,628.3472819617734)
(20,610.3131105673488)
(21,611.4658661169101)
(22,843.8300430912811)
(23,866.4287264873752)
(24,848.2808065693887)
};
\addplot +[green] coordinates {
(1,2.5580175532107733e-07)
(2,2.5514148845243793e-07)
(3,2.547940667449767e-07)
(4,199.9999994097809)
(5,199.99999938113183)
(6,3.440553501051944e-07)
(7,9.396373134640947e-07)
(8,3.238711332150653e-08)
(9,4.185766040322085e-10)
(10,8.327094267276444e-10)
(11,4.3781602861880994e-10)
(12,8.183327582424049e-10)
(13,17.463651858529794)
(14,4.951673091306926e-10)
(15,9.272076567121837e-07)
(16,3.365975064752724e-07)
(17,140.98817831211807)
(18,118.8918762578614)
(19,102.61614806978396)
(20,93.7229222868443)
(21,94.57507777505676)
(22,23.14420812211706)
(23,22.34833780921786)
(24,23.80807809710858)
};
\addplot +[gray] coordinates {
(1,4261.757469227155)
(2,4560.916190813688)
(3,4186.750159818004)
(4,4434.822137447383)
(5,4387.725140833155)
(6,0.0)
(7,0.0)
(8,596.9607243868168)
(9,0.0)
(10,0.0)
(11,0.0)
(12,538.1347759626337)
(13,508.2067899695986)
(14,0.0)
(15,0.0)
(16,0.0)
(17,3437.042007449756)
(18,3375.492336883408)
(19,2805.306730017577)
(20,3088.226617133258)
(21,3756.539165412474)
(22,3807.206777256897)
(23,527.0885724380406)
(24,0.0)
};
%\legend{Demand,Pump,BattTSO,Export}
\end{axis}
\end{tikzpicture}

%% file: dam_levels.tex
\begin{tikzpicture}[scale=0.8]
\begin{axis}[
  legend style={
    legend columns=4,
    at={(1.0,+1.3)},
  },
  x label style={at={(0.55,-0.099)}},
  y label style={at={(-0.04,0.56)}},
  ylabel={Storage (TWh)},
  xlabel={month},
  xtick={1,2,3,4,5,6,7,8,9,10,11,12},
  ymin=0,
  xmin=0.5,
  xmax=12.5,
  ymax=6,
  height=5cm,
  width=8.3cm,
  mark size=0.0pt,
]

% NGNS --------
\addplot +[thick, sharp plot] coordinates {(1,2.55) (2,0.33) (3,0.03) (4,0.15) (5,0.32) (6,1.86) (7,2.59) (8,3.81) (9,5.07) (10,5.44) (11,4.34) (12,3.47) };
% WGWS
\addplot +[dashed, thick, sharp plot] coordinates {(1,2.7) (2,0.39) (3,0.32) (4,0.277) (5,1.26) (6,2.86) (7,4.25) (8,5.25) (9,5.65) (10,5.49) (11,4.34) (12,3.47) };

% WGNS
\addplot +[dotted, thick, sharp plot] coordinates {(1,2.59) (2,0.35) (3,0.006) (4,0.21) (5,0.32) (6,1.93) (7,2.71) (8,4.12) (9,5.43) (10,5.80) (11,4.44) (12,3.47) };

% NGWS
\addplot +[dashdotted, thick, sharp plot] coordinates {(1,2.685) (2,0.33) (3,0.287) (4,0.292) (5,0.989) (6,2.86) (7,4.38) (85.38) (9,5.82) (10,5.8) (11,4.57) (12,3.47) };

%\legend{Historical,Simulated}
\legend{"NGNS", "WGWS", "WGNS", "NGWS"}
\end{axis}
\end{tikzpicture}

%% file: week_producer_red_wgws.tex
\pgfplotsset{compat=1.11,
/pgfplots/ybar legend/.style={
/pgfplots/legend image code/.code={%
\draw[##1,/tikz/.cd,yshift=-0.25em]
(0cm,0cm) rectangle (3pt,0.8em);},
   },}
\begin{tikzpicture}[scale=0.55]
\begin{axis}[ybar stacked,
bar width=1.9230769230769231pt,
xlabel=week of year,ylabel=Production (TWh),xmin=0.1,xmax=53,ymin=0,ymax=3.5,
legend style={legend columns=6,at={(1.45,1.3)},},
]
\addplot +[orange] coordinates {
(1,0.07845599999980382)
(2,0.07870799999980499)
(3,0.07870799999980588)
(4,0.07870799999980656)
(5,0.07870799999980703)
(6,0.07870799999980477)
(7,0.07870799999980417)
(8,0.07870799999978581)
(9,0.07870169999977138)
(10,0.07869959999982812)
(11,0.07868229999988528)
(12,0.07870799999971569)
(13,0.07855520000146743)
(14,0.07848849999838822)
(15,0.0754716000016911)
(16,0.07864679999789069)
(17,0.07842497034499234)
(18,0.07828340393968017)
(19,0.0785114000027723)
(20,0.07800794800358254)
(21,0.07802970000484798)
(22,0.07807177925700005)
(23,0.07847889999700108)
(24,0.07816453773241615)
(25,0.07820279113153439)
(26,0.07831444800609255)
(27,0.0782449000006256)
(28,0.0784698000030261)
(29,0.078276828810146)
(30,0.07772303815545302)
(31,0.07829672827241749)
(32,0.0782609000169325)
(33,0.07821980000502876)
(34,0.07856580000029241)
(35,0.07863198587762879)
(36,0.07868549999759279)
(37,0.07856659999751914)
(38,0.07827690816090999)
(39,0.0787079999975884)
(40,0.07868759999865395)
(41,0.07844382678738783)
(42,0.0787079999988235)
(43,0.07870799999872553)
(44,0.07870799999931974)
(45,0.07870799999980843)
(46,0.07870799999981055)
(47,0.07870799999980606)
(48,0.07870799999979795)
(49,0.07870799999980704)
(50,0.07870799999980661)
(51,0.07870799999980643)
(52,0.07811569999972869)
};
\addplot +[blue] coordinates {
(1,0.23402199239891733)
(2,0.23402199239891022)
(3,0.23402199239891266)
(4,0.234021992398917)
(5,0.2544413635339196)
(6,0.2697661965589083)
(7,0.2697661965589046)
(8,0.269766196558893)
(9,0.25997706818982425)
(10,0.2526909029987024)
(11,0.25253007480329914)
(12,0.2526998039984124)
(13,0.27158350391777925)
(14,0.3872351357205257)
(15,0.379346718973197)
(16,0.38555732214281435)
(17,0.3836852245108679)
(18,0.46164073292302676)
(19,0.4734724447198445)
(20,0.46910109545825784)
(21,0.46874485423683054)
(22,0.5219847983815668)
(23,0.597586939284583)
(24,0.5899118068451578)
(25,0.5852611907514556)
(26,0.5996846178060048)
(27,0.6015508164049936)
(28,0.6107237046859663)
(29,0.6156149067349582)
(30,0.5935090001368237)
(31,0.5571171184887845)
(32,0.5300327517572629)
(33,0.5308141851891129)
(34,0.5354461322349452)
(35,0.49479148076502627)
(36,0.3891729471753882)
(37,0.3814345725519462)
(38,0.3865707471863344)
(39,0.3881539886032278)
(40,0.237413976118721)
(41,0.2366106322252537)
(42,0.23786135888801663)
(43,0.23786575018344822)
(44,0.2534741222354966)
(45,0.265133964208952)
(46,0.26512028575893587)
(47,0.26515167839891324)
(48,0.24269644943886662)
(49,0.18655837703891795)
(50,0.18655837703891615)
(51,0.1865583770389145)
(52,0.18609345529992327)
};
\addplot +[yellow] coordinates {
(1,0.10809678406008809)
(2,0.15671330869028832)
(3,0.1833132874842989)
(4,0.1339839466937172)
(5,0.1750981054149153)
(6,0.27560547885745657)
(7,0.28487522486723066)
(8,0.3502663188010696)
(9,0.4687107324735895)
(10,0.44658908568012207)
(11,0.5767815617289471)
(12,0.5719777127424045)
(13,0.7001810430070197)
(14,0.4727289341796072)
(15,0.4494945844370933)
(16,0.5211841894968)
(17,0.7113001467822374)
(18,0.7097257501518212)
(19,0.7015261827657656)
(20,0.79491585634522)
(21,0.7670055102756487)
(22,0.8400766506936491)
(23,0.5200284683401138)
(24,0.7529745314801592)
(25,0.8569608855360217)
(26,0.8265816354854782)
(27,0.5600279102200826)
(28,0.7418373610144207)
(29,0.7478293094664771)
(30,0.8396855027503577)
(31,0.7835803714704769)
(32,0.7822945136066394)
(33,0.8141788198935593)
(34,0.6937385086857618)
(35,0.4969056813367004)
(36,0.5538584855976008)
(37,0.5657694865643971)
(38,0.52232945711303)
(39,0.3903976321567081)
(40,0.3841875962334478)
(41,0.2268103374008557)
(42,0.36813406165120993)
(43,0.22464017696109687)
(44,0.29440964394528407)
(45,0.22117586380289211)
(46,0.18493626991962459)
(47,0.18132908161862216)
(48,0.09942642115571684)
(49,0.10104980497585254)
(50,0.14376308356526915)
(51,0.1141417094959069)
(52,0.14086508999664082)
};
\addplot +[lime] coordinates {
(1,0.007432179999947695)
(2,0.008637649999946947)
(3,0.006010259999947304)
(4,0.005510879999947448)
(5,0.005819579999947717)
(6,0.008500229999947073)
(7,0.00865111999994693)
(8,0.00614667999994616)
(9,0.006752609999938829)
(10,0.008186059999934927)
(11,0.0073979399999115234)
(12,0.005408379999920916)
(13,0.00836866000001343)
(14,8.447251561529864e-13)
(15,6.770000791942811e-06)
(16,0.00029620952104791485)
(17,1.8772428932825137e-05)
(18,5.35284325691178e-05)
(19,3.0944524605472006e-05)
(20,5.699475621116535e-13)
(21,2.2754416603100246e-05)
(22,5.756715404568082e-13)
(23,5.795817459338811e-13)
(24,5.654773836452724e-13)
(25,6.261000057191781e-05)
(26,5.815331183112243e-13)
(27,9.131000054521667e-05)
(28,5.605000057686722e-05)
(29,5.774700155767377e-13)
(30,0.00012356000055104276)
(31,5.777989953021384e-13)
(32,5.686293254718202e-13)
(33,5.757033356725121e-13)
(34,5.838545442549335e-13)
(35,5.854337682865807e-13)
(36,1.738996060076962e-05)
(37,0.0005588500005282752)
(38,0.00016248440079189526)
(39,2.3520660812857137e-05)
(40,0.004444389995372952)
(41,0.0043327120737203665)
(42,0.004650529997601543)
(43,0.006015809995155507)
(44,0.005609329998198593)
(45,0.004436079999947804)
(46,0.0025223399999482776)
(47,0.003101119999947025)
(48,0.007686389999944644)
(49,0.00790191999994772)
(50,0.005963669999947496)
(51,0.004202139999947314)
(52,0.0052783599999468075)
};
\addplot +[teal] coordinates {
(1,0.5879912110778994)
(2,0.4732314500144001)
(3,0.47780669752020655)
(4,0.5632231328485829)
(5,0.5276288500433657)
(6,0.361222840628111)
(7,0.3838125768634419)
(8,0.49193142051701544)
(9,0.229875020083226)
(10,0.14437419999672918)
(11,0.0608503076243942)
(12,0.21350109004221163)
(13,0.1862179963522295)
(14,0.05208289568022168)
(15,0.1723644049856393)
(16,0.31940649043554037)
(17,0.2715652618124639)
(18,0.3197376508993475)
(19,0.41774844447485104)
(20,0.28579550094276507)
(21,0.3768580187138175)
(22,0.35390519947386573)
(23,0.28511664261418734)
(24,0.21813695736338454)
(25,0.3585311906400483)
(26,0.30901266071606454)
(27,0.28130698509448504)
(28,0.2541726527346745)
(29,0.21835895369250238)
(30,0.25234071263794833)
(31,0.27584292824905193)
(32,0.276777932166433)
(33,0.2778715942033066)
(34,0.35159001271039536)
(35,0.3695185265115229)
(36,0.4399117165946686)
(37,0.44754312178032774)
(38,0.45313029524604315)
(39,0.3983349107782502)
(40,0.25770720700122157)
(41,0.3848391673397913)
(42,0.27695002313353445)
(43,0.14338186485859694)
(44,0.3682910329021378)
(45,0.6253180350562609)
(46,0.553851738775038)
(47,0.34485543795810253)
(48,0.36324890865984355)
(49,0.6337067929289477)
(50,0.5745222544543688)
(51,0.49946832426969295)
(52,0.5288517633053152)
};
\addplot +[cyan] coordinates {
(1,0.12385732080574206)
(2,0.09077137962216417)
(3,0.07798426801073802)
(4,0.14450352902229188)
(5,0.12041427539705261)
(6,0.08569232307434907)
(7,0.0901730150339919)
(8,0.11435020984219652)
(9,0.1285388547915318)
(10,0.07761871641341231)
(11,0.03436682413230314)
(12,0.1348354970853686)
(13,0.07951743248367728)
(14,0.002397113449818101)
(15,0.019415732862515078)
(16,0.04083681268481758)
(17,0.019717224576679397)
(18,0.012554890162230454)
(19,0.011011813834064634)
(20,0.0046282114474718505)
(21,0.010623308082936147)
(22,0.0033705327824925992)
(23,0.014360334187039785)
(24,0.00687480182718041)
(25,0.03853230297461802)
(26,0.007517057638803196)
(27,0.05397788905512828)
(28,0.019314729690261358)
(29,0.005970690940626624)
(30,0.05285566514178267)
(31,0.0030634015381280954)
(32,5.500005408068879e-05)
(33,0.00011000007368574955)
(34,5.500009636386556e-05)
(35,6.670892128225845e-05)
(36,0.0142566259577778)
(37,0.0298333707244761)
(38,0.02258969800291154)
(39,0.009026586748092534)
(40,0.04662260691715835)
(41,0.10857610299130159)
(42,0.013102514807958952)
(43,0.015660395672125633)
(44,0.08382289433369366)
(45,0.16400955916521598)
(46,0.16173734587846197)
(47,0.08594259490972667)
(48,0.08793287275484285)
(49,0.12623039757054427)
(50,0.11749049726392508)
(51,0.12346294321216146)
(52,0.11668576970922)
};
\addplot +[green] coordinates {
(1,0.0001000000000562647)
(2,0.0007683199984110168)
(3,0.0007297042651987515)
(4,0.001016472020633133)
(5,0.0006093231129785688)
(6,0.0009447999946239291)
(7,0.0008983119075046942)
(8,0.0009602480984576326)
(9,0.002341627257646336)
(10,0.0021957983253416613)
(11,0.002120000000744257)
(12,0.0022409880153605225)
(13,0.0023564878303994354)
(14,0.004579905829593133)
(15,0.0031491200002424224)
(16,0.0032689710907748046)
(17,0.0028061466441152735)
(18,0.0020786182975831725)
(19,0.00211108116516296)
(20,0.0017209320779491494)
(21,0.002058680377578918)
(22,0.0023040000013463606)
(23,0.0014711157867517248)
(24,0.0021121599516819665)
(25,0.0027800218866734757)
(26,0.0021623448426518034)
(27,0.001929133468758176)
(28,0.00167747101666993)
(29,0.0024729599940959716)
(30,0.0023040000092854423)
(31,0.0022803318090587332)
(32,0.002217707062506911)
(33,0.0023877301653942896)
(34,0.0021263810901566)
(35,0.0008832622444711508)
(36,0.000981034864914909)
(37,0.002856860340425283)
(38,0.002353705314388214)
(39,0.0022139005627622977)
(40,0.0022071732235705264)
(41,0.0023040000003604964)
(42,0.0009447999986404572)
(43,0.0008601600011835994)
(44,0.0010444800002148528)
(45,0.0005340013699875294)
(46,0.0004930517194170789)
(47,0.0005529600000779117)
(48,0.0006604793542676398)
(49,0.00020000065515881135)
(50,0.00018399997998095412)
(51,0.0011212800156869638)
(52,0.0010291199947943393)
};
\addplot +[brown] coordinates {
(1,0.27836801999979277)
(2,0.2793735379906288)
(3,0.27971999999979474)
(4,0.2808091031649285)
(5,0.28118406419437303)
(6,0.27679057480844926)
(7,0.27837644810295614)
(8,0.27740799999979576)
(9,0.2562358051818411)
(10,0.24617473308374987)
(11,0.228026303365377)
(12,0.27400234971434884)
(13,0.21733471489051018)
(14,2.6473344112191915e-13)
(15,2.6084924968812955e-13)
(16,0.001709699789970363)
(17,3.500359331925584e-05)
(18,0.0026184690935800075)
(19,0.0004114238891346546)
(20,2.557587805020204e-13)
(21,0.0010533953496539043)
(22,2.558994966601757e-13)
(23,2.5652759973565273e-13)
(24,2.556624661598087e-13)
(25,0.0008404798524750806)
(26,2.563878444321197e-13)
(27,0.0015383111512297064)
(28,2.580554153957813e-13)
(29,2.5642405252410397e-13)
(30,2.5784910251039945e-13)
(31,0.00013054456535405696)
(32,2.559706027026236e-13)
(33,2.5612625632082663e-13)
(34,2.564698562367335e-13)
(35,2.56430541779281e-13)
(36,0.0001787565225979399)
(37,2.5995059946330565e-13)
(38,0.005751838420665653)
(39,2.6466849994842396e-13)
(40,0.011453957622738055)
(41,0.030040408834546353)
(42,0.0009407416258552655)
(43,0.002885772697412726)
(44,0.15880002000003207)
(45,0.28189323594763727)
(46,0.28299609620032706)
(47,0.27971593999978966)
(48,0.25694412893100327)
(49,0.2798939600651447)
(50,0.2797199999998094)
(51,0.27971999999980707)
(52,0.27355378000025504)
};
\addplot +[pink] coordinates {
(1,1.8277993529380096e-14)
(2,1.8260638621407727e-14)
(3,1.826325454916006e-14)
(4,1.8283906469821002e-14)
(5,2.8128932819170647e-06)
(6,1.8249881469140484e-14)
(7,1.8264194805057344e-14)
(8,1.8269552513382494e-14)
(9,1.8206589260090414e-14)
(10,1.8151358939603413e-14)
(11,1.4255614117847441e-05)
(12,1.7988805271513833e-14)
(13,1.7937016004826107e-14)
(14,1.7703110921487477e-14)
(15,1.7698639376198093e-14)
(16,1.7727498676715454e-14)
(17,1.7706893643831887e-14)
(18,2.557198102848467e-08)
(19,1.769869558069038e-14)
(20,1.7685072222001547e-14)
(21,1.769314184845535e-14)
(22,1.7693068258898016e-14)
(23,1.7693136675741424e-14)
(24,1.7685046753740682e-14)
(25,1.770370976185043e-14)
(26,1.7697478954341768e-14)
(27,1.7695726383865758e-14)
(28,1.7694492724728837e-14)
(29,1.7688177820691107e-14)
(30,1.7688592909798093e-14)
(31,1.769278233236993e-14)
(32,1.768623486454951e-14)
(33,1.7687360869319688e-14)
(34,1.7692796879331257e-14)
(35,1.7693002140443927e-14)
(36,1.7699799502939963e-14)
(37,1.7703411232367248e-14)
(38,1.7728943513999465e-14)
(39,1.7714951997737338e-14)
(40,0.0006609013998299777)
(41,0.0014358808658988928)
(42,1.7753849100754263e-14)
(43,1.7745647544977713e-14)
(44,1.804132176841626e-14)
(45,3.229813127985037e-05)
(46,0.0012959909651013074)
(47,1.8245419074704556e-14)
(48,1.823818812611023e-14)
(49,1.827838018482785e-14)
(50,1.8274568241017304e-14)
(51,1.8261866289997797e-14)
(52,1.826175652915866e-14)
};
\addplot +[black] coordinates {
(1,0.09153503445132728)
(2,0.09449239378539841)
(3,0.09394632060294378)
(4,0.09594464084582555)
(5,0.09617957002499182)
(6,0.0964676947151207)
(7,0.0961001739707743)
(8,0.09629744503474127)
(9,0.07712275286211547)
(10,0.0756015237976997)
(11,0.08180477748206787)
(12,0.0836633931100399)
(13,0.08573293154172941)
(14,0.06645044108539838)
(15,0.06912308158924779)
(16,0.07752605281593129)
(17,0.07039290605836099)
(18,0.07769415908430205)
(19,0.07907425780176659)
(20,0.07368116735524316)
(21,0.07156656859834182)
(22,0.07175595668332461)
(23,0.07018709553913728)
(24,0.06331103523082081)
(25,0.07205227395408838)
(26,0.07324319624599213)
(27,0.06738276047592412)
(28,0.07070086821294447)
(29,0.06305620499188122)
(30,0.06401984912418078)
(31,0.07036891294281741)
(32,0.06464348381605964)
(33,0.0668448042158848)
(34,0.06598743118154465)
(35,0.06071555747187908)
(36,0.0656346966528001)
(37,0.0621482240895556)
(38,0.06368408292135075)
(39,0.0712781054065493)
(40,0.07253965272221137)
(41,0.08261915330832471)
(42,0.08078081848499004)
(43,0.07591287636366857)
(44,0.07883281634650417)
(45,0.09523593364628127)
(46,0.09597374137191482)
(47,0.06933997176668702)
(48,0.05836416840067981)
(49,0.09291288539605261)
(50,0.09523708615576396)
(51,0.09342507701418062)
(52,0.09653215008109059)
};
\addplot +[lightgray] coordinates {
(1,0.4302319765892801)
(2,0.4707781255333982)
(3,0.4365649062609671)
(4,0.4912419696897828)
(5,0.4821305555320145)
(6,0.44223652519765855)
(7,0.4486765650847248)
(8,0.4456877963081268)
(9,0.44398842570485925)
(10,0.4816936636618596)
(11,0.5027529545541302)
(12,0.5127183949096779)
(13,0.4378227839096269)
(14,0.532033889091113)
(15,0.4748774875290408)
(16,0.6013089850539544)
(17,0.4418171600043691)
(18,0.35451694848930193)
(19,0.35877770254663344)
(20,0.27815211214867747)
(21,0.28455122003402283)
(22,0.2793346305100334)
(23,0.36243724598976385)
(24,0.2869830832350332)
(25,0.3265055026370903)
(26,0.3207832785201759)
(27,0.35159058242374275)
(28,0.34860899725877337)
(29,0.28500705882178146)
(30,0.26812466513940036)
(31,0.35706429040120014)
(32,0.31214189106386986)
(33,0.2970197906370679)
(34,0.3052369344592915)
(35,0.3585330938657823)
(36,0.44963104579705176)
(37,0.4421130503111232)
(38,0.42882167445241137)
(39,0.5561151195458736)
(40,0.6527437027056974)
(41,0.6890927927164704)
(42,0.6660111699437248)
(43,0.7499318702114487)
(44,0.5585941280983586)
(45,0.4540428609442103)
(46,0.4619523470987525)
(47,0.42491601406844876)
(48,0.46298439090395227)
(49,0.4481520497931748)
(50,0.45447858917071277)
(51,0.4384614909426663)
(52,0.47828701286078923)
};
\legend{Biomass,RoR,PV,WindOn,Dam,Pump,BattTSO,Gas,LoadShed,LoadShift\_Down,Import}
\end{axis}
\end{tikzpicture}

%% file: week_producer_red_ngws.tex
\pgfplotsset{compat=1.11,
/pgfplots/ybar legend/.style={
/pgfplots/legend image code/.code={%
\draw[##1,/tikz/.cd,yshift=-0.25em]
(0cm,0cm) rectangle (3pt,0.8em);},
   },}
\begin{tikzpicture}[scale=0.55]
\begin{axis}[ybar stacked,
bar width=1.9230769230769231pt,
xlabel=week of year,xmin=0.1,xmax=53,ymin=0,ymax=3.5,
legend style={legend columns=4,at={(1.0,-0.1)},},
]
\addplot +[orange] coordinates {
(1,0.07845599999791272)
(2,0.0787079999978979)
(3,0.0787079999974694)
(4,0.07866629999809083)
(5,0.0787079999977661)
(6,0.0787079999975491)
(7,0.07870799999639791)
(8,0.0787079999970502)
(9,0.07866599996102396)
(10,0.07867091087472812)
(11,0.07862659996634914)
(12,0.07870799997425468)
(13,0.07810904847017273)
(14,0.07814290001334535)
(15,0.07483259999960302)
(16,0.07854580000902252)
(17,0.07790354027682153)
(18,0.07774858921586271)
(19,0.07778820007504356)
(20,0.07702064810462154)
(21,0.07707799533037274)
(22,0.07679700012622549)
(23,0.07806640001748814)
(24,0.07713258578478653)
(25,0.07749698210422604)
(26,0.0776263001037188)
(27,0.0780182000092963)
(28,0.07778741872452528)
(29,0.07758608137140835)
(30,0.07573380014150273)
(31,0.07745323001423497)
(32,0.07720830014615061)
(33,0.0774309723636189)
(34,0.07800660008303188)
(35,0.07829110004345478)
(36,0.0784363000166705)
(37,0.07817247973512817)
(38,0.07770271658376067)
(39,0.07849220001885769)
(40,0.07802980002299961)
(41,0.07761609998421354)
(42,0.07870369996735585)
(43,0.07863230000154706)
(44,0.07866239994869903)
(45,0.07870799999815382)
(46,0.07870799999822999)
(47,0.07870799999813671)
(48,0.07870799999813947)
(49,0.07870799999820001)
(50,0.07870799999819555)
(51,0.07870799999819601)
(52,0.07797649999872032)
};
\addplot +[blue] coordinates {
(1,0.23402199238883145)
(2,0.2340219923886301)
(3,0.2340219923886082)
(4,0.2338366887787862)
(5,0.254375798928822)
(6,0.26976619654850936)
(7,0.2693539345465627)
(8,0.26930165671000467)
(9,0.2523086028169346)
(10,0.250447492595743)
(11,0.24675743062083316)
(12,0.2526998039374359)
(13,0.25948048253152434)
(14,0.3794249529769969)
(15,0.3748403046791939)
(16,0.387299323636986)
(17,0.3754239447066581)
(18,0.44955194343997396)
(19,0.4449891469418606)
(20,0.40949803352687764)
(21,0.4051390925564262)
(22,0.45904952006936245)
(23,0.5616428618335736)
(24,0.49328154307620337)
(25,0.5528370015191778)
(26,0.5687269127548849)
(27,0.5720925385053189)
(28,0.5754819754930094)
(29,0.5509546797329975)
(30,0.48819966473022186)
(31,0.5063214556067792)
(32,0.4574155161521601)
(33,0.4742248793777697)
(34,0.5170430277112578)
(35,0.4907357936123063)
(36,0.38765808828650494)
(37,0.381188871709256)
(38,0.3796468434154652)
(39,0.3890554483800519)
(40,0.237153522243634)
(41,0.2357842664762302)
(42,0.23790404719015187)
(43,0.23774040600989912)
(44,0.25326778317739695)
(45,0.26515167839009257)
(46,0.2648860898503695)
(47,0.2651516783900213)
(48,0.2426964494297652)
(49,0.1865472135290861)
(50,0.18655837703007022)
(51,0.1865583770300694)
(52,0.18618178270402924)
};
\addplot +[yellow] coordinates {
(1,0.16449253686464133)
(2,0.25623200037541694)
(3,0.2847264637967658)
(4,0.20761390940255905)
(5,0.2825925166374279)
(6,0.4455444693358797)
(7,0.4427222873296974)
(8,0.5339096968737682)
(9,0.6956143366249075)
(10,0.6770890875665543)
(11,0.8503952419717218)
(12,0.839088224546198)
(13,1.02953160110288)
(14,0.6984632747653615)
(15,0.6649881736672767)
(16,0.7917565185881285)
(17,1.0376520749799665)
(18,1.0291716324874154)
(19,1.0200068256124368)
(20,1.1275460861733024)
(21,1.123737761308322)
(22,1.2267116814451362)
(23,0.7727994411722154)
(24,1.0927469123399418)
(25,1.2500586644008902)
(26,1.1973104825941108)
(27,0.8320404429752911)
(28,1.0859615856273297)
(29,1.0917844684089344)
(30,1.2285813689698284)
(31,1.1448975139320892)
(32,1.1257331436622202)
(33,1.1733638277864002)
(34,1.0113195116737423)
(35,0.728906461570644)
(36,0.8334197205681587)
(37,0.8346238071307119)
(38,0.7623324756040232)
(39,0.5846361145431073)
(40,0.5860826781500564)
(41,0.3518686156142334)
(42,0.5538685724778627)
(43,0.3592643549075889)
(44,0.4489093165944996)
(45,0.33462494278724625)
(46,0.29449695333119114)
(47,0.286888166705903)
(48,0.15143131800515497)
(49,0.15636182139697713)
(50,0.22135331701948596)
(51,0.1705575406425898)
(52,0.2147846780241245)
};
\addplot +[lime] coordinates {
(1,0.09176457999904913)
(2,0.10664474999903104)
(3,0.07420665999903098)
(4,0.06803567999905445)
(5,0.07185097999905526)
(6,0.10494832999902265)
(7,0.10578932915649525)
(8,0.07582539999895782)
(9,0.07902791999904103)
(10,0.09192074340364954)
(11,0.08374582382087703)
(12,0.06520487792862976)
(13,0.08289095321469311)
(14,1.7354686176494562e-11)
(15,6.770015905281975e-06)
(16,1.7637363564229053e-11)
(17,1.680330416989927e-11)
(18,0.00019478329587422326)
(19,1.5751053501675766e-11)
(20,1.5040472425766827e-11)
(21,1.528100430204701e-11)
(22,1.520797207192199e-11)
(23,1.605603376626143e-11)
(24,1.4868663514919185e-11)
(25,0.0022363800154223184)
(26,1.5220725314774598e-11)
(27,0.0008234500165176875)
(28,0.0020670500155712287)
(29,1.5008403615401545e-11)
(30,0.004554560013781159)
(31,1.5549975500164418e-11)
(32,1.524802798122703e-11)
(33,1.529853927094739e-11)
(34,1.595365256846272e-11)
(35,1.6474858998131263e-11)
(36,1.7015743104826907e-11)
(37,1.875782861215846e-11)
(38,0.001985996745317655)
(39,1.8358408357093826e-11)
(40,0.007838308730996618)
(41,0.045374194876270875)
(42,0.04327869057485114)
(43,0.008959207049403319)
(44,0.038860970018735216)
(45,0.05476657999930391)
(46,0.031138339999317345)
(47,0.03828631999929774)
(48,0.09490548999927355)
(49,0.09755861999930195)
(50,0.07363726999930122)
(51,0.05188503999930213)
(52,0.06473285999929651)
};
\addplot +[teal] coordinates {
(1,0.6910763739814264)
(2,0.49436226716120557)
(3,0.5710833908644084)
(4,0.6765303027301784)
(5,0.5545968002939976)
(6,0.2793947675915396)
(7,0.3341289302423803)
(8,0.4081887594211277)
(9,0.23254446869949538)
(10,0.1217797389658426)
(11,0.043875840570338526)
(12,0.17579526838185916)
(13,0.11830988896316813)
(14,0.1748030096133764)
(15,0.15168053269995804)
(16,0.22969034924841034)
(17,0.25474523801609394)
(18,0.2977497306926785)
(19,0.388555985568314)
(20,0.2801195100716506)
(21,0.33497949954724066)
(22,0.2810105600509879)
(23,0.40594600967052863)
(24,0.2972290531146459)
(25,0.3506164466355551)
(26,0.32006194942218924)
(27,0.3860021508131737)
(28,0.2534679321525527)
(29,0.24704632410757227)
(30,0.2526578630738933)
(31,0.24263090357019554)
(32,0.2529119645317381)
(33,0.22262223688845842)
(34,0.3440750025323916)
(35,0.4378638004889871)
(36,0.4233188852324443)
(37,0.41907695255794414)
(38,0.40726513466140235)
(39,0.3775700503460909)
(40,0.24401859862584535)
(41,0.4890823368724215)
(42,0.2240103539644038)
(43,0.07636829818399292)
(44,0.3558023027272996)
(45,0.6395550639304176)
(46,0.6124437918273128)
(47,0.35835676106048164)
(48,0.3828101918628237)
(49,0.6277837258430338)
(50,0.5529456338251872)
(51,0.5839705953195382)
(52,0.43518592203372053)
};
\addplot +[cyan] coordinates {
(1,0.11308364721246712)
(2,0.08990355258815654)
(3,0.08175534093648253)
(4,0.14920118180429326)
(5,0.13870543191331972)
(6,0.09347495159749564)
(7,0.11508184255705065)
(8,0.13333513119678342)
(9,0.15627136150391138)
(10,0.09766286985131334)
(11,0.03374411667163535)
(12,0.1362455905987328)
(13,0.12226575772704498)
(14,0.03040377033736259)
(15,0.026852067670530964)
(16,0.013030708955329142)
(17,0.03066756795742324)
(18,0.00532455171335194)
(19,0.013595695356129277)
(20,0.009875823724333906)
(21,0.024919479406229912)
(22,0.04000665160026523)
(23,0.031246528985201815)
(24,0.01126118347969032)
(25,0.04835003719324903)
(26,0.012838239146769572)
(27,0.06046740944664314)
(28,0.032124446520805004)
(29,0.013635418490865782)
(30,0.048439551231442994)
(31,0.015533644904479428)
(32,0.012186615417141802)
(33,0.010991785888405043)
(34,0.006437287536154093)
(35,0.00857586084844379)
(36,0.010813544002971695)
(37,0.04392375773094953)
(38,0.021742499473345517)
(39,0.0006691621291557687)
(40,0.04133646416181376)
(41,0.09823332771904797)
(42,0.014785523379890878)
(43,0.006619329868201065)
(44,0.07181719745320878)
(45,0.16149329200271878)
(46,0.16173307604248377)
(47,0.09132653029157223)
(48,0.11045862323860796)
(49,0.15863312468722682)
(50,0.12399970190755354)
(51,0.14772870387943546)
(52,0.15052417774965723)
};
\addplot +[green] coordinates {
(1,0.00036720923172542315)
(2,0.0010745506098929077)
(3,0.0019017714065442902)
(4,0.0017210664909861075)
(5,0.0023294502576958202)
(6,0.0027347200158483937)
(7,0.0022347200034003208)
(8,0.002204640001997435)
(9,0.005165332894122937)
(10,0.004379259311704807)
(11,0.004763727944520375)
(12,0.005365595588714408)
(13,0.004772107326815758)
(14,0.0038508004053445168)
(15,0.004174878716740162)
(16,0.00389710317421679)
(17,0.004692960479535182)
(18,0.00528592807474711)
(19,0.005090658082202191)
(20,0.0052335021012161675)
(21,0.0049688189977066)
(22,0.00537329014246315)
(23,0.0041178274859803845)
(24,0.004179414380148903)
(25,0.005919042033530222)
(26,0.005928959479808016)
(27,0.005130056184584326)
(28,0.0052484897323601945)
(29,0.006270385540724601)
(30,0.005422656170449792)
(31,0.004452345664572603)
(32,0.004432843919894098)
(33,0.005785664559526832)
(34,0.005535775334186299)
(35,0.0031187904314810906)
(36,0.0029552293573070665)
(37,0.005907315270994692)
(38,0.004724705519424758)
(39,0.0036077991834380825)
(40,0.004094838663777514)
(41,0.0036834194780149187)
(42,0.003460447047452772)
(43,0.0032458436481343942)
(44,0.0020811654751589424)
(45,0.0007804328129449797)
(46,0.0008973327883814591)
(47,0.001292352168376737)
(48,0.0011366400881039729)
(49,0.00039666902319543195)
(50,5.685888462167046e-05)
(51,0.0005408749185646038)
(52,0.0011629622329323052)
};
\addplot +[pink] coordinates {
(1,2.659785423096492e-13)
(2,2.6671706990667135e-13)
(3,2.665714077436061e-13)
(4,0.00041349707740597733)
(5,0.0005237628957216847)
(6,2.659600910756363e-13)
(7,2.656019981823065e-13)
(8,2.663167293761955e-13)
(9,2.673023539749318e-13)
(10,2.6834789218055435e-13)
(11,2.6334389203199004e-13)
(12,2.640936259865549e-13)
(13,2.6350273871148623e-13)
(14,2.610331853494994e-13)
(15,2.6072510825805146e-13)
(16,2.6165335078552685e-13)
(17,2.610415819799364e-13)
(18,0.00019533948399357644)
(19,2.60841058721953e-13)
(20,2.5978215417829083e-13)
(21,2.6015676669495517e-13)
(22,2.6012364343442604e-13)
(23,2.602398254613574e-13)
(24,2.5955172466355643e-13)
(25,2.6050352948753375e-13)
(26,2.601090630651628e-13)
(27,2.602122045552586e-13)
(28,2.601741897582351e-13)
(29,2.5974674284743863e-13)
(30,2.597514786340478e-13)
(31,2.6019676443242e-13)
(32,2.5996144545335114e-13)
(33,2.6006504443633245e-13)
(34,2.605196432442183e-13)
(35,2.6069225178749485e-13)
(36,2.6125185883103235e-13)
(37,2.6132131129765753e-13)
(38,0.0006341892240329384)
(39,2.614464405257315e-13)
(40,0.0009534890525169347)
(41,0.0035723177854212625)
(42,2.6172913767760007e-13)
(43,3.566343138527797e-05)
(44,2.6733855896711597e-13)
(45,3.827042200121362e-05)
(46,0.002949948341412508)
(47,2.7174646230578114e-13)
(48,2.7098188451580584e-13)
(49,5.3551018970032565e-05)
(50,2.699357420851855e-13)
(51,2.701832373544965e-13)
(52,2.6966042773521247e-13)
};
\addplot +[black] coordinates {
(1,0.08916436661822286)
(2,0.09321576400752918)
(3,0.09215230375223393)
(4,0.09539727679144845)
(5,0.0959792366183214)
(6,0.09707855023962793)
(7,0.09641608539557543)
(8,0.09628973178160524)
(9,0.07498712603853902)
(10,0.06739020372281085)
(11,0.07537246568667645)
(12,0.07653530591472124)
(13,0.08455542232787025)
(14,0.06263512385723458)
(15,0.0665941079371568)
(16,0.07107067728694756)
(17,0.0668529194016637)
(18,0.07327399933087395)
(19,0.07671762237780369)
(20,0.07200960399927162)
(21,0.07038063675646374)
(22,0.0699303882921318)
(23,0.06959529299923338)
(24,0.06192833691353011)
(25,0.07048474078227206)
(26,0.06709365063768577)
(27,0.07027855876358655)
(28,0.06895193053162049)
(29,0.06175364360122507)
(30,0.060869412247832703)
(31,0.06976388752572206)
(32,0.06318884436908774)
(33,0.06429068129664502)
(34,0.06448413738315552)
(35,0.05998789137247801)
(36,0.06424264117555936)
(37,0.0620499352085794)
(38,0.06078389904597965)
(39,0.06795406178674734)
(40,0.06759345726564427)
(41,0.07553769214816781)
(42,0.06752794943030069)
(43,0.07310237010639731)
(44,0.07673596728805729)
(45,0.09099306662228558)
(46,0.09210322891741671)
(47,0.0638222084073477)
(48,0.05515271272320601)
(49,0.09073975267553558)
(50,0.09164496010835202)
(51,0.09181869686006355)
(52,0.09619675929365049)
};
\addplot +[lightgray] coordinates {
(1,0.45884709809407437)
(2,0.5185080200302166)
(3,0.44788279009373066)
(4,0.5162152510233252)
(5,0.5315491243697572)
(6,0.49289757405289586)
(7,0.4929846225924495)
(8,0.4867314729528604)
(9,0.4302468145307217)
(10,0.46979680787655886)
(11,0.4799619157838797)
(12,0.5122781708326976)
(13,0.4058071245093619)
(14,0.35464237899179324)
(15,0.37910292930617145)
(16,0.526846247280528)
(17,0.29118005345194126)
(18,0.19330347421773206)
(19,0.2212374272065001)
(20,0.16226810590183474)
(21,0.17394301265072673)
(22,0.15700774206468893)
(23,0.1791792962682862)
(24,0.15653230549607072)
(25,0.14370040021999997)
(26,0.15051943241947877)
(27,0.1790019492986642)
(28,0.19834184884189657)
(29,0.14946412724400332)
(30,0.14749651495908045)
(31,0.2201529156980644)
(32,0.2050744834495319)
(33,0.19285007852721686)
(34,0.18128627465876979)
(35,0.209844765069731)
(36,0.29528052368653623)
(37,0.2867570856624596)
(38,0.33754477340953554)
(39,0.43665428311039595)
(40,0.5525335151599645)
(41,0.5970761562240846)
(42,0.592134753915543)
(43,0.6821509712116353)
(44,0.5433248857768167)
(45,0.49768916801805074)
(46,0.5199865240439518)
(47,0.4755767753790457)
(48,0.4893013470330747)
(49,0.5006535786219087)
(50,0.5450113917165994)
(51,0.47399090360618323)
(52,0.5973531398339235)
};
%\legend{Biomass,RoR,PV,WindOn,Dam,Pump,BattTSO,LoadShed,LoadShift\_Down,Import}
\end{axis}
\end{tikzpicture}

%% file: week_consumer_red_wgws.tex
\pgfplotsset{compat=1.11,
/pgfplots/ybar legend/.style={
/pgfplots/legend image code/.code={%
\draw[##1,/tikz/.cd,yshift=-0.25em]
(0cm,0cm) rectangle (3pt,0.8em);},
   },}
\begin{tikzpicture}[scale=0.55]
\begin{axis}[ybar stacked,
bar width=1.9230769230769231pt,
xlabel=week of year,ylabel=Consumption (TWh),xmin=0.1,xmax=53,ymin=0,ymax=3.5,
legend style={legend columns=5,at={(1.2,1.23)},},
]

\addplot +[red] coordinates {
(1,1.5090729999999999)
(2,1.534803)
(3,1.430831)
(4,1.6742679999999999)
(5,1.671389)
(6,1.5085469999999999)
(7,1.559362)
(8,1.73267)
(9,1.458141)
(10,1.4864229999999998)
(11,1.516119)
(12,1.765277)
(13,1.6483599999999998)
(14,1.354506)
(15,1.317035)
(16,1.8004339999999999)
(17,1.621266)
(18,1.540651)
(19,1.6203459999999998)
(20,1.390447)
(21,1.4952029999999998)
(22,1.553564)
(23,1.436051)
(24,1.415387)
(25,1.78273)
(26,1.701309)
(27,1.4810269999999999)
(28,1.596339)
(29,1.453009)
(30,1.5008969999999997)
(31,1.5904969999999998)
(32,1.4461389999999998)
(33,1.4866119999999998)
(34,1.5312459999999999)
(35,1.4436799999999999)
(36,1.599218)
(37,1.589308)
(38,1.5181559999999998)
(39,1.6070559999999998)
(40,1.5392489999999999)
(41,1.7225009999999998)
(42,1.599451)
(43,1.360598)
(44,1.524659)
(45,1.7829679999999999)
(46,1.7137559999999998)
(47,1.317528)
(48,1.31647)
(49,1.556752)
(50,1.570508)
(51,1.4115529999999998)
(52,1.570552)
};
\addplot +[cyan] coordinates {
(1,6.837705348288042e-13)
(2,0.0055167594813499)
(3,0.0020700046177186144)
(4,0.0002600000007017388)
(5,0.0001653000006884139)
(6,0.009862395670447196)
(7,0.010533088621119624)
(8,0.017956082275450695)
(9,0.10246768926314738)
(10,0.04946394873763131)
(11,0.08607181405336552)
(12,0.06385887562057482)
(13,0.08347851784400627)
(14,0.06975967706153953)
(15,0.08007068190064785)
(16,0.0375711548452192)
(17,0.05928920045281188)
(18,0.08181967451569418)
(19,0.07753555851046846)
(20,0.10088344210487066)
(21,0.09141235237095023)
(22,0.0995368028786033)
(23,0.09361834634358181)
(24,0.12003988681242217)
(25,0.10239421966397479)
(26,0.08456808213734272)
(27,0.12376634952662534)
(28,0.10429986707370252)
(29,0.09813272300141561)
(30,0.16210993604833768)
(31,0.09057541875875773)
(32,0.1117585574224513)
(33,0.09953644964324854)
(34,0.07831694984123865)
(35,0.07259747802196466)
(36,0.06716287188516751)
(37,0.08685890097320224)
(38,0.07023313991519951)
(39,0.042994424162384845)
(40,0.07091943374442496)
(41,0.020934845506499444)
(42,0.022392699578224935)
(43,0.08932666720887729)
(44,0.08473070150687362)
(45,0.002910642002866077)
(46,0.0005992000006599523)
(47,0.0012017377239728214)
(48,0.030146469774033045)
(49,5.2000000680787956e-05)
(50,7.548413955650356e-05)
(51,6.965990677513191e-13)
(52,0.013585925498886806)
};
\addplot +[green] coordinates {
(1,0.0001085069445569329)
(2,0.0008336805556140877)
(3,0.0012084464592370134)
(4,0.0010852779476781184)
(5,0.00037066295944112655)
(6,0.0009166666667126702)
(7,0.0012166666667080648)
(8,0.0012166666668154327)
(9,0.0022975538010657458)
(10,0.0026258680561698296)
(11,0.0023003472230449977)
(12,0.0022885805608750262)
(13,0.002700000000924608)
(14,0.004552849207467743)
(15,0.003833283143674402)
(16,0.003275655371642874)
(17,0.002900000001217832)
(18,0.0026552950109333504)
(19,0.0018908198559807656)
(20,0.0022666666681255107)
(21,0.0022511419845562706)
(22,0.002500000001908848)
(23,0.0015962627292079137)
(24,0.002183333334958303)
(25,0.0031250237487436532)
(26,0.00234629432348872)
(27,0.0020932437808846498)
(28,0.0018201725444887947)
(29,0.002266666668486591)
(30,0.002764631689268571)
(31,0.002371615898676431)
(32,0.002661103590854234)
(33,0.0023765907226778386)
(34,0.002229934365747085)
(35,0.0012500000031443626)
(36,0.0009875000393911673)
(37,0.002874378373097899)
(38,0.0024397714599335017)
(39,0.0026246993215816895)
(40,0.002589139422094641)
(41,0.00249999999999747)
(42,0.0006085069439912138)
(43,0.0012865856321875942)
(44,0.0011967477021236858)
(45,0.00032789800026184087)
(46,0.00036985925825532874)
(47,0.0006000000000473411)
(48,0.0010333333333882946)
(49,0.0001000000000539313)
(50,0.00040000000002289886)
(51,0.0012166666666517092)
(52,0.0007166666667092618)
};
\addplot +[black] coordinates {
(1,0.09153503445134366)
(2,0.09449239378543209)
(3,0.09394632060295303)
(4,0.09594464084582854)
(5,0.09617957002503384)
(6,0.09646769471515046)
(7,0.09610017397077435)
(8,0.09629744503476154)
(9,0.07712275286217361)
(10,0.07560152379772164)
(11,0.08180477748208527)
(12,0.08366339311004253)
(13,0.08573293154178879)
(14,0.06645044108543173)
(15,0.06912308158924993)
(16,0.07752605281595891)
(17,0.07039290605836795)
(18,0.0776941590843025)
(19,0.07907425780177514)
(20,0.07368116735523669)
(21,0.07156656859836645)
(22,0.07175595668335845)
(23,0.0701870955391381)
(24,0.06331103523082449)
(25,0.07205227395410037)
(26,0.07324319624601332)
(27,0.06738276047591192)
(28,0.0707008682129552)
(29,0.06305620499191507)
(30,0.06401984912417767)
(31,0.07036891294283981)
(32,0.06464348381607171)
(33,0.06684480421592429)
(34,0.06598743118156454)
(35,0.06071555747189509)
(36,0.06563469665278837)
(37,0.06214822408956291)
(38,0.06368408292141386)
(39,0.07127810540656676)
(40,0.07253965272222027)
(41,0.08261915330832628)
(42,0.08078081848498529)
(43,0.07591287636366358)
(44,0.07883281634628073)
(45,0.09523593364628138)
(46,0.09597374137192623)
(47,0.06933997176670044)
(48,0.0583641684006865)
(49,0.09291288539608202)
(50,0.09523708615579346)
(51,0.09342507701508479)
(52,0.09653215008114477)
};
\addplot +[gray] coordinates {
(1,0.3393739779622753)
(2,0.2518503241953207)
(3,0.34074966484501235)
(4,0.2574057478749426)
(5,0.25411196718981144)
(6,0.2801409067676586)
(7,0.27282570312651744)
(8,0.28338212116916894)
(9,0.3122156006115065)
(10,0.19970994336378378)
(11,0.13903136053332243)
(12,0.21466776030118648)
(13,0.24739930451979586)
(14,0.10072784768207137)
(15,0.17318745374633945)
(16,0.1109346699954755)
(17,0.2259147102361866)
(18,0.3160840484288245)
(19,0.3438290595564655)
(20,0.4187245476506581)
(21,0.40008094713648634)
(22,0.42344678822036375)
(23,0.3282140371275553)
(24,0.3975476582885111)
(25,0.35942773199783007)
(26,0.3558326665557196)
(27,0.32337124451179194)
(28,0.3524017267863686)
(29,0.40012231879152066)
(30,0.4208945762344065)
(31,0.37393168013768413)
(32,0.421222034715448)
(33,0.41207687980209734)
(34,0.35496588507110405)
(35,0.2818032614993467)
(36,0.25932513054366446)
(37,0.2696346329273433)
(38,0.30915789692262086)
(39,0.17029853556898833)
(40,0.06337153804899803)
(41,0.01655001572879949)
(42,0.02485099352313584)
(43,0.00873854773671063)
(44,0.1921672022945593)
(45,0.30907735860692287)
(46,0.27888840703900036)
(47,0.3449430892078826)
(48,0.25263823808154356)
(49,0.30549730300503436)
(50,0.270404987321592)
(51,0.3130745982790448)
(52,0.2239054589815978)
};
\legend{Demand,Pump,BattTSO,LoadShift\_Up,Export}
\end{axis}
\end{tikzpicture}

%% file: week_consumer_red_ngws.tex
\pgfplotsset{compat=1.11,
/pgfplots/ybar legend/.style={
/pgfplots/legend image code/.code={%
\draw[##1,/tikz/.cd,yshift=-0.25em]
(0cm,0cm) rectangle (3pt,0.8em);},
   },}
\begin{tikzpicture}[scale=0.55]
\begin{axis}[ybar stacked,
bar width=1.9230769230769231pt,
xlabel=week of year,xmin=0.1,xmax=53,ymin=0,ymax=3.5,
legend style={legend columns=4,at={(1.0,-0.1)},},
]
\addplot +[red] coordinates {
(1,1.5090729999999999)
(2,1.534803)
(3,1.430831)
(4,1.6742679999999999)
(5,1.671389)
(6,1.5085469999999999)
(7,1.559362)
(8,1.73267)
(9,1.458141)
(10,1.4864229999999998)
(11,1.516119)
(12,1.765277)
(13,1.6483599999999998)
(14,1.354506)
(15,1.317035)
(16,1.8004339999999999)
(17,1.621266)
(18,1.540651)
(19,1.6203459999999998)
(20,1.390447)
(21,1.4952029999999998)
(22,1.553564)
(23,1.436051)
(24,1.415387)
(25,1.78273)
(26,1.701309)
(27,1.4810269999999999)
(28,1.596339)
(29,1.453009)
(30,1.5008969999999997)
(31,1.5904969999999998)
(32,1.4461389999999998)
(33,1.4866119999999998)
(34,1.5312459999999999)
(35,1.4436799999999999)
(36,1.599218)
(37,1.589308)
(38,1.5181559999999998)
(39,1.6070559999999998)
(40,1.5392489999999999)
(41,1.7225009999999998)
(42,1.599451)
(43,1.360598)
(44,1.524659)
(45,1.7829679999999999)
(46,1.7137559999999998)
(47,1.317528)
(48,1.31647)
(49,1.556752)
(50,1.570508)
(51,1.4115529999999998)
(52,1.570552)
};
\addplot +[cyan] coordinates {
(1,0.0011894732810362138)
(2,0.020864942415597762)
(3,0.00805425228332245)
(4,0.01009672557179471)
(5,0.012446286438537036)
(6,0.027305362451990146)
(7,0.02910321680161744)
(8,0.022057239686370272)
(9,0.1503103308591133)
(10,0.10468522578871399)
(11,0.15471435453196722)
(12,0.07775527312039504)
(13,0.16123556552538038)
(14,0.08454893889896153)
(15,0.08981708627530982)
(16,0.04478206711398815)
(17,0.07097201112769795)
(18,0.07632784037973568)
(19,0.0836552960177343)
(20,0.122670614503665)
(21,0.11758683707968738)
(22,0.1273089738041571)
(23,0.10111711458869274)
(24,0.141296904922971)
(25,0.10662407984643253)
(26,0.0915174836530167)
(27,0.11284275951591964)
(28,0.10652461403840982)
(29,0.1100828861044619)
(30,0.17805007845552212)
(31,0.10152612228959099)
(32,0.12742666547351134)
(33,0.11281867953963473)
(34,0.07561567871191205)
(35,0.06719935387447103)
(36,0.0644628005744862)
(37,0.08458742100327502)
(38,0.07788946283293352)
(39,0.047738987676839255)
(40,0.0766671171484291)
(41,0.022132399310037434)
(42,0.019939994214799653)
(43,0.06295074467191861)
(44,0.08419831995205211)
(45,7.863189803167588e-12)
(46,0.0010066489853476155)
(47,0.0056726061164262525)
(48,0.02911788435360043)
(49,0.0014704022024159549)
(50,5.275532742489213e-05)
(51,7.810417378617393e-12)
(52,0.007157939442089019)
};
\addplot +[green] coordinates {
(1,0.00040000000331927375)
(2,0.0011644095592042642)
(3,0.002847955372789399)
(4,0.0018100196372931878)
(5,0.0024517126380266033)
(6,0.002950000007988977)
(7,0.0025333333348766677)
(8,0.002483333338252551)
(9,0.005000000013730008)
(10,0.005200000011099528)
(11,0.0053255208578979795)
(12,0.00498871063595112)
(13,0.00601140116449293)
(14,0.003393517468604678)
(15,0.005251833127448431)
(16,0.003876620359736757)
(17,0.005064687082927841)
(18,0.006178175705264732)
(19,0.00469038429806331)
(20,0.006078020616272969)
(21,0.005825520887399588)
(22,0.005830411633104548)
(23,0.0038170871271836163)
(24,0.005073764051591178)
(25,0.0065348034778201)
(26,0.006433333398654245)
(27,0.005169719744624733)
(28,0.006091723284053094)
(29,0.006187483552007853)
(30,0.005721837846300905)
(31,0.005219487645416129)
(32,0.005200000056450926)
(33,0.005941840325510219)
(34,0.0059479045348829265)
(35,0.0035533879184012503)
(36,0.0033333333720341595)
(37,0.00602566155799123)
(38,0.0047763025559359905)
(39,0.004234027887538561)
(40,0.004957202410024056)
(41,0.003996765927257547)
(42,0.0037548253552933336)
(43,0.003521965758106474)
(44,0.0022582090705541525)
(45,0.0006778919965893358)
(46,0.0003092668612401951)
(47,0.001402291938296965)
(48,0.002066666668259935)
(49,2.4278634080156048e-12)
(50,0.00024566278408083727)
(51,0.0006000000018113188)
(52,0.0006618947913114048)
};
\addplot +[black] coordinates {
(1,0.08916436661822345)
(2,0.09321576400753451)
(3,0.0921523037522531)
(4,0.09539727679144842)
(5,0.09597923661793924)
(6,0.0970785502396292)
(7,0.09641608539558395)
(8,0.09628973178160727)
(9,0.07498712603860379)
(10,0.06739020372282384)
(11,0.07537246568667662)
(12,0.07653530591472331)
(13,0.08455542232788549)
(14,0.06263512385723635)
(15,0.06659410793716063)
(16,0.07107067728695045)
(17,0.0668529194016652)
(18,0.07327399933087463)
(19,0.07671762237780369)
(20,0.0720096039992735)
(21,0.07038063675646342)
(22,0.06993038829213201)
(23,0.06959529299923417)
(24,0.06192833691352994)
(25,0.07048474078227271)
(26,0.06709365063768599)
(27,0.07027855876358768)
(28,0.06895193053162191)
(29,0.06175364360122577)
(30,0.060869412247833134)
(31,0.06976388752572485)
(32,0.06318884436908828)
(33,0.06429068129664585)
(34,0.06448413738315936)
(35,0.05998789137247974)
(36,0.06424264117556076)
(37,0.06204993520857933)
(38,0.060783899045978085)
(39,0.0679540617867488)
(40,0.0675934572656449)
(41,0.0755376921481603)
(42,0.06752794943030463)
(43,0.07310237010639821)
(44,0.07673596728805737)
(45,0.090993066622283)
(46,0.09210322891741167)
(47,0.06382220840736584)
(48,0.05515271272320568)
(49,0.0907397526755386)
(50,0.09164496010835198)
(51,0.09181869686006512)
(52,0.09619675929366348)
};
\addplot +[gray] coordinates {
(1,0.321446964485517)
(2,0.22262278117575518)
(3,0.3325532018265086)
(4,0.24605913209544808)
(5,0.22894486621706434)
(6,0.22866664667888462)
(7,0.2500051162880576)
(8,0.23099418412805603)
(9,0.31639350615504885)
(10,0.1954386846463495)
(11,0.14571182196038865)
(12,0.21736454803198646)
(13,0.2855599971558663)
(14,0.27728263075349413)
(15,0.2643743373527481)
(16,0.18197336343666073)
(17,0.3749626816747325)
(18,0.43536895653649627)
(19,0.46257225854256073)
(20,0.5523660744990592)
(21,0.526150301845337)
(22,0.5592530600771927)
(23,0.49201316373357357)
(24,0.5706053287119172)
(25,0.5353260707979158)
(26,0.5337524588846039)
(27,0.5145367179890845)
(28,0.5215254097857065)
(29,0.5674621152551372)
(30,0.5664170629884949)
(31,0.5141993994710755)
(32,0.556197201764245)
(33,0.5518969255416653)
(34,0.530893896298802)
(35,0.4429038302887894)
(36,0.36486815722121135)
(37,0.3697291872550603)
(38,0.392757569247305)
(39,0.211656042165205)
(40,0.13116789525301478)
(41,0.15366056979253015)
(42,0.12500026894758684)
(43,0.025945663881654462)
(44,0.18161049214933941)
(45,0.24916153635635485)
(46,0.2521681403639)
(47,0.2709836859382453)
(48,0.20379350863322093)
(49,0.2484739019129352)
(50,0.21146413226967184)
(51,0.2817870353843948)
(52,0.14953018834312465)
};
%\legend{Demand,Pump,BattTSO,LoadShift\_Up,Export}
\end{axis}
\end{tikzpicture}

%% file: hour6840plus24_producer_red_wgws.tex
\pgfplotsset{compat=1.11,
/pgfplots/ybar legend/.style={
/pgfplots/legend image code/.code={%
\draw[##1,/tikz/.cd,yshift=-0.25em]
(0cm,0cm) rectangle (3pt,0.8em);},
   },}
\begin{tikzpicture}[scale=0.54]
\begin{axis}[ybar stacked,
bar width=4.166666666666667pt,
xlabel=hour of day, ylabel=Power (MW),xmin=0.1,xmax=25,ymin=0,ymax=30000,
legend style={legend columns=6,at={(1.45,1.3)},},
]
\addplot +[orange] coordinates {
(1,468.49999999314645)
(2,468.49999999314684)
(3,468.49999999314673)
(4,468.4999999931439)
(5,468.49999999314065)
(6,468.2999999786224)
(7,447.8999999814077)
(8,427.4999999691437)
(9,405.59999997576614)
(10,425.39999996276174)
(11,445.79999997549555)
(12,466.19999997249255)
(13,466.4000000263123)
(14,468.49999999463284)
(15,468.4999999932389)
(16,468.49999999315645)
(17,468.4999999930389)
(18,468.49999999304544)
(19,468.4999999929949)
(20,468.4999999929968)
(21,468.4999999929967)
(22,468.49999999299826)
(23,468.4999999930614)
(24,468.4999999929843)
};
\addplot +[blue] coordinates {
(1,1416.0955199539121)
(2,1416.0955199539137)
(3,1416.0955199539173)
(4,1416.0955199539108)
(5,1416.095519953902)
(6,1416.095519954493)
(7,1416.0955199542748)
(8,1416.095519956152)
(9,1248.5474399943082)
(10,1264.3511844688016)
(11,1310.3666239952317)
(12,1338.3351519961138)
(13,1336.4805439875718)
(14,1410.1094559630524)
(15,1416.0955199567622)
(16,1416.095519955855)
(17,1416.0955199529112)
(18,1416.0955199529278)
(19,1416.0955199526743)
(20,1416.0955199526768)
(21,1416.0955199526743)
(22,1416.095519952679)
(23,1416.095519953751)
(24,1416.0955199545356)
};
\addplot +[yellow] coordinates {
(1,31.37642641428955)
(2,30.235899402797532)
(3,31.033936775969806)
(4,32.74347457321786)
(5,30.22512299052022)
(6,38.390840415191086)
(7,474.0956354974217)
(8,1026.230512348605)
(9,1674.8310704817904)
(10,2617.0951978920484)
(11,2900.4132347947652)
(12,2578.72085878991)
(13,2091.399579429161)
(14,2702.8629056976642)
(15,2028.275220402601)
(16,1470.5404663750721)
(17,393.2911546711251)
(18,29.184001598612706)
(19,25.787914880766234)
(20,27.99377831808807)
(21,29.672235108212902)
(22,29.81441422008111)
(23,29.46522690571798)
(24,30.196727443043134)
};
\addplot +[lime] coordinates {
(1,22.969999978214762)
(2,29.64999997823065)
(3,28.669999978228887)
(4,42.36999997825464)
(5,27.399999978223153)
(6,11.279999990269282)
(7,15.4599999909078)
(8,18.519999992666897)
(9,18.149999999300686)
(10,22.159999999704333)
(11,13.239999999902842)
(12,15.619999999931714)
(13,18.649999999920958)
(14,47.77999999636091)
(15,46.53999999287822)
(16,32.509999992528854)
(17,17.089999949784215)
(18,18.659999949789217)
(19,20.38052590440611)
(20,19.910517379903187)
(21,23.860538919010956)
(22,23.860493470316722)
(23,11.819999987329608)
(24,18.329999991054997)
};
\addplot +[teal] coordinates {
(1,146.1099329394978)
(2,143.4950300412692)
(3,141.42280709468145)
(4,139.2621329145104)
(5,137.15058235665026)
(6,740.7000020692294)
(7,1068.8240096305392)
(8,5509.020850715021)
(9,5622.899196117014)
(10,5995.860811500031)
(11,6927.832581475783)
(12,6932.546128366787)
(13,6955.012062007529)
(14,6281.202793771389)
(15,6244.269327233975)
(16,2386.2300202280885)
(17,1.3612118308355128e-06)
(18,1.361620552944654e-06)
(19,1.528568613102446e-06)
(20,1.5285989599128504e-06)
(21,1.5285600520895064e-06)
(22,1.5288008817923942e-06)
(23,355.0000003647669)
(24,1263.6000011878525)
};
\addplot +[cyan] coordinates {
(1,3.349980171640453e-08)
(2,3.349986859423031e-08)
(3,3.3500364440828e-08)
(4,3.34962024935183e-08)
(5,3.3492938349133806e-08)
(6,3.587556978763211e-08)
(7,4.555644289461268e-08)
(8,4.402783634503202e-08)
(9,4838.107998675402)
(10,5383.397813239345)
(11,5338.798664986006)
(12,5464.683383304143)
(13,5815.4579291690525)
(14,2698.856779185179)
(15,46.78511991482809)
(16,3.829908805945104e-08)
(17,2.5676605737437848e-08)
(18,2.5682734010274748e-08)
(19,2.3842346766561324e-08)
(20,2.38420597527678e-08)
(21,2.384140589125105e-08)
(22,2.384175760668257e-08)
(23,2.802396624121002e-08)
(24,39.95459713138829)
};
\addplot +[green] coordinates {
(1,2.23649054407397e-09)
(2,2.23669174472982e-09)
(3,2.23708717635076e-09)
(4,2.23638544910798e-09)
(5,2.23633094835942e-09)
(6,2.23899093802273e-09)
(7,2.05902876331422e-09)
(8,3.22258472852269e-09)
(9,52.944904914133)
(10,81.4370219740136)
(11,82.107914915205)
(12,82.860125932594)
(13,84.6500322641066)
(14,3.90806048053527e-10)
(15,2.69005883948218e-09)
(16,3.13481860049026e-09)
(17,2.94025469699745e-09)
(18,3.47581300321776e-09)
(19,2.9625438015416e-09)
(20,3.3020360572419e-09)
(21,3.67243637613178e-09)
(22,4.10220181615748e-09)
(23,8.46921318267942e-09)
(24,9.25112440503666e-09)
};
\addplot +[brown] coordinates {
(1,1.8718167289800933e-09)
(2,1.871712733784367e-09)
(3,1.8714420036832232e-09)
(4,1.871622009638907e-09)
(5,1.8722803955992608e-09)
(6,1.846093270354183e-09)
(7,209.04000000230883)
(8,838.3000000051563)
(9,1695.640000004232)
(10,2061.0000000028044)
(11,2164.9999999998918)
(12,2164.9999999999386)
(13,2164.9999999998554)
(14,1389.1895856286103)
(15,531.8495856299745)
(16,1.848156462370171e-09)
(17,1.866958810891945e-09)
(18,1.860054556638949e-09)
(19,1.858623495939675e-09)
(20,1.8579511889070164e-09)
(21,1.858927729520695e-09)
(22,1.8579099112639258e-09)
(23,1.8472324791932203e-09)
(24,1.8242466047230242e-09)
};
\addplot +[pink] coordinates {
(1,1.05581399418401e-10)
(2,1.05586177368761e-10)
(3,1.05594062298076e-10)
(4,1.05577816195655e-10)
(5,1.05547742816064e-10)
(6,1.05561094980847e-10)
(7,1.05592713185801e-10)
(8,1.06024476127685e-10)
(9,16.8900138918603)
(10,131.6086815226)
(11,378.469730657303)
(12,475.545758738278)
(13,313.383129157479)
(14,1.06610008372982e-10)
(15,1.06070142754925e-10)
(16,1.0580220287979e-10)
(17,1.05464981660654e-10)
(18,1.0546812349639e-10)
(19,1.05527224763485e-10)
(20,1.05540834646427e-10)
(21,1.05526911808048e-10)
(22,1.05551279162099e-10)
(23,1.05580407469095e-10)
(24,1.05597037179307e-10)
};
\addplot +[black] coordinates {
(1,50.1808039093456)
(2,50.139326550386)
(3,50.1056564641635)
(4,50.0508884465907)
(5,49.9949227907761)
(6,185.602604657982)
(7,307.081860172767)
(8,782.668515894967)
(9,1203.08637603968)
(10,1080.0408725666)
(11,1688.15021724898)
(12,1829.1347803218)
(13,1826.347173567)
(14,984.476775720145)
(15,832.270355845415)
(16,766.161999552711)
(17,41.4392474524536)
(18,41.4450415600759)
(19,28.7927504161388)
(20,28.8150564063358)
(21,28.7922536498084)
(22,28.8818124314132)
(23,36.7336718118244)
(24,52.2801706723642)
};
\addplot +[lightgray] coordinates {
(1,4195.720566570583)
(2,4215.005690917455)
(3,4254.044184355946)
(4,4103.420671783398)
(5,3816.498458958141)
(6,4371.999999999031)
(7,4371.9999999989595)
(8,4371.999999999024)
(9,4731.999999999866)
(10,4731.999999999878)
(11,4556.851032004744)
(12,4521.6489530964045)
(13,4545.4175503890465)
(14,4731.999999999116)
(15,4371.999999999027)
(16,4371.999999999021)
(17,3255.999183083615)
(18,3220.185726127328)
(19,4376.820351017833)
(20,4569.679909987524)
(21,4363.735028301599)
(22,4577.495997568722)
(23,4371.999999997049)
(24,4371.999999998961)
};
\legend{Biomass,RoR,PV,WindOn,Dam,Pump,BattTSO,Gas,LoadShed,LoadShift\_Down,Import}
\end{axis}
\end{tikzpicture}

%% file: hour6840plus24_producer_red_ngws.tex
\pgfplotsset{compat=1.11,
/pgfplots/ybar legend/.style={
/pgfplots/legend image code/.code={%
\draw[##1,/tikz/.cd,yshift=-0.25em]
(0cm,0cm) rectangle (3pt,0.8em);},
   },}
\begin{tikzpicture}[scale=0.54]
\begin{axis}[ybar stacked,
bar width=4.166666666666667pt,
xlabel=hour of day,xmin=0.1,xmax=25,ymin=0,ymax=30000,
legend style={legend columns=4,at={(1.0,-0.1)},},
]
\addplot +[orange] coordinates {
(1,468.4999998697628)
(2,468.49999986973944)
(3,468.4999998697199)
(4,468.4999998696875)
(5,468.4999998696169)
(6,466.79999982918577)
(7,446.39999988695433)
(8,425.9999999242346)
(9,405.6000000347338)
(10,405.59999999869325)
(11,405.6000000265039)
(12,405.60000008021217)
(13,425.99999996244713)
(14,404.10000000984473)
(15,425.9999999076597)
(16,446.3999998800152)
(17,466.79999981928825)
(18,468.49999985636873)
(19,468.4999998543432)
(20,468.4999998543911)
(21,468.4999998543846)
(22,468.4999998544052)
(23,468.4999998621722)
(24,468.4999998622969)
};
\addplot +[blue] coordinates {
(1,1416.095519117409)
(2,1416.095519117372)
(3,1416.0955191173603)
(4,1416.095519117399)
(5,1416.0955191174373)
(6,1416.0955190356765)
(7,1416.0955190854086)
(8,1416.0955191277212)
(9,1248.5474398690578)
(10,1268.6600877372923)
(11,1310.3666238605497)
(12,1314.5890238967545)
(13,1338.3351519310231)
(14,1203.3547678638724)
(15,1416.0955191297653)
(16,1416.095519085923)
(17,1416.0955189873669)
(18,1416.0955189875028)
(19,1416.095518977845)
(20,1416.095518977863)
(21,1416.0955189778554)
(22,1416.0955189778426)
(23,1416.0955190661446)
(24,1416.095519085411)
};
\addplot +[yellow] coordinates {
(1,39.984992429915025)
(2,39.07148301406589)
(3,39.53534802903005)
(4,40.85311221956891)
(5,38.32441720187263)
(6,48.28571871413565)
(7,645.9022539970794)
(8,1497.9524906452918)
(9,2521.306562252587)
(10,3937.1379649670857)
(11,4356.491781521767)
(12,4065.0870240658264)
(13,3576.547175676206)
(14,4011.8448557419206)
(15,3584.336805341631)
(16,2341.8125262678514)
(17,603.0551163936996)
(18,37.33859483825551)
(19,32.19712168398157)
(20,34.93760915715754)
(21,36.666153767403536)
(22,36.90029072116121)
(23,36.738509672178395)
(24,37.86703437008892)
};
\addplot +[lime] coordinates {
(1,267.3100005222584)
(2,344.7300005237446)
(3,332.9500005235333)
(4,492.1100005257444)
(5,317.9100005231553)
(6,130.6799921954957)
(7,170.10000518511134)
(8,227.97000011850596)
(9,223.8499999795793)
(10,274.35999996146745)
(11,163.33999998319413)
(12,192.61999999669067)
(13,230.24999999752984)
(14,588.9699999951365)
(15,571.8799987658933)
(16,357.480005570051)
(17,188.40000078148375)
(18,205.77000078366083)
(19,236.0699990879851)
(20,230.32999908851014)
(21,276.3699990859199)
(22,276.36999908622965)
(23,130.76000248144527)
(24,201.69000517156724)
};
\addplot +[teal] coordinates {
(1,3443.45456652354)
(2,3323.9953530158477)
(3,3325.422823725916)
(4,3119.223536244449)
(5,3169.7712092762154)
(6,1.8326858133663694e-05)
(7,222.06062219690182)
(8,4258.641914544193)
(9,5653.746772245737)
(10,6167.688032222326)
(11,6921.4575175483815)
(12,6961.828633071452)
(13,7227.190591784565)
(14,6146.539009156192)
(15,3705.027456571712)
(16,1265.5703861380985)
(17,1.1730037038722437e-05)
(18,1.173027681228627e-05)
(19,1.1490530850357399e-05)
(20,1.1490493390568281e-05)
(21,1.1490494974990094e-05)
(22,1.1490497024646581e-05)
(23,2.2632272444599043e-05)
(24,370.79568133364984)
};
\addplot +[cyan] coordinates {
(1,2.867265797632083e-06)
(2,2.866498290754557e-06)
(3,2.8656754722907016e-06)
(4,2.8663161102183767e-06)
(5,2.866509181033426e-06)
(6,7.04283330966704e-07)
(7,1.0647466296900422e-06)
(8,684.2858434762212)
(9,5075.100974446787)
(10,5197.154535488505)
(11,5654.334395636701)
(12,5748.15063465119)
(13,5829.016266803194)
(14,3917.903855883077)
(15,2.110943842410939e-06)
(16,1.0647338575389699e-06)
(17,5.717292342541739e-07)
(18,5.717351801372828e-07)
(19,5.560118483411743e-07)
(20,5.559753987999098e-07)
(21,5.560266584070268e-07)
(22,5.559605783919654e-07)
(23,8.056885556325157e-07)
(24,1.0644551734413961e-06)
};
\addplot +[green] coordinates {
(1,1.109168326440332e-07)
(2,1.109225375524022e-07)
(3,1.10940835041351e-07)
(4,1.109693161277087e-07)
(5,1.110153536829681e-07)
(6,5.49779388099101e-08)
(7,8.09377244016382e-08)
(8,9.78737850209296e-08)
(9,123.4703736696022)
(10,169.0195350641615)
(11,168.5878434407623)
(12,170.340045586278)
(13,100.14111713415609)
(14,36.4410851072529)
(15,1.0728941722409541e-07)
(16,1.1772039089963479e-07)
(17,1.078854431247434e-07)
(18,1.313941581575089e-07)
(19,9.20489783618623e-08)
(20,9.67010665574335e-08)
(21,1.0049071215762699e-07)
(22,1.03902130196809e-07)
(23,1.134719890456288e-06)
(24,2.599011676464825e-06)
};
\addplot +[pink] coordinates {
(1,1.55299694061128e-09)
(2,1.55279790830148e-09)
(3,1.55263783475176e-09)
(4,1.55249407517323e-09)
(5,1.55186876872433e-09)
(6,1.55066144834976e-09)
(7,1.55124570444596e-09)
(8,1.56425217148161e-09)
(9,263.693181249528)
(10,333.59942602194)
(11,629.750124641708)
(12,718.155649835142)
(13,578.420982327365)
(14,14.5251734871357)
(15,1.56763414002315e-09)
(16,1.56015717197234e-09)
(17,1.54719065729468e-09)
(18,1.54682638587733e-09)
(19,1.54939243903446e-09)
(20,1.55022397720347e-09)
(21,1.5496782841323e-09)
(22,1.55017785466634e-09)
(23,1.55082594033999e-09)
(24,1.55101192982767e-09)
};
\addplot +[black] coordinates {
(1,119.657385325774)
(2,118.292270852104)
(3,116.695700583976)
(4,119.075176660548)
(5,117.583709157308)
(6,273.691354462853)
(7,374.889883593606)
(8,810.258213560075)
(9,1265.00216168722)
(10,1200.35252574112)
(11,1558.12374994712)
(12,1740.18324148575)
(13,1739.8350344367)
(14,1044.57806997479)
(15,909.807297932798)
(16,545.776891876331)
(17,13.826517258574)
(18,13.8560083719618)
(19,20.0052317697563)
(20,20.0463335191367)
(21,20.0041547500852)
(22,20.0876856533245)
(23,25.4457285171315)
(24,27.1659557161091)
};
\addplot +[lightgray] coordinates {
(1,2459.9999999780584)
(2,2459.999999978061)
(3,2459.99999997806)
(4,2459.999999978103)
(5,2459.9999999780957)
(6,4549.167807766604)
(7,4419.8350086640585)
(8,4653.571853709807)
(9,4731.99999999765)
(10,4731.9999999964875)
(11,4611.921458793413)
(12,4587.458925204056)
(13,4536.537018503858)
(14,4393.330129030914)
(15,4608.499321820662)
(16,4516.434434372486)
(17,3045.725745695846)
(18,3314.6204734253097)
(19,4394.655660052501)
(20,4562.641732696744)
(21,4345.126218360284)
(22,4577.548531739984)
(23,4565.351778718487)
(24,4477.383184227526)
};
%\legend{Biomass,RoR,PV,WindOn,Dam,Pump,BattTSO,LoadShed,LoadShift\_Down,Import}
\end{axis}
\end{tikzpicture}

%% file: hour6840plus24_consumer_red_wgws.tex
\pgfplotsset{compat=1.11,
/pgfplots/ybar legend/.style={
/pgfplots/legend image code/.code={%
\draw[##1,/tikz/.cd,yshift=-0.25em]
(0cm,0cm) rectangle (3pt,0.8em);},
   },}
\begin{tikzpicture}[scale=0.54]
\begin{axis}[ybar stacked,
bar width=4.166666666666667pt,
xlabel=hour of day,ylabel=Power (MW),xmin=0.1,xmax=25,ymin=0,ymax=30000,
legend style={legend columns=5,at={(1.2,1.23)},},
]
\addplot +[red] coordinates {
(1,6070.999999999999)
(2,6090.000000000002)
(3,6123.999999999999)
(4,5982.0)
(5,5671.000000000001)
(6,5927.0000000000055)
(7,6516.999999999999)
(8,13619.000000000005)
(9,21024.00000000002)
(10,23261.000000000004)
(11,25497.999999999993)
(12,25660.00000000001)
(13,25399.000000000004)
(14,20182.000000000004)
(15,14941.999999999995)
(16,10085.999999999998)
(17,4690.000000000001)
(18,4415.0)
(19,5264.0)
(20,5478.0)
(21,5257.999999999996)
(22,5652.000000000005)
(23,5929.999999999999)
(24,6069.999999999996)
};
\addplot +[cyan] coordinates {
(1,86.50000007251016)
(2,86.50000007250677)
(3,86.50000007250567)
(4,86.5000000725308)
(5,86.50000007264511)
(6,86.50000002877027)
(7,86.50000002726581)
(8,86.50000002245967)
(9,250.30000000447941)
(10,260.1007234462229)
(11,106.30000000571579)
(12,106.30000000218567)
(13,106.30000000328288)
(14,364.88263746798094)
(15,86.50000002150097)
(16,86.5000000227331)
(17,395.69137018928205)
(18,272.85345814487135)
(19,86.50000043308924)
(20,86.5000004328554)
(21,86.50000043312157)
(22,86.50000043222623)
(23,86.50000003892455)
(24,86.50000002680238)
};
\addplot +[green] coordinates {
(1,0.000144063091226799)
(2,0.000143619304487424)
(3,0.000143057829012952)
(4,0.000143718521112713)
(5,0.000144988469349683)
(6,99.8100415434566)
(7,99.9999999747655)
(8,5.66417934037188e-09)
(9,1.17416379990816e-11)
(10,1.1753280964642e-11)
(11,1.17511306758034e-11)
(12,1.17465761139533e-11)
(13,1.17524053382025e-11)
(14,4.42643310176499e-10)
(15,1.23589590911335e-08)
(16,7.92544951113875e-09)
(17,1.01222951986396e-08)
(18,6.80943910698101e-09)
(19,1.06220466184012e-08)
(20,7.91611276149857e-09)
(21,6.46923354842116e-09)
(22,5.52707642833049e-09)
(23,3.27756452006142e-09)
(24,3.1941615510925e-09)
};
\addplot +[black] coordinates {
(1,173.45310566075)
(2,176.621323182908)
(3,179.371961523173)
(4,183.942543889277)
(5,188.364461997929)
(6,897.804539315094)
(7,1246.99702525424)
(8,324.835398897285)
(9,234.397000046879)
(10,273.250859540693)
(11,202.730000052369)
(12,103.995140560348)
(13,112.898000142117)
(14,168.09565849933)
(15,604.204052442871)
(16,379.538006112211)
(17,506.723736292987)
(18,506.216832422635)
(19,985.877063274818)
(20,966.494783156671)
(21,986.15557704757)
(22,806.148238748531)
(23,666.65583690642)
(24,1147.90198918421)
};
\addplot +[gray] coordinates {
(1,0.0)
(2,0.0)
(3,0.0)
(4,0.0)
(5,0.0)
(6,221.2543862177296)
(7,359.9999999940131)
(8,359.999999997704)
(9,0.0)
(10,0.0)
(11,0.0)
(12,0.0)
(13,0.0)
(14,0.0)
(15,353.8810764944037)
(16,359.9999999981782)
(17,0.0)
(18,0.0)
(19,0.0)
(20,0.0)
(21,0.0)
(22,0.0)
(23,6.458582075797059)
(24,356.5550271341003)
};
\legend{Demand,Pump,BattTSO,LoadShift\_Up,Export}
\end{axis}
\end{tikzpicture}

%% file: hour6840plus24_consumer_red_ngws.tex
\pgfplotsset{compat=1.11,
/pgfplots/ybar legend/.style={
/pgfplots/legend image code/.code={%
\draw[##1,/tikz/.cd,yshift=-0.25em]
(0cm,0cm) rectangle (3pt,0.8em);},
   },}
\begin{tikzpicture}[scale=0.54]
\begin{axis}[ybar stacked,
bar width=4.166666666666667pt,
xlabel=hour of day,xmin=0.1,xmax=25,ymin=0,ymax=30000,
legend style={legend columns=4,at={(1.0,-0.1)},},
]
\addplot +[red] coordinates {
(1,6070.999999999999)
(2,6090.000000000002)
(3,6123.999999999999)
(4,5982.0)
(5,5671.000000000001)
(6,5927.0000000000055)
(7,6516.999999999999)
(8,13619.000000000005)
(9,21024.00000000002)
(10,23261.000000000004)
(11,25497.999999999993)
(12,25660.00000000001)
(13,25399.000000000004)
(14,20182.000000000004)
(15,14941.999999999995)
(16,10085.999999999998)
(17,4690.000000000001)
(18,4415.0)
(19,5264.0)
(20,5478.0)
(21,5257.999999999996)
(22,5652.000000000005)
(23,5929.999999999999)
(24,6069.999999999996)
};
\addplot +[cyan] coordinates {
(1,86.50000027518368)
(2,86.50000027518276)
(3,86.50000027520335)
(4,86.50000027514183)
(5,86.50000027512368)
(6,86.50000031927681)
(7,86.50000026978269)
(8,86.50000026196881)
(9,250.3000000716247)
(10,228.00000007418086)
(11,122.00000007161259)
(12,106.30000005522899)
(13,89.99952794833867)
(14,1335.0000006307769)
(15,86.50000024053062)
(16,86.50000026955308)
(17,86.50000042239132)
(18,86.50000042535139)
(19,86.50000041829061)
(20,86.50000041826111)
(21,86.50000041826183)
(22,86.50000041824367)
(23,86.50000029724505)
(24,86.50000026989544)
};
\addplot +[green] coordinates {
(1,3.06005983474936e-07)
(2,3.06026767632074e-07)
(3,3.05962125772095e-07)
(4,3.0595655490934997e-07)
(5,3.05837378630854e-07)
(6,199.9999998288958)
(7,169.9703194274545)
(8,1.219541819264339e-07)
(9,4.726943647331051e-10)
(10,3.3279519523232497e-10)
(11,3.12779428243312e-10)
(12,4.2086224093489895e-10)
(13,4.75501479308037e-10)
(14,4.73525026055961e-10)
(15,3.1729706402138097e-07)
(16,2.88141904813108e-07)
(17,4.30681233467684e-07)
(18,2.60470822664191e-07)
(19,1.453546149318551e-06)
(20,8.65299568635187e-07)
(21,6.66755405675462e-07)
(22,5.582868050465619e-07)
(23,9.47349337124247e-08)
(24,9.08015249055404e-08)
};
\addplot +[black] coordinates {
(1,65.3059538666751)
(2,65.7051538831744)
(3,65.7034668402701)
(4,66.3779199398418)
(5,66.6740986889964)
(6,671.22041094442)
(7,921.812974059339)
(8,269.275834819938)
(9,238.017465360637)
(10,196.572107123053)
(11,159.973495327158)
(12,137.713177816282)
(13,93.273810608376)
(14,244.586945618192)
(15,193.146401133432)
(16,717.069763816399)
(17,957.402910497461)
(18,954.68060800677)
(19,1217.02354169193)
(20,1168.05120415259)
(21,1218.26205585924)
(22,1057.00203720563)
(23,626.391562500666)
(24,842.997383070729)
};
\addplot +[gray] coordinates {
(1,1992.19651229511)
(2,1928.479474886411)
(3,1882.9959273832)
(4,1980.979427072562)
(5,2164.010758832102)
(6,0.0)
(7,0.0)
(8,0.0)
(9,0.0)
(10,0.0)
(11,0.0)
(12,0.0)
(13,0.0)
(14,0.0)
(15,0.0)
(16,0.0)
(17,0.0)
(18,0.0)
(19,0.0)
(20,0.0)
(21,0.0)
(22,0.0)
(23,0.0)
(24,0.0)
};
%\legend{Demand,Pump,BattTSO,LoadShift\_Up,Export}
\end{axis}
\end{tikzpicture}

%% file: bare_conf.bbl
\begin{thebibliography}{10}

\bibitem{Palmintier2016}
B.~S. Palmintier and M.~D. Webster, ``{Impact of Operational Flexibility on
  Electricity Generation Planning With Renewable and Carbon Targets},'' {\em
  IEEE Transactions on Sustainable Energy}, vol.~7, pp.~672--684, apr 2016.

\bibitem{Ma2013}
J.~Ma, V.~Silva, R.~Belhomme, D.~S. Kirschen, and L.~F. Ochoa, ``{Evaluating
  and Planning Flexibility in Sustainable Power Systems},'' {\em IEEE
  Transactions on Sustainable Energy}, vol.~4, pp.~200--209, jan 2013.

\bibitem{Maranon-Ledesma2019}
H.~Mara{\~{n}}{\'{o}}n-Ledesma and A.~Tomasgard, ``{Analyzing Demand Response
  in a Dynamic Capacity Expansion Model for the European Power Market},'' {\em
  Energies}, vol.~12, p.~2976, aug 2019.

\bibitem{Gils2016}
H.~C. Gils, ``{Economic potential for future demand response in Germany –
  Modeling approach and case study},'' {\em Applied Energy}, vol.~162,
  pp.~401--415, jan 2016.

\bibitem{GJORGIEV2022118193}
B.~Gjorgiev, J.~Garrison, {\em et~al.}, ``Nexus-e: A platform of interfaced
  high-resolution models for energy-economic assessments of future electricity
  systems,'' {\em Applied Energy}, vol.~307, p.~118193, 2022.

\bibitem{Raycheva2020}
E.~Raycheva, J.~Garrison, C.~Schaffner, and G.~Hug, ``High resolution
  generation expansion planning considering flexibility needs: the case of
  {S}witzerland in 2030,'' in {\em MEDPOWER 2020}, pp.~416--421, 2020.

\bibitem{MORALESESPANA2022122544}
G.~Morales-España, R.~Martínez-Gordón, and J.~Sijm, ``Classifying and
  modelling demand response in power systems,'' {\em Energy}, vol.~242,
  p.~122544, 2022.

\bibitem{Schwele2020}
A.~Schwele, J.~Kazempour, and P.~Pinson, ``Do unit commitment constraints
  affect generation expansion planning? {A} scalable stochastic model,'' {\em
  Energy Systems}, vol.~11, 5 2020.

\bibitem{TRONDLE20201929}
T.~Tröndle, J.~Lilliestam, S.~Marelli, and S.~Pfenninger, ``Trade-offs between
  geographic scale, cost, and infrastructure requirements for fully renewable
  electricity in europe,'' {\em Joule}, vol.~4, no.~9, pp.~1929--1948, 2020.

\bibitem{Fortenbacher2018}
P.~Fortenbacher, T.~Demiray, and C.~Schaffner, ``{Transmission network
  reduction method using nonlinear optimization},'' in {\em 2018 Power Systems
  Computation Conference (PSCC)}, pp.~1--7, IEEE, jun 2018.

\bibitem{ENTSOE2018c}
{European Network of Transmission System Operators for Electricity (ENTSO-E)},
  ``{European transmission grid dataset for the TYNDP 2018},'' 2018.

\bibitem{ENTSOE2018b}
{ENTSO-E}, ``{Cross-border physical flow},'' 2018.

\bibitem{ENTSOE-ERAA2021}
{European Network of Transmission System Operators for Electricity (ENTSO-E)},
  ``{European Resource Adequacy Assessment 2021 Edition: Net Transfer
  Capacities},'' 2021.

\bibitem{ENTSOE-MAF2020}
{European Network of Transmission System Operators for Electricity (ENTSO-E)},
  ``{Mid-term Adequacy Forecast 2020 Edition: Net Transfer Capacities},'' 2020.

\bibitem{Pickering2022}
B.~Pickering, F.~Lombardi, and S.~Pfenninger, ``Diversity of options to reach
  carbon-neutrality across the entire european energy system,'' {\em SSRN
  Electronic Journal}, 2022.

\bibitem{IEA2020b}
{International Energy Agency}, {Nuclear Energy Agency}, and {Organisation for
  Economic Cooperation and Development}, ``{Projected Costs of Generating
  Electricity 2020 Edition},'' tech. rep., Organisation for Economic
  Co-operation and Development, Paris, France, 2020.

\bibitem{ENTSO-TYNDP2020}
{European Network of Transmission System Operators for Electricity (ENTSO-E)},
  ``{Ten-Year Network Development Plan 2020},'' 2021.

\bibitem{Kiani2021}
A.~Kiani, M.~Lejeune, C.~Li, J.~Patel, and P.~Feron, ``{Liquefied synthetic
  methane from ambient CO2 and renewable H2 - A technoeconomic study},'' {\em
  Journal of Natural Gas Science and Engineering}, vol.~94, p.~104079, oct
  2021.

\bibitem{Swissgrid2025}
Swissgrid, ``{2025 transmission grid topology for Switzerland},'' 2017.

\bibitem{BFE2018a}
{Bundesamt für Energie (BFE)}, ``{Statistik der Wasserkraftanlagen der
  Schweiz},'' 2018.

\bibitem{BFE2019b}
{Bundesamt für Energie (BFE)}, ``{Kernkraftwerke Geodatenmodell
  Dokumentation},'' 2019.

\bibitem{BFE2019a}
{Bundesamt für Energie (BFE)}, ``{Windenergieanlagen Geodatenmodell
  Dokumentation},'' 2019.

\bibitem{BFE2019c}
{Bundesamt für Energie (BFE)}, ``{Kehrichtverbrennungsanlagen Geodatenmodell
  Dokumentation},'' 2019.

\bibitem{Abrell2019b}
J.~Abrell, P.~Eser, J.~B. Garrison, J.~Savelsberg, and H.~Weigt, ``{Integrating
  economic and engineering models for future electricity market evaluation: A
  Swiss case study},'' {\em Energy Strategy Reviews}, vol.~25, pp.~86--106, aug
  2019.

\bibitem{BFE_EP2050plus}
{Prognos AG}, {TEP Energy GmbH}, {Infras AG}, {Ecoplan AG}, and {Bundesamt fur
  Energie (BFE)}, ``{Energieperspektiven 2050+}.'' 2021.

\bibitem{BFE2019_Hydro}
{Bundesamt fur Energie (BFE)}, ``{Wasserkraftpotenzial der Schweiz: Abschatzung
  des Ausbaupotenxials der Wasserkraftnutzung im Rahmen der Energiestrategie
  2050},'' tech. rep., Ittigen, 2019.

\bibitem{Bauer2019}
C.~Bauer, B.~Cox, T.~Heck, and X.~Zhang, ``{Potential, costs and environmental
  assessment of electricity generation technologies: An update of potentials
  and electricity generation costs},'' tech. rep., Paul Scherrer Institut
  (PSI), 2019.

\bibitem{Bauer2017}
C.~Bauer {\em et~al.}, ``{Potentials, costs and environmental assessment of
  electricity generation technologies.},'' tech. rep., PSI, WSL, ETHZ, EPFL,
  Villigen PSI, Switzerland, 2017.

\bibitem{Kober2019}
T.~Kober, C.~Bauer, C.~Bach, M.~Beuse, G.~Georges, M.~Held, S.~Heselhaus,
  P.~Korba, L.~Kung, A.~Malhotra, S.~Moebus, D.~Parra, J.~Roth, M.~Rudisuli,
  T.~Schildhauer, T.~Schmidt, T.~Schmidt, M.~Schreiber, F.~{Segundo Sevilla},
  B.~Steffen, and S.~Teske, ``{Perspectives of Power-to-X technologies in
  Switzerland}.'' 2019.

\bibitem{CentIvDistIv_PowerTech2021}
E.~Raycheva, X.~Han, C.~Schaffner, and G.~Hug, ``Coordinated generation
  expansion planning for transmission and distribution systems,'' in {\em 2021
  IEEE Madrid PowerTech}, pp.~1--6, 2021.

\bibitem{Prognos}
A.~Kirchner {\em et~al.}, ``Energieperspektiven 2050+. {K}urzbericht.,'' tech.
  rep., Swiss Federal Office of Energy (SFOE), 2020.

\bibitem{Meteoswiss}
{Meteoswiss}, ``{Measurement Values}.''

\bibitem{PNNL2020}
K.~Mongird, ``2020 grid energy storage technology cost and performance
  assessment,'' tech. rep., Pacific Northwest National Laboratory, 2020.

\bibitem{Garrison2014AGU}
J.~Garrison, ``A grid-level unit commitment assessment of high wind penetration
  and utilization of compressed air energy storage in {ERCOT},'' 2014.

\bibitem{bynum2021pyomo}
M.~L. Bynum, G.~A. Hackebeil, W.~E. Hart, C.~D. Laird, B.~L. Nicholson, J.~D.
  Siirola, J.-P. Watson, and D.~L. Woodruff, {\em Pyomo--optimization modeling
  in python}, vol.~67.
\newblock Springer Science \& Business Media, third~ed., 2021.

\bibitem{gurobi}
{Gurobi Optimization, LLC}, ``{Gurobi Optimizer Reference Manual},'' 2022.

\bibitem{Ulbig2015}
A.~Ulbig and G.~Andersson, ``Analyzing operational flexibility of electric
  power systems,'' {\em International Journal of Electrical Power \& Energy
  Systems}, vol.~72, pp.~155--164, 11 2015.

\end{thebibliography}
